\documentclass[12pt,letter]{article}

\usepackage{amsmath} 
\usepackage{amsthm}
\usepackage{amstext}
\usepackage{amssymb}
\usepackage{enumerate}
\usepackage{array}
\usepackage{color}
\usepackage{bbm}
\usepackage{natbib}
\usepackage{graphicx}
\usepackage{caption}
\usepackage[title]{appendix}
\usepackage{url} % not crucial - just used below for the URL 
\usepackage{fancyhdr}
\usepackage{authblk}

\newcommand{\blind}{1}

\def\bs{\boldsymbol}

\newtheorem{thm}{Theorem}[section]
\newtheorem{assumption}[thm]{Assumption}
\newtheorem{rem}{Remark}
\newtheorem{lemma}{Lemma}

\def\boxit#1{\vbox{\hrule\hbox{\vrule\kern6pt
			\vbox{\kern6pt#1\kern6pt}\kern6pt\vrule}\hrule}}

\allowdisplaybreaks
\begin{document}
	\sloppy
	
	\def\spacingset#1{\renewcommand{\baselinestretch}%
		{#1}\small\normalsize} \spacingset{1.2}
	
	\if1\blind
	{
		\title{\bf Extreme Continuous Treatment Effects: Measures, Estimation and Inference}% in Continuous Treatment Effect Models} 
		%\title{\bf Nonparametric Extreme Treatment Effect} 
		\author[1]{Wei Huang \thanks{wei.huang@unimelb.edu.au}}
		\author[2]{Shuo Li \thanks{shuoli@tjufe.edu.cn}}
		\author[1]{Liuhua Peng \thanks{liuhua.peng@unimelb.edu.au}}
		\affil[1]{School of Mathematics and Statistics, University of Melbourne, Australia}
		\affil[2]{School of Statistics, Tianjin University of Finance and Economics, Tianjin, China}
		\date{}
		\maketitle
	} \fi
	
	\if0\blind
	{
		\bigskip
		\bigskip
		\bigskip
		\begin{center}
			{\LARGE \bf Nonparametric Extreme Treatment Effect}
		\end{center}
		\begin{center}
			Author 1 and Author 2
		\end{center}
		\medskip
	} \fi
	
	%\bigskip
	\begin{abstract}
		This paper concerns estimation and inference for treatment effects in deep tails of the counterfactual distribution of unobservable potential outcomes corresponding to a continuously valued treatment. We consider two measures for the deep tail characteristics: the extreme quantile function and the tail mean function defined as the conditional mean beyond a quantile level. Then we define the extreme quantile treatment effect (EQTE) and the extreme average treatment effect (EATE), which can be identified through the commonly adopted unconfoundedness condition and estimated with the aid of extreme value theory. Our limiting theory is for the EQTE and EATE processes indexed by a set of quantile levels and hence facilitates uniform inference. 
		%We propose a novel bootstrap procedure for finite-sample implementation. 
		%Our theory is focused on the continuous treatment effect model, but can be readily extended to the discrete model which is technically easier to handle. 
		Simulations suggest that our method works well in finite samples and an empirical application illustrates its practical merit.
	\end{abstract}
	
	\vfill
	
	\newpage
	\spacingset{1.5} % DON'T change the spacing!
	\section{Introduction}
	This paper is devoted to measuring, estimation, and inference for treatment effects on deep tails of the potential outcome distributions, namely the extreme treatment effect (ETE). For a tail level $\alpha$, we consider two measures for tail characteristics: the $\alpha$th-quantile and the $\alpha$th-tail mean defined as the conditional mean beyond the $\alpha$th-quantile. Then for $\alpha$ that is close to 0 or 1, we define the extreme quantile treatment effect (EQTE) and the extreme average treatment effect (EATE), which are the differences of their corresponding measures at different treatment status.  Traditionally, analysis of treatment effects is based on measures from the average perspective such as the average treatment effect (ATE); see, for example, \citet[][]{Abadie2018ARE} for a comprehensive review. Due to its ability to capture heterogeneity, the quantile treatment effect (QTE) has been enjoying much popularity since the pioneering work of \citet{Firpo2007ECA}. The EQTE and EATE considered in this paper are, respectively, conceptual extensions of the traditional QTE and ATE to extreme situations.
	
	In this paper, we employ the commonly adopted unconfoundedness condition to identify the counterfactual potential outcome distributions. Once the potential outcome distributions are identified, so are the  the EQTE and EATE. Then, assuming that the potential outcome distributions have Pareto-type tails, we construct estimators of the EQTE and EATE based on tail approximations from extreme value theory (EVT). When deriving the limiting results, we work under the double asymptotics, namely $N \to \infty$ and $\alpha=\alpha_N \to 1$, where $N$ is the sample size and $\alpha$ is the quantile level. In addition, our limiting results are for EQTE and EATE process indexed by a set of $\alpha$.%, thereby facilitating uniform inference over multiple quantile tail levels. 

 	This paper is related to some other papers that consider estimation and inference for extremal conditional quantiles. \citet{Chernozhukov2005AOS} and \citet{Chernozhukov2011RES} use EVT to estimate extremal conditional quantiles in a linear model. Our context is different form theirs because we consider the unconditional quantiles and do not impose any parametric structure on the quantile function. In addition, \citet{Li2022ET} propose EVT-based estimators for the conditional quantile and conditional tail mean (that is, conditional value-at-risk and conditional expected shortfall) in a semiparametric dynamic model. We note that the counterfactual nature of the treatment effect model distinguishes the current paper from \citet{Chernozhukov2005AOS}, \citet{Chernozhukov2011RES}, and \citet{Li2022ET}. Finally, this paper is related to the recent paper of \citet{Zhang2018AOS}, which builds on  \citet{Chernozhukov2005AOS} and \citet{Chernozhukov2011RES} and consider the EQTE in the binary treatment effect model. A main methodological difference between our method and \citet{Zhang2018AOS} lies in that we adopt EVT in the estimation but \citet{Zhang2018AOS} estimate extreme quantiles using the conventional quantile regression approach.

	The main contributions of this paper are summarized as follows. First, apart from the conventional quantile function, we propose to add the tail mean function into the conceptual framework, thereby enriching the toolkit for treatment effect analysis and providing empirical researchers more choices. Unlike the quantile function which can only tell us the extreme situation, the tail mean function goes further to tell us what the situation would be once things are more extreme than the quantile. In addition, both the EQTE and the EATE can capture heterogeneity. Second, to the best of our knowledge, we are the first to consider ETE in the continuous treatment effect model. 
	The estimation of the counterfactual distribution is adopted from the recent developments in the continuous treatment effect model \citep[e.g.,][]{ai2021estimation}.
	However, our theoretical investigation is distinct from \citep[e.g.,][]{ai2021estimation} as they consider a fixed tail level $\alpha$. 
	%Though focusing on the case of continuous treatments, our theory can be readily extended to the discrete treatment effect model which is technically easier to handle. 
	Third, a theoretical highlight of this paper is that our limiting results are for the EQTE and EATE processes and hence facilitates uniform inference over multiple  tail levels. Compared with inference at a single tail level, uniform inference can provide an overall picture of the tail feature. 
	%Fourth, we propose a novel bootstrap procedure for finite-sample implementation. The bootstrap procedure is based on uniform expansions of the EQTE and EATE estimators and employs the idea of multiplier bootstrap. 
	%Simulations and a real data example confirm that the bootstrap procedure delivers favorable finite-sample performance.
	The proposed approaches have been justified by both simulations and a real data example.
	
	The rest of this paper is organized as follows. Section \ref{sec:framework} introduces the framework. Section \ref{sec:F} describes the estimation procedure of the survival function along with its asymptotic 
	property. 
	Section \ref{sec:quantile} studies both intermediate and extreme quantile, as well as the tail mean.
	The estimation and inference of the EQTE and EATE is established in Section \ref{sec:TE}.
	%The main limiting theory is given in Section \ref{sec:asymptotic}. %Section \ref{sec:bootstrap} introduces the bootstrap procedure. %Section \ref{sec:test} discusses testing for treatment effects. 
	Simulations and real data applications are contained in Section \ref{sec:Numerical_study}.

	\section{Basic Framework}\label{sec:framework}
	We consider the continuously valued treatment situation where the observed treatment variable is denoted by $T$ with probability density function $f_T(t)$ and support $\mathcal{T}\subset \mathbb{R}$. Let $Y(t)$ denotes the potential outcome if one was treated at level $t$ for $t\in\mathcal{T}$. In practice, each individual can only receive one treatment level $T$ and we only observe the corresponding outcome $Y:=Y(T)$. We are also given a vector of covariates $\bs{X}\in\mathbb{R}^r$, with $r$ a positive integer, that is related to both $T$ and $Y(t)$ for $t\in\mathcal{T}$.
	
	The goal of this paper is to estimate and infer the EQTE and EATE based on the extreme quantile function of $Y(t)$, for $t\in\mathcal{T}$, of level $\alpha\in(0,1)$, from a random sample $\{T_i,\boldsymbol{X}_i,Y_i\}_{i=1}^N$, where we allow $\alpha=\alpha_N\rightarrow 1$ as $N\rightarrow \infty$.
	
	Specifically, we call $\bar{F}_t(y):=\mathbb{P}\{Y(t)\geq y\}$ the counterfactual survival function of $Y(t)$. 
	The $\alpha$-quantile function and the $\alpha$-level tail mean of $Y(t)$ are defined by
	\begin{eqnarray}\label{Def:qtCTM}
	    q_t(\alpha) := \inf\{y: \bar{F}_t(y)\leq 1-\alpha\} \quad \text{and} \quad \mathrm{TM}_{t}(\alpha) := \mathbb{E}\left\{Y(t)\mid Y(t)>q_{t}(\alpha)\right\}\,,
	\end{eqnarray}
	respectively.
	Then, with $\alpha=\alpha_N\to1$ as $N\to\infty$, we are ready to define the EQTE and EATE for any two fixed treatments $t_1\neq t_2 \in \mathcal{T}$, which are 
	\begin{eqnarray}
	    \mathrm{EQTE}_{t_1,t_2}(\alpha_N) = \frac{q_{t_1}(\alpha_N)}{q_{t_2}(\alpha_N)}\,,
	\end{eqnarray}
	and
	\begin{eqnarray}
	    \mathrm{EATE}_{t_1,t_2}(\alpha_N) = \frac{\mathrm{TM}_{t_1}(\alpha_N)}{\mathrm{TM}_{t_2}(\alpha_N)}\,,
	\end{eqnarray}
	respectively. It is worth mentioning that when $\alpha_N\to1$, $q_{t}(\alpha_N)\to\infty$ for any fixed $t\in\mathcal{T}$. Thus, we consider comparing the quantile and tail mean in the relative way.
	
	In particular, we are interested in the following heavy tail distribution of $Y(t)$, which is commonly employed in the extreme value theory literature (see \citealp[Theorem 1.1.8]{deHaan2006}): 
	\begin{assumption}\label{assump:extreme_domain}
		The function $\bar{F}_{t}(y)$ is continuously differentiable in $y$ and third-order continuously differentiable with respect to $t$.
		In addition, there exists a function $\gamma(t)>0$ such that for any $t\in\mathcal{T}$,
		\begin{equation}
			\lim_{y\to\infty}\frac{yf_{Y(t)}(y)}{\bar{F}_{t}(y)} = \frac{1}{\gamma(t)},
		\end{equation}
		where $f_{Y(t)}(\cdot)$ is the density function of $Y(t)$.
	\end{assumption} 
	Assumption \ref{assump:extreme_domain} implies that $F_{t}(y)$ is in the domain of attraction of the Fr\'{e}chet distribution \citep[Theorem 1.1.11]{deHaan2006}, or equivalently,
	\begin{equation}\label{eq:tail_index_F}
		\lim_{y\to\infty}\frac{\bar{F}_{t}(zy)}{\bar{F}_{t}(y)}=z^{-1/\gamma(t)}, \text{~~for any~} z>0,
	\end{equation}
	which means that $\bar{F}_{t}(y)$ is regularly varying with index $-1/\gamma(t)$. $\gamma(\cdot)$ is an unknown function of $t$ referred to as the extreme-value index.
	
	Note that $Y(t)$ is never observed simultaneously for all $t\in\mathcal{T}$ or even on a dense subset of $\mathcal{T}$ for any individual, but only at a particular $T$ level of treatment. Thus, in order to identify $\bar{F}_t$ and $\gamma$ from the observed data, we adopt the following \emph{unconfoundedness} assumption, which is imposed in most of the treatment effect literature (see e.g. \citealp{rosenbaum1983central}, \citealp{Hirano03}, \citealp{chernozhukov2013inference}, \citealp{donald2014testing} and \citealp{Ai_Linton_Motegi_Zhang_cts_treat}):
	\begin{assumption}[Unconfoundedness]\label{UnconfoundAssump}
		For all $t\in\mathcal{T}$, given $\boldsymbol{X}$, $T$ is independent of $Y(t)$, i.e. $T \perp \{Y(t)\}_{t\in\mathcal{T}} | \bs{X}$.
	\end{assumption}
	Under Assumption \ref{UnconfoundAssump},  for any suitable function $\phi(\cdot)$ and every fixed $t\in\mathcal{T}$, we have
	\begin{align}
		\mathbb{E}[\phi\{Y(t)\}]
		= &\mathbb{E}(\mathbb{E}[\phi\{Y(t)\}|\bs{X}])
		= \mathbb{E}(\mathbb{E}[\phi\{Y(t)\}|\bs{X}, T=t])
		= \mathbb{E}[\mathbb{E}\{\phi(Y)|\bs{X},T=t\}]\notag \\
		=& \int_{\mathcal{X}}\int_{\mathcal{Y}}\phi(y)f_{Y|X,T}(y|\bs{x},t)f_X(\bs{x})\,dy\,d\bs{x}\nonumber\\
		%=& \int_{\mathcal{X}}\int_{\mathcal{Y}}\frac{f_{T}(t)}{f_{T|X}(t|\bs{x})}\phi(y)f_{Y|X,T}(y|\bs{x},t)f_{X|T}(\bs{x}|t)\,dy\,d\bs{x}\nonumber\\
		=& \int_{\mathcal{X}}\int_{\mathcal{Y}}\frac{f_{T}(t)}{f_{T|X}(t|\bs{x})}\phi(y)f_{Y,X|T}(y,\bs{x}|t)\,dy\,d\bs{x}\nonumber\\
		=& \mathbb{E}\{\pi_0(T,\bs{X})\phi(Y)|T=t\}\,,\label{PC}
	\end{align}
	where 
	\begin{equation}
		\pi_0(t,\bs{x}) := \frac{f_T(t)}{f_{T|X}(t|\bs{x})}\,,\label{pi0def}
	\end{equation}
	with $f_T$ and $f_{T|X}$ the density function of $T$ and the conditional density of $T$ given $\bs{X}$, respectively. The function $\pi_0(t,\boldsymbol{x})$ is called the \emph{stabilized weights} (see e.g.~\citealp{Ai_Linton_Motegi_Zhang_cts_treat}). 
	
    Thus, the survival function $\bar{F}_{t}(y)$ can be identified from $(T,\boldsymbol{X}, Y)$ by 
    \begin{equation}\label{Identify:Fbar}
        \bar{F}_{t}(y) = \mathbb{E}\{\pi_0(T,\boldsymbol{X})\mathbbm{1}(Y>y)|T=t\}\,.
    \end{equation}
    In the following section, we introduce a nonparametric estimator of $\bar{F}_t$ and its asymptotic behaviour corresponding to the extreme quantiles, where $\alpha=\alpha_N\rightarrow 1$ as $N\rightarrow \infty$. From there, we can form estimators for the index function $\gamma$, the extreme quantile function $q_t(\alpha)$ and the test of EQTE and EATE.

	\section{Estimator of the Survival Function and the Large Sample Properties}\label{sec:F}
	
	Suppose we have a consistent estimator of the stabilized weight $\pi_0$, denoted by $\hat{\pi}$ (examples of such an estimator can be found in the treatment effect literature, e.g.~\citealp{ai2021estimation}). We can apply nonparametric regression method to \eqref{Identify:Fbar} with $\pi_0$ replaced by $\hat{\pi}$ to estimate $\bar{F}_t$. 
	It is noticed that $\mathbb{E}\{\pi_0(T, \bs{X})|T=t\}=1$ for any $T\in\mathcal{T}$, then to ensure our estimated survival function lies in the unit interval $[0,1]$, we normalize $\bar{F}_{t}(y)$ as
	\begin{equation}\label{eq:normalization}
		\bar{F}_{t}(y) = \frac{\mathbb{E}\{\pi_0(T,\bs{X})\mathbbm{1}(Y> y)|T=t\}}{\mathbb{E}\{\pi_0(T, \bs{X})|T=t\}}.
	\end{equation}
	
	We then estimate $\bar{F}_{t}(y)$ by plugging $\hat{\pi}$ into \eqref{eq:normalization} and employing kernel methods on both its numerator and denominator. Specifically, let $K(\cdot)$ be a prespecified kernel function on $\mathbb{R}$ and $K_h(t) = K(t/h)$ for $t,h\in \mathbb{R}$, we define
	\begin{eqnarray}
		\notag \widehat{\bar{F}}_{t,h}(y) = \frac{\sum_{i=1}^{N}\hat{\pi}(T_i,\mathbf{X}_i)\mathbbm{1}(Y_i>y)K_h(T_i-t)}{\sum_{i=1}^{N}\hat{\pi}(T_i,\mathbf{X}_i)K_h(T_i-t)}
	\end{eqnarray}
	as the estimator of $\bar{F}_{t}(y)$, where $h$ is the bandwidth. 
	
	In order to estimate and infer the quantiles of $F_{t}(y)$ in the tail region, where $\alpha=\alpha_N\to1$ as $N\to\infty$, we need to study the asymptotics of $\widehat{F}_{t}(y)$ for $y=y_N\to y_{F}(t)$ as $N\to\infty$, where $y_{F}(t)=\sup\{y: F_{t}(y)<1\}$ is the right endpoint of $F_{t}(y)$. 
	Under Assumption \ref{assump:extreme_domain}, $y_{F}(t)=\infty$.
	We make the following assumptions.
	
	\begin{assumption}
		\label{as:suppX} (i) The support $\mathcal{X}$ of the control variables  $%
		\boldsymbol{X}$ is a compact subset of $\mathbb{R}^{r}$.  
		(ii)The support  $\mathcal{T}$ of the treatment variable $T$ is a compact subset of $\mathbb{R}$. Without loss of generality, we assume $\mathcal{T}=[0,1]$.
		(iii) There exist two positive constants $\eta_{1}$ and $\eta_{2}$ such 
		that  
		\begin{equation*}
			0<\eta_{1}\leq\pi_{0}(t,\boldsymbol{x})\leq\eta_{2}<\infty\ ,\ \forall ~  (t,%
			\boldsymbol{x})\in\mathcal{T}\times\mathcal{X}\ .  
		\end{equation*}
	\end{assumption}
	
	\begin{assumption}\label{assump:f_T}
		(i) For any given $(t,\mathbf{x})\in\mathcal{T}\times\mathcal{X}$, the conditional distribution function $F_{Y|T,\mathbf{X}}(y|t,\mathbf{x})$ is continuous in $y$. (ii) For any $y$, $F_{Y|T,\mathbf{X}}(y|t,\mathbf{x})$ is continuously differentiable with respect to $(t,\mathbf{x})\in\mathcal{T}\times\mathcal{X}$. (iii) The density function $f_{T}(t)$ is third-order continuously differentiable. In addition, there exist two positive constants $\eta_{3}$ and $\eta_{4}$ such that
		\begin{equation}
			\notag 0<\eta_{3}\leq f_{T}(t) \leq \eta_{4}<\infty, ~ \forall ~t\in\mathcal{T}.
		\end{equation}
	\end{assumption}
	
	\begin{assumption}\label{asaump:K} 
		$K(\cdot)$ is a bounded univariate kernel function, symmetric around
		the  origin, with support $\mathcal{K}$ included in the unit interval $[-1,1]$, and satisfying (i) $\int_{\mathcal{K}} K(u)du=1$; (ii) $\int_{\mathcal{K}}
		u^{2}K(u)du=\kappa  _{21}\in(0,\infty)$; (iii) $\int_{\mathcal{K}}
		K^{2}(u)du=\kappa_{02}<\infty$.
	\end{assumption}
	
	Note that many consistent estimator of $\pi_0$ are available in the literature. We here assume the following high-level condition on $\hat{\pi}$:
	\begin{assumption}\label{assump:pi}
		The estimator $\hat{\pi}(\cdot,\cdot)$ of $\pi_{0}(\cdot,\cdot)$ satisfies
		\begin{equation}
			\sup_{1\leq i\leq N}\left|\hat{\pi}(T_i,\mathbf{X}_i)-\pi_{0}(T_i,\mathbf{X}_i)\right| = O_{p}(\delta_N),
		\end{equation}
		where $\delta_N\to 0$ at a certain rate.
	\end{assumption}
	For example, $\hat{\pi}$ can be a nonparametric general empirical likelihood estimator satisfying an expanding set of moment equations introduced by \citet{ai2021estimation}. In that case, $\delta_N = L_N^{-s+1/2}+N^{-1/2}L_N$, where $L_N=o(N^{1/2})\rightarrow\infty$, with $s=+\infty$ if both $T$ and $\boldsymbol{X}$ are discrete, and if some components of $(T,\boldsymbol{X})$ are continuous, then $s>1/2$ relates to the smoothness of $\pi_0(t,\boldsymbol{x})$ in terms of the continuous components.
	
	Denote $B(t, r)=\{t_0: |t_0-t|\leq r\}$ as the ball with center $t$ and radius $r$. The following quantity is used to measure the discrepancy between the extremes of the counterfactual survival distributions of $Y$ at neighboring points of $t$. Define the oscillation of $\bar{F}_{t}(y)$ above a high level $y_N$ as
	\begin{equation}\label{eq:oscillation}
		\omega_t(y_N, h) = \sup_{y\geq y_N, t^{\prime}\in B(t, h)}\frac{1}{\log y}\left|\log\left\{\frac{\bar{F}_{t^{\prime}}(y)}{\bar{F}_{t}(y)}\right\}\right|.
	\end{equation}
	
	\begin{rem}
		In order to have a clear insight into the order of $\omega_{t}(y_N,h)$, consider that $F_{t}(y)$ belongs to the generalized Hall class of heavy-tailed distribution \citep{Hall1982}, that is, 
		\begin{eqnarray}\label{eq:Hall_class}
			\bar{F}_{t}(y)= 1-F_{t}(y) = y^{-1/\gamma(t)}\left\{c_0(t)+c_1(t)y^{-\beta(t)}+o(y^{-\beta(t)})\right\},
		\end{eqnarray}
		where $\gamma(t)$, $c_0(t)$, $c_1(t)$ and $\beta(t)$ are functions of $t$ with $\gamma(t)>0$, $c_0(t)>0$ and $\beta(t)>0$. It is clear that the $\bar{F}_{t}(y)$ defined in \eqref{eq:Hall_class} satisfies Assumption \ref{assump:extreme_domain} with index function $\gamma(\cdot)$.  Then
		\begin{eqnarray}
			\notag & & \omega_t(y_N, h) \\%& = & \sup_{y\geq y_N, t^{\prime}\in B(t, h)}\frac{1}{\log y}\left|\log\left\{\frac{\bar{F}_{t^{\prime}}(y)}{\bar{F}_{t}(y)}\right\}\right| \\
			\notag & = & \sup_{y\geq y_N, t^{\prime}\in B(t, h)}\frac{1}{\log y}\left|\left\{\frac{1}{\gamma(t)}-\frac{1}{\gamma(t^{\prime})}\right\}\log(y)+\log\left\{\frac{c_0(t^{\prime})+c_1(t^{\prime})y^{-\beta(t^{\prime})}+o(y^{-\beta(t^{\prime})})}{c_0(t)+c_1(t)y^{-\beta(t)}+o(y^{-\beta(t)})}\right\}\right|.
		\end{eqnarray}
		Thus, $\omega_t(y_N, h)=O(h)$ under suitable (Lipschitz) conditions on $\gamma(t)$, $c_0(t)$, $c_1(t)$ and $\beta(t)$.
	\end{rem}
	
	We have the following conditions on the rates of $h$ and $y_N$.
	\begin{assumption}\label{asaump:y_N} 
		Let $y_N\to \infty$, $h\to 0$ as $N\to\infty$. In addition, for any fixed $t\in\mathcal{T}$, $\omega_t(y_{N},h)\log(y_N)\to0$, $Nh\bar{F}_{t}(y_N)\to\infty$, $Nh^3\bar{F}_{t}(y_N)\to0$, $\sqrt{Nh\bar{F}_{t}(y_N)}\omega_{t}(y_N,h)\log(y_N)\to0$ and $Nh\bar{F}_{t}(y_N)\delta_N^2\to 0$. 
	\end{assumption}
	
	Let $\Delta$ be any fixed constant that is greater or equal to $1$ throughout the paper. Let $\Rightarrow$ denote week convergence of stochastic processes.
	
	\begin{lemma}\label{lemma:01}
		Suppose Assumptions \ref{assump:extreme_domain}, \ref{UnconfoundAssump}, \ref{as:suppX}, \ref{assump:f_T}, \ref{asaump:K}, \ref{assump:pi} and \ref{asaump:y_N} hold. 
		Denote for $\upsilon\in[1, \Delta]$,
		\begin{equation}
			\mathcal{E}_{t,h}(\upsilon) = \sqrt{Nh\bar{F}_{t}(\upsilon y_N)}\left\{\frac{\widehat{\bar{F}}_{t,h}(\upsilon y_N)}{\bar{F}_{t}(\upsilon y_N)}-1\right\}.
		\end{equation}
		If there exists a function $\varpi_{t}^{F}(\upsilon_1,\upsilon_2)$ such that
		\begin{eqnarray}
		    \notag \varpi_{t}^{F}(\upsilon_1,\upsilon_2)=\lim_{y_N\to\infty}\frac{ \mathbb{E}\left[\left\{\pi_{0}(T,\mathbf{X})\right\}^2\mathbbm{1}(Y>\upsilon_1 y_N)\mathbbm{1}(Y>\upsilon_2 y_N)|T=t\right]}{\sqrt{\bar{F}_{t}(\upsilon_1y_N)\bar{F}_{t}(\upsilon_2y_N)}},
		\end{eqnarray}
		let $\Omega_{t}^{F}(\upsilon_1,\upsilon_2)=\kappa_{02}\{f_{T}(t)\}^{-1}\varpi_{t}^{F}(\upsilon_1,\upsilon_2)$.
		Then, for any fixed $t_1\neq t_2\in\mathcal{T}$,
		\begin{eqnarray}
			\notag {\mathcal{E}_{t_1,h}(\upsilon) \choose \mathcal{E}_{t_2,h}(\upsilon)} \Rightarrow \Psi_{t_1,t_2}^{F}(\upsilon)
		\end{eqnarray}
		in $\ell^{\infty}([1,\Delta])$ as $N\to\infty$, where $\Psi_{t_1,t_2}^{F}(\upsilon)$ is a centered Gaussian process with covariance function
		\begin{eqnarray}
			\Omega_{t_1, t_2}^{F}(\upsilon_{1}, \upsilon_{2}) = \left( 
			\begin{array}{cc}
				\Omega_{t_1}^{F}(\upsilon_1,\upsilon_2) & 0 \\
				0 & \Omega_{t_2}^{F}(\upsilon_1,\upsilon_2)
			\end{array}
			\right).
		\end{eqnarray}
	\end{lemma}
	
	Lemma \ref{lemma:01} establishes joint limiting behavior of $\mathcal{E}_{t_1,h}(\upsilon)$ and $\mathcal{E}_{t_2,h}(\upsilon)$ for different $t_1$ and $t_2$, and it implies that $\mathcal{E}_{t_1,h}(\upsilon)$ and $\mathcal{E}_{t_2,h}(\upsilon)$ are asymptotically independent for any $t_1\neq t_2\in\mathcal{T}$.
	Moreover, the convergence is established uniformly over a tail region, which facilitates uniform inference.
	Specifically, when $\upsilon=1$,
	\begin{eqnarray}
	    \notag \mathcal{E}_{t,h}(1) = \sqrt{Nh\bar{F}_{t}(y_N)}\left\{\frac{\widehat{\bar{F}}_{t,h}(y_N)}{\bar{F}_{t}(y_N)}-1\right\} \to N(0, \Omega_{t}^{F}(1,1))
	\end{eqnarray}
	in distribution, where $\Omega_{t}^{F}(1,1)=\lim_{y_N\to\infty}\frac{ \mathbb{E}\left[\left\{\pi_{0}(T,\mathbf{X})\right\}^2\mathbbm{1}(Y>y_N)|T=t\right]}{\bar{F}_{t}(y_N)}$.
	
	It follows from the proof of Lemma \ref{lemma:01} that uniformly in $\upsilon\in[1,\Delta]$,
	\begin{eqnarray}
		\notag \notag \mathcal{E}_{t,h}(\upsilon) = \sum_{i=1}^{N}\psi_{t,h}(T_i, \mathbf{X}_i, Y_i, \upsilon, y_N) + o_{p}(1),
	\end{eqnarray}
	where
	\begin{eqnarray}
		\psi_{t,h}(T_i, \mathbf{X}_i, Y_i, \upsilon, y_N) =  \sqrt{\frac{h\bar{F}_{t}(\upsilon y_N)}{N}}\left\{\frac{\pi_0(T_i,\mathbf{X}_i)\mathbbm{1}(Y_i>\upsilon y_N)K_h\left(T_i-t\right)}{\bar{F}_{t}(\upsilon y_N)f_{T}(t)}-1\right\}
	\end{eqnarray}
	for $i=1,\ldots,N$.
	Let $\hat{f}_{T,h}(t) = N^{-1}\sum_{i=1}^{N}K_h\left(T_i-t\right)$ be the kernel estimator of $f_{T}(t)$. 
	For any $t\in\mathcal{T}$ and $\upsilon_1,\upsilon_2\in[1,\Delta]$, we propose estimating the covariance function $\Omega_{t}^{F}(\upsilon_1,\upsilon_2)$ using
	\begin{eqnarray}\label{eq:Omega_est}
		\hat{\Omega}_{t}^{F}(\upsilon_1,\upsilon_2) = \sum_{i=1}^{N} \hat{\psi}_{t,h}(T_i, \mathbf{X}_i, Y_i, \upsilon_1, y_N) \hat{\psi}_{t,h}(T_i, \mathbf{X}_i, Y_i, \upsilon_2, y_N),
	\end{eqnarray}
	where
	\begin{eqnarray}
		\notag  \hat{\psi}_{t,h}(T_i, \mathbf{X}_i, Y_i, \upsilon, y_N) = \sqrt{\frac{h\widehat{\bar{F}}_{t,h}(\upsilon y_N)}{N}}\left\{\frac{\hat{\pi}(T_i,\mathbf{X}_i)\mathbbm{1}(Y_i>\upsilon y_N)K_h\left(T_i-t\right)}{\widehat{\bar{F}}_{t,h}(\upsilon y_N)\hat{f}_{T,h}(t)}-1\right\}.
	\end{eqnarray}
	
	The following theorem shows that $\hat{\Omega}_{t}^{F}(\upsilon_1,\upsilon_2)$ is consistent for $\Omega_{t}^{F}(\upsilon_1,\upsilon_2)$.
	
	\begin{thm}\label{theo:Omega_est}
		Under the conditions assumed in Lemma \ref{lemma:01}, for any fixed $t\in\mathcal{T}$ and $\upsilon_1,\upsilon_2\in[1,\Delta]$, $\hat{\Omega}_{t}^{F}(\upsilon_1,\upsilon_2)/\Omega_{t}^{F}(\upsilon_1,\upsilon_2) \to 1$
		in probability as $N\to\infty$.
	\end{thm}
	
	\section{Estimate and Infer the Quantile Function and the Tail Mean Function}\label{sec:quantile}
	Depending on the rate of the quantile level $\alpha=\alpha_N$, compared to $N$, the estimation of the quantile function $q_t(\alpha_N)$ differs. Specifically, the case where $Nh(1-\alpha_N)\to\infty$ is called the intermediate case and the one where $Nh(1-\alpha_N)\to c\in[0,\infty)$ is called the extreme case. We will see in the following subsection that, following from the convergence rate of our survival function estimator in Lemma~\ref{lemma:01}, the intermediate quantile can be consistently estimated from $\widehat{\bar{F}}_{t,h}$ using the definition in \eqref{Def:qtCTM}. However, the extreme quantile estimation require estimation of the index function $\gamma$ and making use of~\eqref{eq:tail_index_F}.
	
	\subsection{Intermediate Quantile}\label{sec:intermediate}
	
	In this section, we consider estimating the $\alpha$-quantile of $Y(t)$, that is $q_{t}(\alpha)$, in the intermediate case when $\alpha=\alpha_N\to1$ and $Nh(1-\alpha_N)\to\infty$. Recall $q_{t}(\alpha)=\inf\left\{z: \bar{F}_{t}(z)\leq1-\alpha\right\}$, we considering using the $\alpha_N$-quantile of $1-\widehat{\bar{F}}_{t,h}(y)$ to estimate $q_{t}(\alpha_N)$, that is, define
	\begin{equation}\label{eq:intermediate_quantile_est}
		\hat{q}_{t, h}(\alpha_N) = \inf\{z: 1-\widehat{\bar{F}}_{t,h}(z)\geq\alpha_N\} = \inf\{z: \widehat{\bar{F}}_{t,h}(z)\leq1-\alpha_N\}.
	\end{equation}
	
	We need to following second-order condition on the survival function $\bar{F}_{t}(y)$.
	
	\begin{assumption}\label{assump:second_order}
		There exist functions $\rho(t)< 0$  and $A_t(\cdot)$ which is eventually positive or negative with $\lim_{y\to\infty}A_t(y)=0$, such that
		\begin{equation}\label{eq:second_order}
			\lim_{y\to\infty}\frac{\bar{F}_{t}(zy)/\bar{F}_{t}(y)-z^{-1/\gamma(t)}}{A_{t}(1/\bar{F}_{t}(y))} = z^{-1/\gamma(t)}\frac{z^{\rho(t)/\gamma(t)}-1}{\gamma(t)\rho(t)} \text{~~for any~} z>0.
		\end{equation}
	\end{assumption} 
	
	According to Theorem 2.3.9 of \citet{deHaan2006}, \eqref{eq:second_order} is equivalent to
	\begin{eqnarray}
		\notag \lim_{y\to\infty}\frac{q_{t}\{1-1/(zy)\}/q_{t}(1-1/y)-z^{\gamma(t)}}{A_{t}(y)} = z^{\gamma(t)}\frac{z^{\rho(t)}-1}{\rho(t)} \text{~~for any~} z>0.
	\end{eqnarray}
	If $F_{t}(y)$ belongs to the generalized Hall class of heavy-tailed distribution as in \eqref{eq:Hall_class}, then $\bar{F}_{t}(y)$ satisfies Assumption \ref{assump:second_order} with index function $\gamma(\cdot)$, $\rho(t)=-\beta(t)\gamma(t)$, and $A_t(y)\asymp y^{\rho(t)}$.
	
	\begin{assumption}\label{asaump:alpha_N} 
		Let $\alpha_N\to \infty$, $h\to 0$, $Nh(1-\alpha_N)\to\infty$, $Nh^3(1-\alpha_N)\to0$, $Nh(1-\alpha_N)\delta_N^2\to 0$ as $N\to\infty$. In addition, for any fixed $t\in\mathcal{T}$, $\omega_t(q_{t}(\alpha_N),h)\log(q_{t}(\alpha_N))\to0$,  $\sqrt{Nh(1-\alpha_N)}\omega_t(q_{t}(\alpha_N),h)\log(q_{t}(\alpha_N))\to0$ and $\sqrt{Nh(1-\alpha_N)}A_{t}(1/(1-\alpha_N))\to0$.
	\end{assumption}
	
	The asymptotic property of $\hat{q}_{t, h}(\alpha_N)$ is given in the following theorem. Let $\Delta_0<1$ be any positive fixed number.
	
	\begin{thm}\label{theorem:quantile_intermediate}
		Suppose Assumptions \ref{assump:extreme_domain}, \ref{UnconfoundAssump}, \ref{as:suppX}, \ref{assump:f_T}, \ref{asaump:K}, \ref{assump:pi}, \ref{assump:second_order} and \ref{asaump:alpha_N} hold. 
		Denote for $\upsilon\in[\Delta_0, 1]$,
		\begin{equation}
			\mathcal{Q}_{t,h}(\upsilon) = \sqrt{Nh\upsilon(1-\alpha_N)}\left\{\frac{\hat{q}_{t, h}\{1-\upsilon(1-\alpha_N)\}}{q_{t}(1-\upsilon(1-\alpha_N))}-1\right\}.
		\end{equation}
		If there exists a function $\varpi_{t}^{Q}(\upsilon_1,\upsilon_2)$ such that
		\begin{eqnarray}
		    \notag \varpi_{t}^{Q}(\upsilon_1,\upsilon_2)=\lim_{\alpha_N\to 1}\frac{ \mathbb{E}\left[\left\{\pi_{0}(T,\mathbf{X})\right\}^2\mathbbm{1}[Y>q_{t}\{1-\upsilon_1(1-\alpha_N)\}]\mathbbm{1}[Y>q_{t}\{1-\upsilon_2(1-\alpha_N)\}]|T=t\right]}{\sqrt{\upsilon_1\upsilon_2}(1-\alpha_N)},
		\end{eqnarray}
		let $\Omega_{t}^{Q}(\upsilon_1,\upsilon_2)=\left\{\gamma(t)\right\}^2\kappa_{02} \{f_{T}(t)\}^{-1}\varpi_{t}^{Q}(\upsilon_1,\upsilon_2)$.
		Then, for any fixed $t_1\neq t_2\in\mathcal{T}$, 
		\begin{eqnarray}
			\notag {\mathcal{Q}_{t_1,h}(\upsilon) \choose \mathcal{Q}_{t_2,h}(\upsilon)} \Rightarrow \Psi_{t_1,t_2}^{Q}(\upsilon),
		\end{eqnarray}
		in $\ell^{\infty}([\Delta_0,1])$ as $N\to\infty$, where $\Psi_{t_1,t_2}^{Q}(\upsilon)$ is a centered Gaussian process with covariance function
		\begin{eqnarray}
			\Omega_{t_1, t_2}^{Q}(\upsilon_{1}, \upsilon_{2}) = \left( 
			\begin{array}{cc}
				\Omega_{t_1}^{Q}(\upsilon_1,\upsilon_2) & 0 \\
				0 & \Omega_{t_2}^{Q}(\upsilon_1,\upsilon_2)
			\end{array}
			\right).
		\end{eqnarray}
	\end{thm}
	Theorem \ref{theorem:quantile_intermediate} indicates that $\hat{q}_{t_1,h}(\alpha_N)$ is asymptotic independent of $\hat{q}_{t_2,h}(\alpha_N)$ for any $t_1\neq t_2$, which is essential in studying the EQTE.
	According to the proof of Theorem \ref{theorem:quantile_intermediate}, 
	\begin{eqnarray}
		\notag \sup_{\upsilon\in[\Delta_0,1]}\left|\mathcal{Q}_{t,h}(\upsilon)-\sum_{i=1}^{N}\gamma(t)\psi_{t,h}^{Q}(T_i, \mathbf{X}_i, Y_i, \upsilon, \alpha_N)\right| = o_{p}(1).
	\end{eqnarray}
	where
	\begin{eqnarray}
		\notag \psi_{t,h}^{Q}(T_i, \mathbf{X}_i, Y_i, \upsilon, \alpha_N) = \sqrt{\frac{h\upsilon(1-\alpha_N)}{N}}\left\{\frac{\pi_0(T_i,\mathbf{X}_i)\mathbbm{1}[Y_i> q_{t}\{1-\upsilon(1-\alpha_N)\}]K_h\left(T_i-t\right)}{\upsilon(1-\alpha_N)f_{T}(t)}-1\right\}.
    \end{eqnarray}
	It is straightforward to estimate $\Omega_{t}^{Q}(\upsilon_{1}, \upsilon_{2})$ using
	\begin{eqnarray}\label{def:Omegahat}
		\hat{\Omega}_{t}^{Q}(\upsilon_{1}, \upsilon_{2}) = \{\hat{\gamma}(t)\}^2\sum_{i=1}^{N} \hat{\psi}_{t,h}^{Q}(T_i, \mathbf{X}_i, Y_i, \upsilon_1, \alpha_N) \hat{\psi}_{t,h}^{Q}(T_i, \mathbf{X}_i, Y_i, \upsilon_2, \alpha_N),
	\end{eqnarray}
	where 
	\begin{eqnarray}
		\notag  \hat{\psi}_{t,h}^{Q}(T_i, \mathbf{X}_i, Y_i, \upsilon, \alpha_N) = \sqrt{\frac{h\upsilon(1-\alpha_N)}{N}}\left\{\frac{\hat{\pi}(T_i,\mathbf{X}_i)\mathbbm{1}[Y_i>\hat{q}_{t,h}\{1-\upsilon(1-\alpha_N)\}]K_h\left(T_i-t\right)}{\upsilon(1-\alpha_N)\cdot N^{-1}\sum^N_{j=1}K_h(T_j-t)}-1\right\},
	\end{eqnarray}
    and $\hat{\gamma}(t)$ is a consistent estimator of $\gamma(t)$ which will be provided in the next section. 
	%The consistency of $\hat{\Omega}_{1,t}(\upsilon_{1}, \upsilon_{2})$ directly follows from Theorem \ref{theo:Omega_est}.\footnote{need more technique!}
	The following theorem shows that $\hat{\Omega}_{t}^{Q}(\upsilon_1,\upsilon_2)$ is consistent for $\Omega_{t}^{Q}(\upsilon_1,\upsilon_2)$.
	
	\begin{thm}\label{theo:Omega_Q_est}
		Under the conditions assumed in Theorem \ref{theorem:quantile_intermediate}. In addition, assume that $\hat{\gamma}(t)$ is a consistent estimator of $\gamma(t)$. Then, for any fixed $t\in\mathcal{T}$ and $\upsilon_1,\upsilon_2\in[\Delta_0,1]$, $\hat{\Omega}_{t}^{Q}(\upsilon_1,\upsilon_2)/\Omega_{t}^{Q}(\upsilon_1,\upsilon_2) \to 1$
		in probability as $N\to\infty$.
	\end{thm}

	\subsection{Estimation of the Tail Index Function $\gamma(\cdot)$}\label{sec:gamma_est}
	
	It is necessary to derive a consistent estimator for $\gamma(t)$ in order to estimate the variance of $\mathcal{Q}_{t,h}(\upsilon)$. Following the technique of \citet{Daouia2013}, we employ the following Pickands \citep{Pickands1975} type estimator:
	\begin{eqnarray}\label{eq:pickands}
		\hat{\gamma}^{P}_{k_N}(t) = \frac{1}{\log(2)}\log\left(\frac{\hat{q}_{t, h}\{1-k_N/(4N)\}-\hat{q}_{t, h}\{1-k_N/(2N)\}}{\hat{q}_{t, h}\{1-k_N/(2N)\}-\hat{q}_{t, h}(1-k_N/N)}\right),
	\end{eqnarray}
	where $k_N$ is a sequence of diverging integers as $N\to\infty$ that satisfies Assumption \ref{asaump:alpha_N} with $1-\alpha_N=k_N/N$. The following theorem shows that $\hat{\gamma}_{k_N}^{P}(t)$ is consistent for $\gamma(t)$ and is asymptotic normal.
	
	\begin{thm}\label{theorem:gamma_estimate}
		Suppose Assumptions \ref{assump:extreme_domain}, \ref{UnconfoundAssump}, \ref{as:suppX}, \ref{assump:f_T}, \ref{asaump:K}, \ref{assump:pi}, \ref{assump:second_order} and \ref{asaump:alpha_N} with $1-\alpha_N=k_N/N$ hold. 
		Assume there exists a function $\varpi_{t}^{Q}(\upsilon_1,\upsilon_2)$ such that
		\begin{eqnarray}
		    \notag \varpi_{t}^{Q}(\upsilon_1,\upsilon_2)=\lim_{k_N\to \infty}\frac{ \mathbb{E}\left[\left\{\pi_{0}(T,\mathbf{X})\right\}^2\mathbbm{1}[Y>q_{t}(1-\upsilon_1k_N/N)]\mathbbm{1}[Y>q_{t}(1-\upsilon_2k_N/N)]|T=t\right]}{\sqrt{\upsilon_1\upsilon_2}k_N/N}.
		\end{eqnarray}
		For any fixed $t\in\mathcal{T}$, let $\Gamma_{P}(t)$ be a Gaussian random variance with mean $0$ and variance
		\begin{eqnarray}
			\notag & & \frac{\{\gamma(t)\}^2\kappa_{02}\{f_{T}(t)\}^{-1} }{\{\log(2)\}^2\{2^{\gamma(t)}-1\}^2}                       
			\left\{2^{2\gamma(t)+2}\varpi_{t}^{Q}(4, 4)-2^{\gamma(t)+5/2}(1+2^{\gamma(t)})\varpi_{t}^{Q}(4, 2)\right. \\
			\notag & & ~~~~~~~~~~~~~~~~~~~~~~~~~~~~~\left. +2^{\gamma(t)+2}\varpi_{t}^{Q}(4, 1)+2(1+2^{\gamma(t)})^2\varpi_{t}^{Q}(2,2) \right. \\
			\notag & & ~~~~~~~~~~~~~~~~~~~~~~~~~~~~~~~\left.-2\sqrt{2}(1+2^{\gamma(t)})\varpi_{t}^{Q}(2, 1)+ \varpi_{t}^{Q}(1, 1) \right\}.
		\end{eqnarray}
		Then for any fixed $t_1\neq t_2 \in \mathcal{T}$,
		\begin{eqnarray}
			\notag {\sqrt{k_Nh}\left\{\hat{\gamma}_{k_N}^{P}(t_1)-\gamma(t_1)\right\} \choose \sqrt{k_Nh}\left\{\hat{\gamma}_{k_N}^{P}(t_2)-\gamma(t_2)\right\}} \to {\Gamma_{P}(t_1) \choose \Gamma_{P}(t_2)},
		\end{eqnarray}
		in distribution, where $\Gamma_{P}(t_1)$ and $\Gamma_{P}(t_2)$ are independent.
	\end{thm}

	There are several disadvantages of the Pickands type estimator $\hat{\gamma}_{k_N}^{P}(t)$. First, the variance of $\hat{\gamma}_{k_N}^{P}(t)$ is huge when the underlying true value of the tail index function $\gamma(t)$ is large \citep{Daouia2011}. Second, the value of $\hat{\gamma}_{K_N}^{P}(t)$ may be negative for some real datasets. Thus, we consider a kernel version of the Hill estimator \citep{Hill1975} defined as
	\begin{eqnarray}\label{eq:hill}
		\hat{\gamma}_{k_N}^{H}(t) = \frac{1}{\sum_{j=1}^{J}\log(1/\upsilon_j)}\sum_{j=1}^{J}\log\left(\frac{\hat{q}_{t, h}(1-\upsilon_j k_N/N)}{\hat{q}_{t, h}(1-k_N/N)}\right),
	\end{eqnarray}
	where $1=\upsilon_1>\upsilon_2>\cdots>\upsilon_J>0$ is a decreasing list of $J$ weights. We have the next theorem that establishes the asymptotic distribution of $\hat{\gamma}_{k_N}^{H}(t)$.
	
	\begin{thm}\label{theorem:hill}
		Suppose the conditions assumed in Theorem \ref{theorem:gamma_estimate} hold.
		For any fixed $t\in\mathcal{T}$, let $\Gamma_{H}(t)$ be a Gaussian random variance with mean $0$ and variance
		\begin{eqnarray}
			\notag & & \frac{1}{\left\{\sum_{j=1}^{J}\log(1/\upsilon_j)\right\}^2}\Theta_{J}^{\top}\Sigma_{J}(t)\Theta_{J}.
		\end{eqnarray}
		where $\Theta_{J} = \left(1-J, \upsilon_{2}^{-1/2},\ldots,\upsilon_{J}^{-1/2}\right)^{\top}$ and $\Sigma_{J}(t)=(\sigma_{j_1,j_2}(t))_{1\leq j_1,j_2\leq J}$ with
		\begin{eqnarray}
			\notag \sigma_{j_1,j_2}(t) = \Omega_{t}^{Q}(\upsilon_{j_1},\upsilon_{j_2}) = \left\{\gamma(t)\right\}^2\kappa_{02} \{f_{T}(t)\}^{-1}\varpi_{t}^{Q}(\upsilon_{j_1},\upsilon_{j_2}),
		\end{eqnarray}
		where $\varpi_{t}^{Q}(\upsilon_1,\upsilon_2)=\lim_{k_N\to \infty}\frac{ \mathbb{E}\left[\left\{\pi_{0}(T,\mathbf{X})\right\}^2\mathbbm{1}[Y>q_{t}(1-\upsilon_1k_N/N)]\mathbbm{1}[Y>q_{t}(1-\upsilon_2k_N/N)]|T=t\right]}{\sqrt{\upsilon_1\upsilon_2}k_N/N}$.
		Then for any fixed $t_1\neq t_2\in \mathcal{T}$,
		\begin{eqnarray}
			\notag {\sqrt{k_Nh}\left\{\hat{\gamma}_{k_N}^{H}(t_1)-\gamma(t_1)\right\} \choose \sqrt{k_Nh}\left\{\hat{\gamma}_{k_N}^{H}(t_2)-\gamma(t_2)\right\}} \to {\Gamma_{H}(t_1) \choose \Gamma_{H}(t_2)},
		\end{eqnarray}
		in distribution, where $\Gamma_{H}(t_1)$ and $\Gamma_{H}(t_2)$ are independent.
	\end{thm}
	
	With $\hat{\Omega}_{t}^{Q}(\upsilon_1,\upsilon_2)$ a consistent estimator of $\Omega_{t}^{Q}(\upsilon_1,\upsilon_2)$ (e.g.~\eqref{def:Omegahat} with $\alpha_N$ replaced by $1-k_N/N$), the variance of $\Gamma_{H}(t)$ can be estimated straightforwardly by plug-in method.

	\subsection{Extreme Quantile}
	
	In this section, we study estimating $q_{t}(\alpha_N)$ in the extreme case when $\alpha_N\to1$ and $Nh(1-\alpha_N)\to c\in[0,\infty)$. In this case, the quantile level is far beyond the observed sample. Thus, we employ the commonly used extrapolation technique in extreme value theory. Theorem 1.1.11 of \citet{deHaan2006} implies that
	\begin{eqnarray}\label{eq:tail_index_quantile}
		\lim_{\alpha_N\to1}\frac{q_{t}\{1-z(1-\alpha_N)\}}{q_{t}(\alpha_N)}=z^{-\gamma(t)}, \text{~~for any~} 0<z<1.
	\end{eqnarray}
	This suggests that when $\alpha_N$ and $\beta_N$ satisfy $Nh(1-\alpha_N)\to c$ and $Nh(1-\beta_N)\to\infty$, we may write
	\begin{eqnarray}\label{eq:quantileratio}
		q_{t}(\alpha_N) \approx q_{t}(\beta_N)\left(\frac{1-\beta_N}{1-\alpha_N}\right)^{\gamma(t)}.
	\end{eqnarray}
	When $\beta_N$ is an intermediate sequence, we can use $\hat{q}_{t, h}(\beta_N)$ defined in Section \ref{sec:intermediate} to estimate $q_{t}(\beta_N)$. Let $k_N$ satisfies $k_N\to\infty$ and $k_N/N\to0$,  we set $\beta_N=1-k_N/N$ such that $Nh(1-\beta_N)=k_Nh\to\infty$. Then the extreme quantile $q_{t}(\alpha_N)$ is estimated via
	\begin{eqnarray}\label{eq:extreme_estimator}
		\hat{q}_{t, h}^{E}(\alpha_N) = \hat{q}_{t, h}(1-k_N/N)\left(\frac{k_N}{N(1-\alpha_N)}\right)^{\hat{\gamma}(t)},
	\end{eqnarray}
	where $\hat{\gamma}(t)$ is any consistent estimator of $\gamma(t)$. 
	
	The following theorem offers the asymptotic result for $\hat{q}_{t, h}^{E}(\alpha_N)$.
	
	\begin{thm}\label{theo:extreme_quantile}
		Suppose Assumptions \ref{assump:extreme_domain}, \ref{UnconfoundAssump}, \ref{as:suppX}, \ref{assump:f_T}, \ref{asaump:K}, \ref{assump:pi} and \ref{assump:second_order}  hold. Let $\alpha_N\to 1$, $h\to 0$, $k_N\to\infty$ such that $Nh(1-\alpha_N)\to c\in[0,\infty)$, $k_N/N\to0$, $k_Nh\to\infty$, $k_Nh^3\to0$ and $k_Nh\delta_N^2\to 0$ as $N\to\infty$. In addition, for any fixed $t\in\mathcal{T}$, assume that
		$\omega_t(q_{t}(1-k_N/N),h)\log(k_N/N)\to0$ $\sqrt{k_Nh}\omega_t(q_{t}(1-k_N/N),h)\log(k_N/N)\to0$, $\sqrt{k_Nh}A_{t}(N/k_N)\to0$
		and $\sqrt{k_Nh}/\log\left(k_N/\{N(1-\alpha_N)\}\right)\to\infty$. 
		Moreover, $\hat{\gamma}(t)$ is an estimator of $\gamma(t)$ that satisfies for any $t_1\neq t_2\in\mathcal{T}$,
		\begin{eqnarray}
			\notag {\sqrt{k_Nh}\left\{\hat{\gamma}(t_1)-\gamma(t_1)\right\} \choose \sqrt{k_Nh}\left\{\hat{\gamma}(t_2)-\gamma(t_2)\right\}} \to {\Gamma(t_1) \choose \Gamma(t_2)},
		\end{eqnarray}
		in distribution, where $\Gamma(t_1)$ and $\Gamma(t_2)$ are independent.
		in distribution, where $\Gamma$ is a non-degenerate distribution.
		Denote for $\upsilon\in[\Delta_0,1]$,
		\begin{equation}
			\mathcal{Q}_{t,h}^{E}(\upsilon) = \frac{\sqrt{k_Nh}}{\log\left(k_N/\{N\upsilon(1-\alpha_N)\}\right)}\left\{\frac{\hat{q}_{t, h}^{E}(1-\upsilon\{1-\alpha_N)\}}{q_{t}\{1-\upsilon(1-\alpha_N)\}}-1\right\}.
		\end{equation}
		Then, for any fixed $t_1\neq t_2\in\mathcal{T}$,
		\begin{eqnarray}
			\notag {\mathcal{Q}^{E}_{t_1,h}(\upsilon) \choose \mathcal{Q}^{E}_{t_2,h}(\upsilon)} \Rightarrow {\Gamma(t_1) \choose \Gamma(t_2)},
		\end{eqnarray}
		in $\ell^{\infty}([\Delta_0,1])$.
	\end{thm}
	Theorem \ref{theo:extreme_quantile} indicates that the extreme quantile estimator has the same asymptotic as the estimator of the tail index, but with a slower convergence rate.
	
	\subsection{Extreme Tail Mean}
	
	In order to estimate the EATE, it is necessary to estimate the tail mean of $Y(t)$.
    According to Proposition 4.1 of \citet{Pan2013SPL},
    \begin{eqnarray}
        \notag \lim_{\alpha_N\to1}\frac{\mathrm{TM}_{t}(\alpha_N)}{q_{t}(\alpha_N)} = \frac{1}{1-\gamma(t)}.
    \end{eqnarray}
    This motivates an estimator of $\mathrm{TM}_{t}(\alpha_N)$ as
    \begin{eqnarray}
        \widehat{\mathrm{TM}}_{t,h}(\alpha_N) = \frac{\hat{q}_{t,h}^{E}(\alpha_N)}{1-\hat{\gamma}(t)},
    \end{eqnarray}
    where $\hat{\gamma}(t)$ is an estimator of $\gamma(t)$.
    
    The asymptotic result of $\widehat{\mathrm{TM}}_{t}(\alpha)$ is given in the following theorem.
    
    \begin{thm}\label{theo:extreme_TM}
		Suppose the conditions of Theorem \ref{theo:extreme_quantile} hold. In addition, assume $\gamma(t)<1/2$ for any $t\in\mathcal{T}$.
		Denote for $\upsilon\in[\Delta_0,1]$,
		\begin{equation}
			\mathcal{M}_{t,h}(\upsilon) = \frac{\sqrt{k_Nh}}{\log\left(k_N/\{N\upsilon(1-\alpha_N)\}\right)}\left\{\frac{\widehat{\mathrm{TM}}_{t, h}\{1-\upsilon(1-\alpha_N)\}}{\mathrm{TM}_{t}\{1-\upsilon(1-\alpha_N)\}}-1\right\}.
		\end{equation}
		Then, for any fixed $t_1\neq t_2\in\mathcal{T}$,
		\begin{eqnarray}
			\notag {\mathcal{M}_{t_1,h}(\upsilon) \choose \mathcal{M}_{t_2,h}(\upsilon)} \Rightarrow {\Gamma(t_1) \choose \Gamma(t_2)},
		\end{eqnarray}
		in $\ell^{\infty}([\Delta_0,1])$, where $\Gamma(t_1)$ and $\Gamma(t_2)$ are independent.
	\end{thm}
	
	The condition that $\gamma(t)<1/2$ for any $t\in\mathcal{T}$ is necessary to ensure the finite variance of $Y(t)$, which is commonly used in the tail mean literature \citet{Li2022ET}.
	
	\section{Estimate and Infer the EQTE and EATE}\label{sec:TE}
	
    With $\hat{q}_{t,h}(\alpha_N)$ and $\widehat{\mathrm{TM}}_{t,h}(\alpha_N)$ as the estimators of the quantile and tail mean, we can estimate the EQTE and EATE directly via plug-in method. For $t_1\neq t_2\in\mathcal{T}$, we estimate $\mathrm{EQTE}_{t_1,t_2}(\alpha_N)$ and $\mathrm{EATE}_{t_1,t_2}(\alpha_N)$ via
    \begin{eqnarray}
        \widehat{\mathrm{EQTE}}_{t_1,t_2}(\alpha_N) = \frac{\hat{q}_{t_1,h}(\alpha_N)}{\hat{q}_{t_2,h}(\alpha_N)}
    \end{eqnarray}
    and
    \begin{eqnarray}
        \widehat{\mathrm{EATE}}_{t_1,t_2}(\alpha_N) = \frac{\widehat{\mathrm{TM}}_{t_1,h}(\alpha_N)}{\widehat{\mathrm{TM}}_{t_2,h}(\alpha_N)}.
    \end{eqnarray}
    
    For the EQTE, we still consider two case scenarios: the intermediate case when $Nh(1-\alpha_N)\to\infty$ and the extreme case when $Nh(1-\alpha_N)\to c\in[0,\infty)$. 
    For the intermediate case, the quantile estimator in $\widehat{\mathrm{EQTE}}_{t_1,t_2}(\alpha_N)$ is $\hat{q}_{t,h}(\alpha_N)$ defined in \eqref{eq:intermediate_quantile_est}.
    For the extreme case, the corresponding quantile estimator is $\hat{q}_{t,h}^{E}(\alpha_N)$ defined in \eqref{eq:extreme_estimator}, and we denote $\widehat{\mathrm{EQTE}}_{t_1,t_2}^{E}(\alpha_N) = \hat{q}_{t_1,h}^{E}(\alpha_N)/\hat{q}_{t_2,h}^{E}(\alpha_N)$. Based on the previous theoretical results, we have the following theorem on the asymptotic properties of $\widehat{\mathrm{EQTE}}_{t_1,t_2}(\alpha_N)$ and $\widehat{\mathrm{EATE}}_{t_1,t_2}(\alpha_N)$.
	
	\begin{thm}\label{theo:EQTE_EATE}
	    (i) Under the conditions assumed in Theorem \ref{theorem:quantile_intermediate} in which $Nh(1-\alpha_N)\to\infty$, we have
	    \begin{eqnarray}
	        \sqrt{Nh\upsilon(1-\alpha_N)}\left\{\frac{\widehat{\mathrm{EQTE}}_{t_1,t_2}\{1-\upsilon(1-\alpha_N)\}}{\mathrm{EQTE}_{t_1,t_2}\{1-\upsilon(1-\alpha_N)\}} - 1 \right\} \Rightarrow \Psi_{t_1,t_2}^{\mathrm{QTE}}(\upsilon)
	    \end{eqnarray}
	    in $\ell^{\infty}[\Delta_0,1]$, where $\Psi_{t_1,t_2}^{\mathrm{QTE}}(\upsilon)$ is a centered process with covariance function
	    \begin{eqnarray}
	        \notag \Omega_{t_1}^{Q}(\upsilon_1,\upsilon_2)+\Omega_{t_2}^{Q}(\upsilon_1,\upsilon_2).
	    \end{eqnarray}
	    (ii) Under the conditions assumed in Theorem \ref{theo:extreme_quantile} in which $Nh(1-\alpha_N)\to c\in[0,\infty)$, we have  
	    \begin{eqnarray}
	        \frac{\sqrt{k_Nh}}{\log(k_N/\{N\upsilon(1-\alpha_N)\})}\left\{\frac{\widehat{\mathrm{EQTE}}_{t_1,t_2}^{E}\{1-\upsilon(1-\alpha_N)\}}{\mathrm{EQTE}_{t_1,t_2}^{E}\{1-\upsilon(1-\alpha_N)\}} - 1 \right\} \Rightarrow \Gamma(t_1)+\Gamma(t_2)
	    \end{eqnarray}
	    in $\ell^{\infty}[\Delta_0,1]$, where $\Gamma(t_1)$ and $\Gamma(t_2)$ are independent. \\
	    (iii) Under the conditions assumed in Theorem \ref{theo:extreme_TM} in which $Nh(1-\alpha_N)\to c\in[0,\infty)$, we have 
	    \begin{eqnarray}
	        \frac{\sqrt{k_Nh}}{\log(k_N/\{N\upsilon(1-\alpha_N)\})}\left\{\frac{\widehat{\mathrm{EATE}}_{t_1,t_2}(1-\upsilon\{1-\alpha_N)\}}{\mathrm{EATE}_{t_1,t_2}\{1-\upsilon(1-\alpha_N)\}} - 1 \right\} \Rightarrow \Gamma(t_1)+\Gamma(t_2)
	    \end{eqnarray}
	    in $\ell^{\infty}[\Delta_0,1]$, where $\Gamma(t_1)$ and $\Gamma(t_2)$ are independent.
	\end{thm}
	
	According to Theorem \ref{theo:EQTE_EATE}, $\widehat{\mathrm{EQTE}}_{t_1,t_2}^{E}\{1-\upsilon(1-\alpha_N)\}$ and $\widehat{\mathrm{EATE}}_{t_1,t_2}\{1-\upsilon(1-\alpha_N)\}$ converges to a normal distribution uniformly over $\upsilon\in[\Delta_0,1]$. Thus, if $\hat{\sigma}^2_{\Gamma}(t)$ is a consistent estimator of the variance of $\Gamma(t)$ (e.g.~as the one defined below Theorem~\ref{theorem:hill}), we can construct $1-\alpha_0$ simultaneous confidence bands for $\mathrm{EQTE}_{t_1,t_2}^{E}(\varrho)$ and $\mathrm{EATE}_{t_1,t_2}(\varrho)$ over $\varrho\in[\Delta_0\alpha_N,\alpha_N]$ as
	\begin{eqnarray}
	    \notag & & \left[\widehat{\mathrm{EQTE}}_{t_1,t_2}^{E}(\varrho)\exp\left(-z_{1-\alpha_0/2}\sqrt{\hat{\sigma}^2_{\Gamma}(t_1)+\hat{\sigma}^2_{\Gamma}(t_2)}\frac{\log(k_N/\{N(1-\varrho)\})}{\sqrt{k_Nh}}\right),\right. \\
		\notag & & ~~~~~~~~\left.\widehat{\mathrm{EQTE}}_{t_1,t_2}^{E}(\varrho)\exp\left(z_{1-\alpha_0/2}\sqrt{\hat{\sigma}^2_{\Gamma}(t_1)+\hat{\sigma}^2_{\Gamma}(t_2)}\frac{\log(k_N/\{N(1-\varrho)\})}{\sqrt{k_Nh}}\right)\right],
	\end{eqnarray}
	and
	\begin{eqnarray}
	    \notag & & \left[\widehat{\mathrm{EATE}}_{t_1,t_2}(\varrho)\exp\left(-z_{1-\alpha_0/2}\sqrt{\hat{\sigma}^2_{\Gamma}(t_1)+\hat{\sigma}^2_{\Gamma}(t_2)}\frac{\log(k_N/\{N(1-\varrho)\})}{\sqrt{k_Nh}}\right),\right. \\
		\notag & & ~~~~~~~~\left.\widehat{\mathrm{EATE}}_{t_1,t_2}(\varrho)\exp\left(z_{1-\alpha_0/2}\sqrt{\hat{\sigma}^2_{\Gamma}(t_1)+\hat{\sigma}^2_{\Gamma}(t_2)}\frac{\log(k_N/\{N(1-\varrho)\})}{\sqrt{k_Nh}}\right)\right],
	\end{eqnarray}
	respectively. Here $z_{1-\alpha_0}$ is the $(1-\alpha_0)$-quantile of $N(0,1)$.
	The consistent estimators of the variances of the Pickands estimator $\hat{\gamma}_{k_N}^{P}(t)$ and the Hill estimator $\hat{\gamma}_{k_N}^{H}(t)$ can be obtained based on the discussions in Section \ref{sec:gamma_est}.
	
	Numerical results indicate that the simultaneous confidence band for $\mathrm{EATE}_{t_1,t_2}(\varrho)$ given above tends to be under covered. To adjust its finite-sample performance, we recall the expansion of $\widehat{\mathrm{TM}}_{t,h}(\alpha_N)$ in the proof of Theorem \ref{theo:extreme_TM}, that is,
	\begin{eqnarray}
	   \notag & & \frac{\sqrt{k_Nh}}{\log\left(k_N/\{N\upsilon(1-\alpha_N)\}\right)}\left\{\frac{\widehat{\mathrm{TM}}_{t, h}\{1-\upsilon(1-\alpha_N)\}}{\mathrm{TM}_{t}\{1-\upsilon(1-\alpha_N)\}}-1\right\} \\
	   \notag & = & \left\{1+\frac{1}{\log\left(k_N/\{N\upsilon(1-\alpha_N)\}\right)\{1-\gamma(t)\}}\right\}\sqrt{k_Nh}\{\hat{\gamma}(t)-\gamma(t)\} + o_{p}(1),
	\end{eqnarray}
	where $\frac{1}{\log\left(k_N/\{N\upsilon(1-\alpha_N)\}\right)\{1-\gamma(t)\}}\to0$ under the conditions of Theorem \ref{theo:extreme_TM}. However, this term could make a non-negligible contribution of the variance of $\widehat{\mathrm{TM}}_{t,h}(\alpha_N)$ in practice. Thus, we would like to recall it when estimating the variance of $\widehat{\mathrm{TM}}_{t,h}(\alpha_N)$. As a consequence, we adjust the simultaneous confidence band for $\mathrm{EATE}_{t_1,t_2}(\varrho)$ as
	\begin{eqnarray}
	    \notag & & \left[\widehat{\mathrm{EATE}}_{t_1,t_2}(\varrho)\exp\left(-z_{1-\alpha_0/2}\sqrt{\hat{\sigma}^2_{\mathrm{TM}}(t_1)+\hat{\sigma}^2_{\mathrm{TM}}(t_2)}\frac{\log(k_N/\{N(1-\varrho)\})}{\sqrt{k_Nh}}\right),\right. \\
		\notag & & ~~~~~~~~\left.\widehat{\mathrm{EATE}}_{t_1,t_2}(\varrho)\exp\left(z_{1-\alpha_0/2}\sqrt{\hat{\sigma}^2_{\mathrm{TM}}(t_1)+\hat{\sigma}^2_{\mathrm{TM}}(t_2)}\frac{\log(k_N/\{N(1-\varrho)\})}{\sqrt{k_Nh}}\right)\right],
	\end{eqnarray}
	where 
	\begin{eqnarray}
	    \notag \hat{\sigma}^2_{\mathrm{TM}}(t) = \left\{1+\frac{1}{\log\left(k_N/\{N(1-\varrho)\}\right)\{1-\hat{\gamma}(t)\}}\right\}^2\hat{\sigma}^2_{\Gamma}(t)
	\end{eqnarray}
	Simulation results in Section \ref{sec:simulation} show that this finite-sample adjustment improves the performance of the simultaneous confidence band for $\mathrm{EATE}_{t_1,t_2}(\varrho)$ significantly.

	\section{Numerical Study}\label{sec:Numerical_study}
	In this section, we apply our method in \eqref{eq:extreme_estimator} with $\widehat{\gamma}^H_{kN}(t)$ in \eqref{eq:hill} to estimate the extreme quantiles and exam the performance of the simultaneous confidence bands for the EQTE and EATE from two simulation models and a real data example. For estimating the extreme quantile, we compare our method to the naive one, which estimate the extreme quantile directly from \citeauthor{ai2021estimation}'s (\citeyear{ai2021estimation}) nonparametric estimator of the distribution function of $Y^*(t)$. That is,
	\begin{equation}
		\widehat{q}^{\text{Naive}}_{t,h}(\alpha_N) = \inf\{z: 1-\widehat{\bar{F}}_{t,h}(z)\geq\alpha_N\} = \inf\{z: \widehat{\bar{F}}_{t,h}(z)\leq1-\alpha_N\}.\label{Def:NaiveQuantiles}
	\end{equation}
	Recall that our estimator in \eqref{eq:extreme_estimator} requires an estimation of the weighting function $\pi_0$ and selection of tuning parameters $h$ and $k_N$. 
	
	Recall from \eqref{pi0def} that $\pi_0$ is a ratio of two densities. A straightforward way to estimate it is to estimate the two densities separately and then form a ratio estimator. However, it is known in the literature of treatment effect that such a ratio estimator is very unstable (see e.g.~\citealp{ai2021estimation,huang2021unified}). Thus, we adopt \citeauthor{ai2021estimation}'s (\citeyear{ai2021estimation}) estimator of $\pi_0$ and their cross-validation method to select $h$ and form the estimator $\widehat{\bar{F}}_{t,h}(y)$.
	
	We set $v_j =1 - (j-1)/J$ for $j=1,\ldots,J$ and select $k_N$ based on \eqref{eq:quantileratio}. Specifically, we first sort the $Y_i$'s in an descending order, denoted by $Z_1>\ldots> Z_N$. Let $\mathcal{K}$ be a candidate set of $k_N$. For each $\ell \in \mathcal{K}$, we calculate the distance
	$$
	D_{\ell} := \sup_{t\in\mathcal{T}}\max_{i=1,\ldots,\ell-1}\left|\frac{\widehat{\bar{F}}_{t,h}(Z_i)}{\widehat{\bar{F}}_{t,h}(Z_{\ell})} - \left(\frac{Z_{\ell}}{Z_i}\right)^{1/\widehat{\gamma}^H_{\ell}} \right|\,.
	$$
	Then we set
	$$
	k_N:= \min_{\ell\in\mathcal{K}} D_{\ell}\,.
	$$
	
	To ensure our $k_N$ and $h$ satisfy the conditions in Theorem~\ref{theo:extreme_quantile}, we set the candidate set of $h$ as $[N^{-1/3},h_{\text{ROT}}]$, where $h_{\text{ROT}}=1.06\sqrt{\text{var}(T)}\times N^{-1/5}$ is the rule of thumb bandwidth for kernel and $\mathcal{K}:=[J, \lfloor 0.2\times N^{0.95} \rfloor]$, where $\lfloor a \rfloor$ represents the floor function that returns the largest integer less than or equal to the real number $a$.

	\subsection{Simulation Studies}\label{sec:simulation}
	We consider the following two models. 
	\begin{eqnarray*}
		&&\textbf{DGP1} \quad \mathbb{P}\{Y(t)>y|X=x\} = (0.9+0.2x)\cdot y^{-1/\gamma(t)},\\
		&&\textbf{DGP2} \quad \mathbb{P}\{Y(t)>y|X=x\} = \Big\{-0.2x+1+y^{1/\gamma(t)}\Big\}^{-1},
	\end{eqnarray*}
	where the tail index function
	$$
	\gamma(t) = \frac{1}{4}+\frac{\sin(2\pi t)}{20}\,.
	$$
	For both models, we simulate the treatment data $T_i = 0.1(0.3+0.4X_i) -0.03 + 0.96\epsilon_i$, where $X_i$ and $\epsilon_i$ are independent a standard uniform random variable, and the observed outcome $Y_i=Y(T_i)$, for $i=1,\ldots,N$, where $N=500, 1000$ and 2000.
	
\begin{figure}
		\centering
		\includegraphics[width = .45\textwidth]{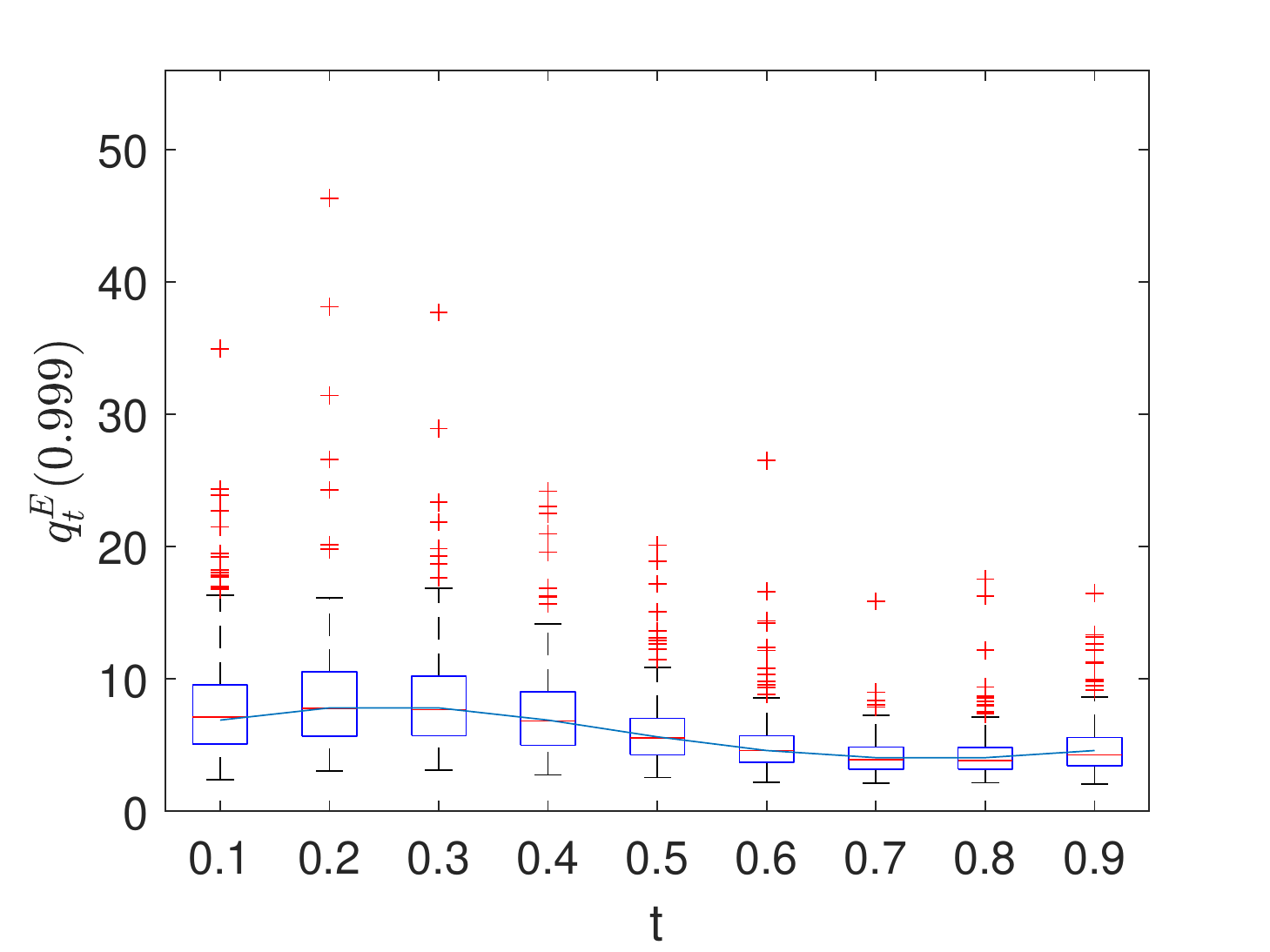}
		\includegraphics[width = .45\textwidth]{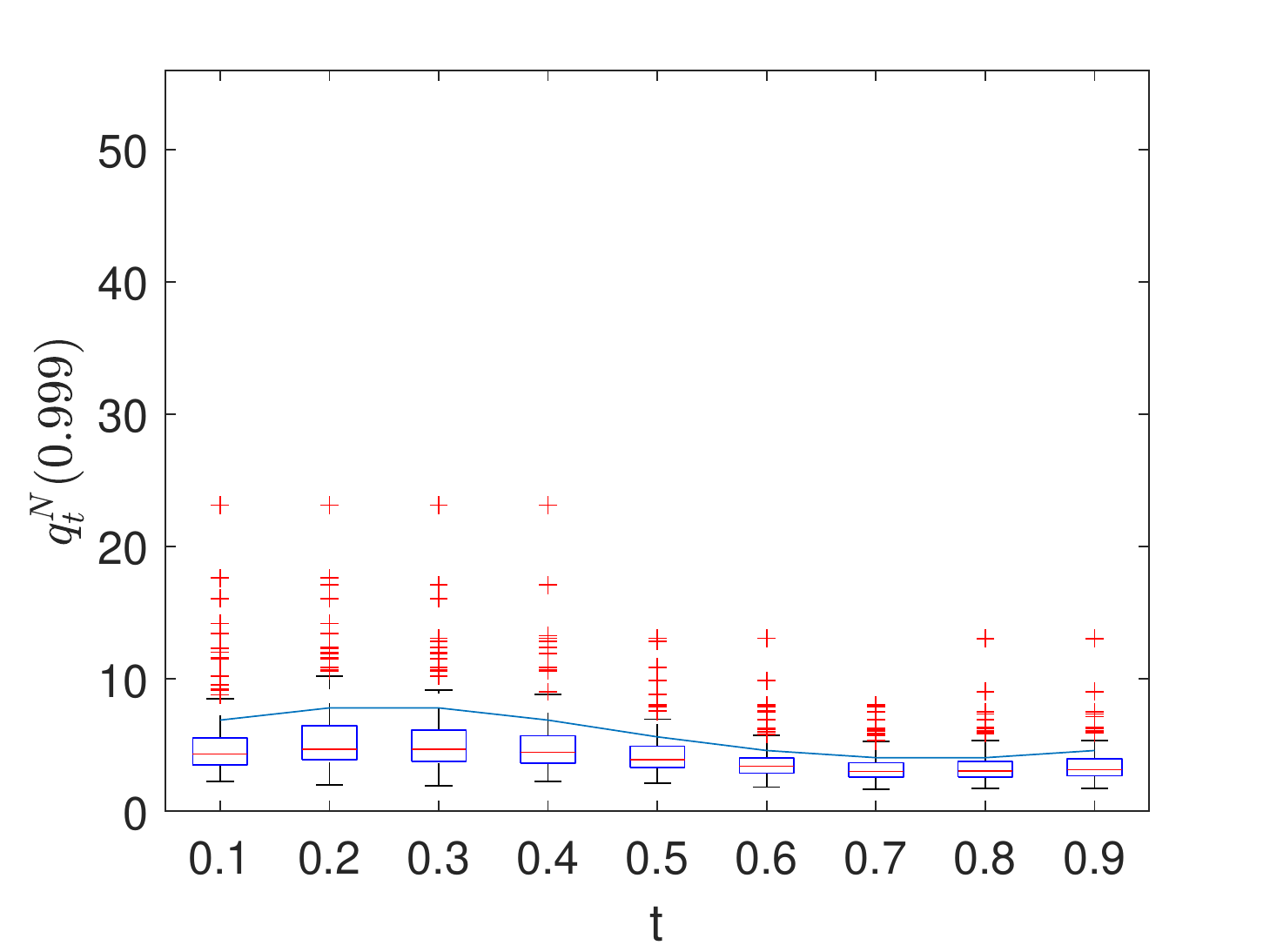}
		
		\includegraphics[width = .45\textwidth]{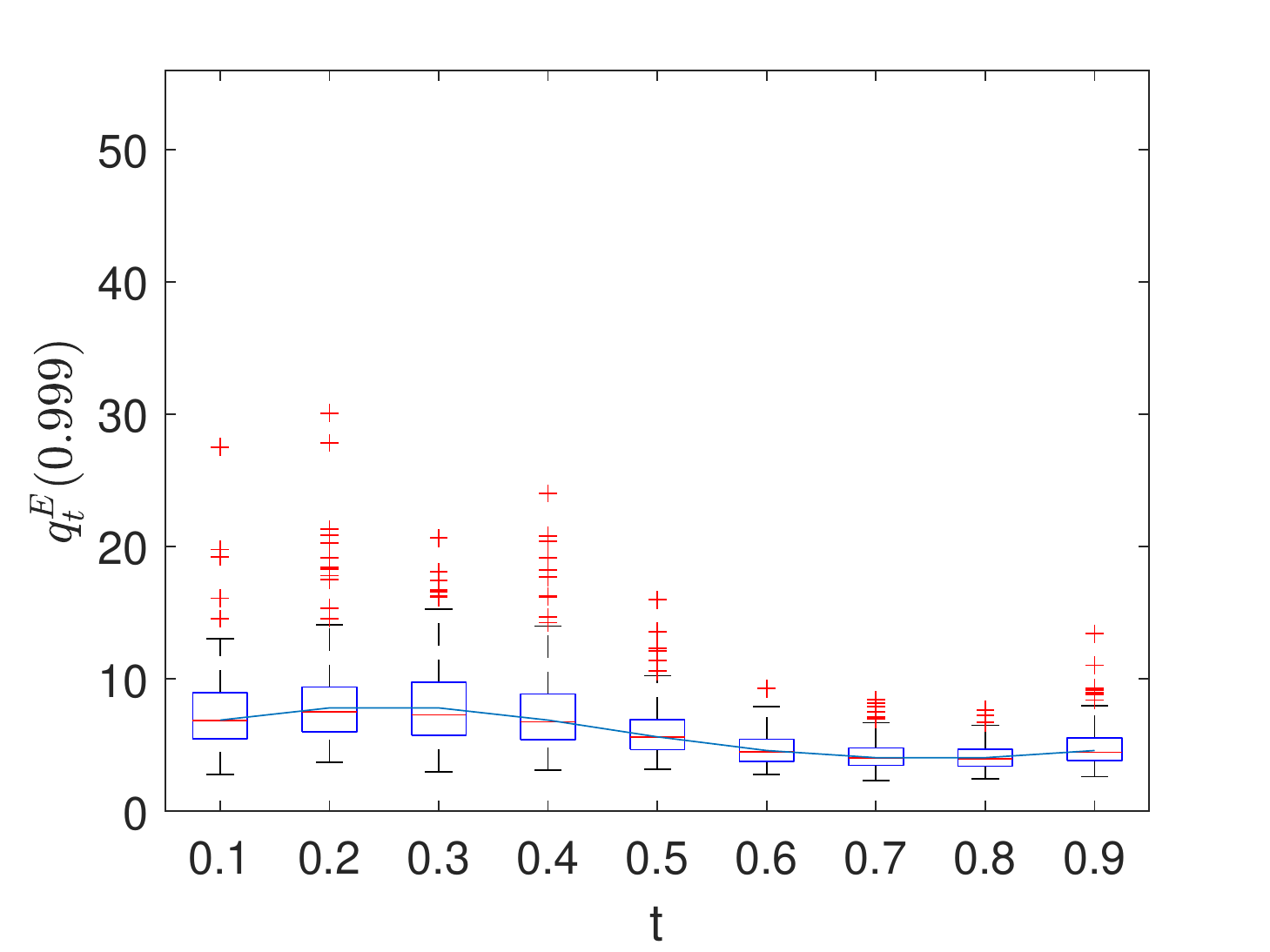}
		\includegraphics[width = .45\textwidth]{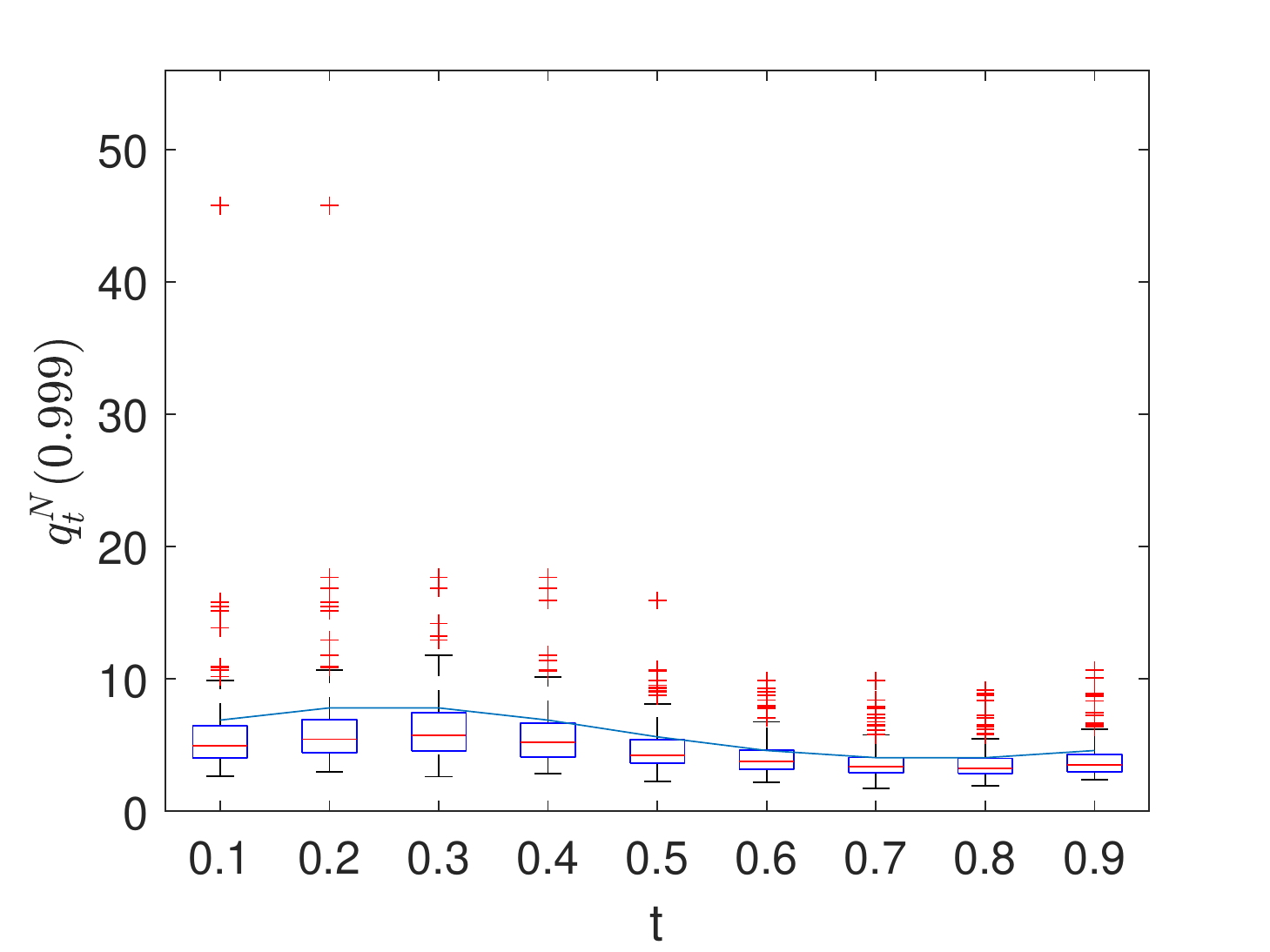}
		
		\includegraphics[width = .45\textwidth]{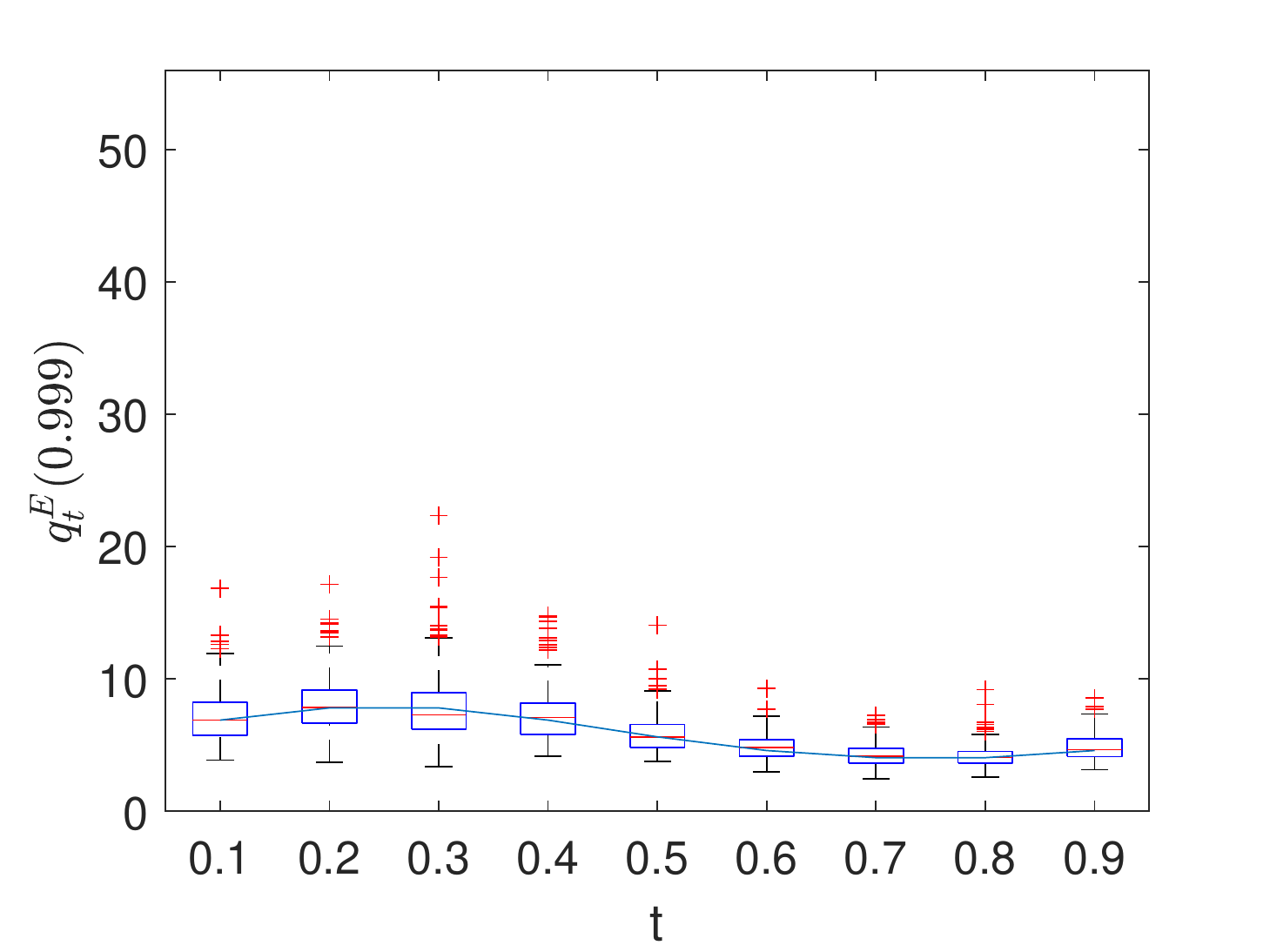}
		\includegraphics[width = .45\textwidth]{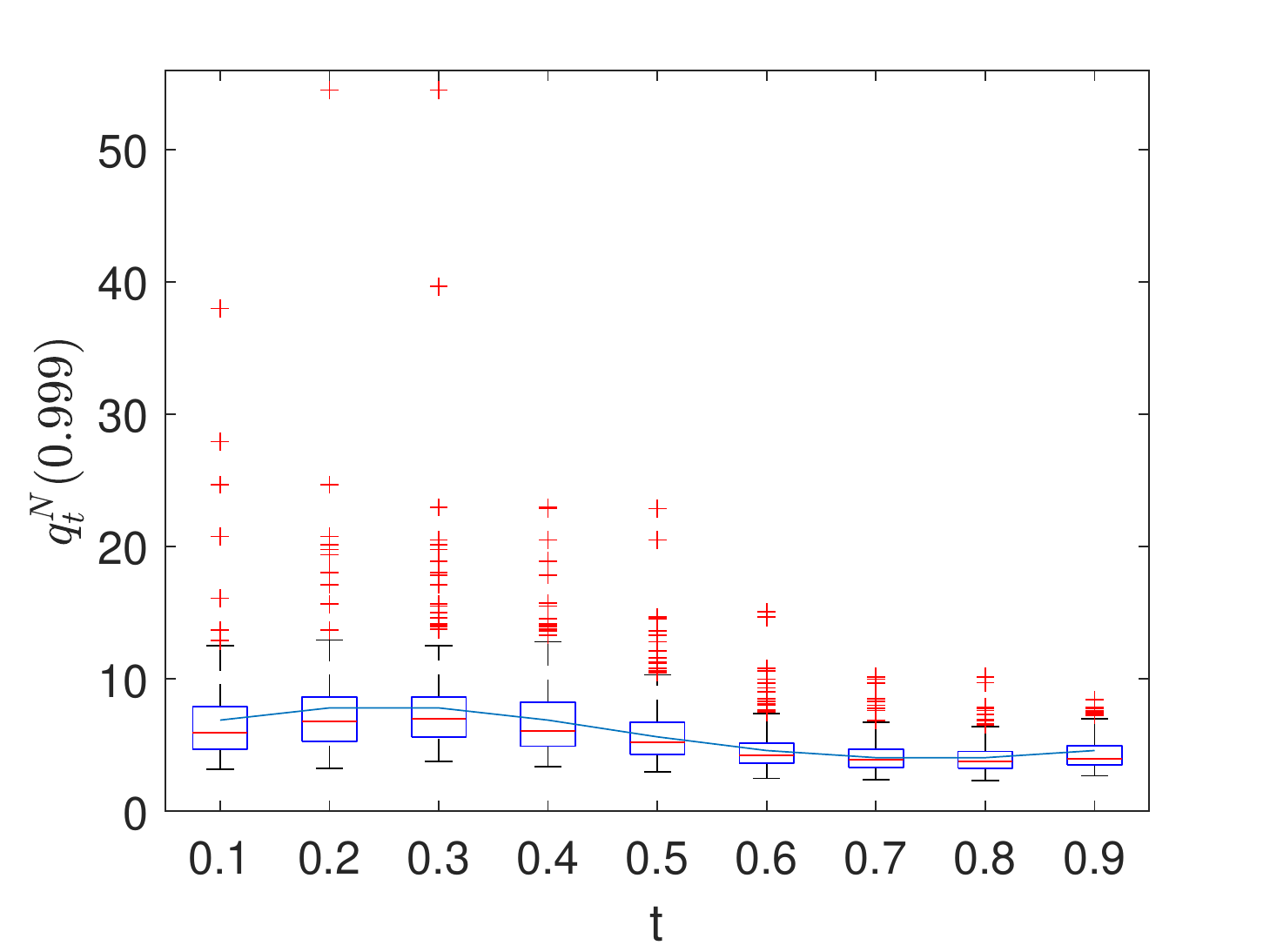}
		
		\caption{The boxplots of the extreme quantiles of $\alpha_N=0.999$ calculated using our method $\widehat{q}_t^E$ (left) and the naive method (right) from DGP1 with $N=500$ (row 1), 1000 (row 2), 2000 (row 3). }\label{fig:ExtremeQuantileDGP1}
	\end{figure}

\begin{figure}
		\centering
		\includegraphics[width = .45\textwidth]{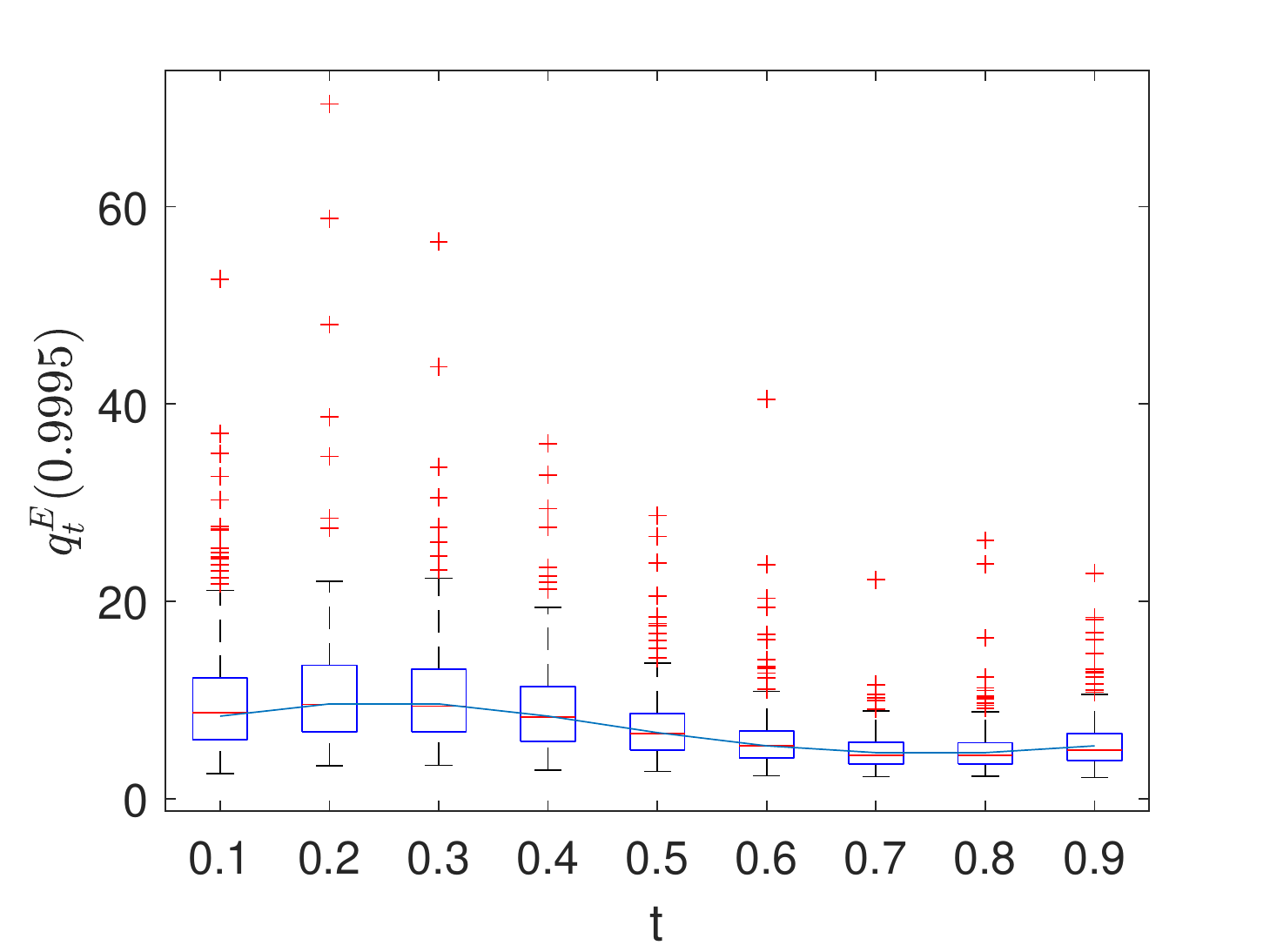}
		\includegraphics[width = .45\textwidth]{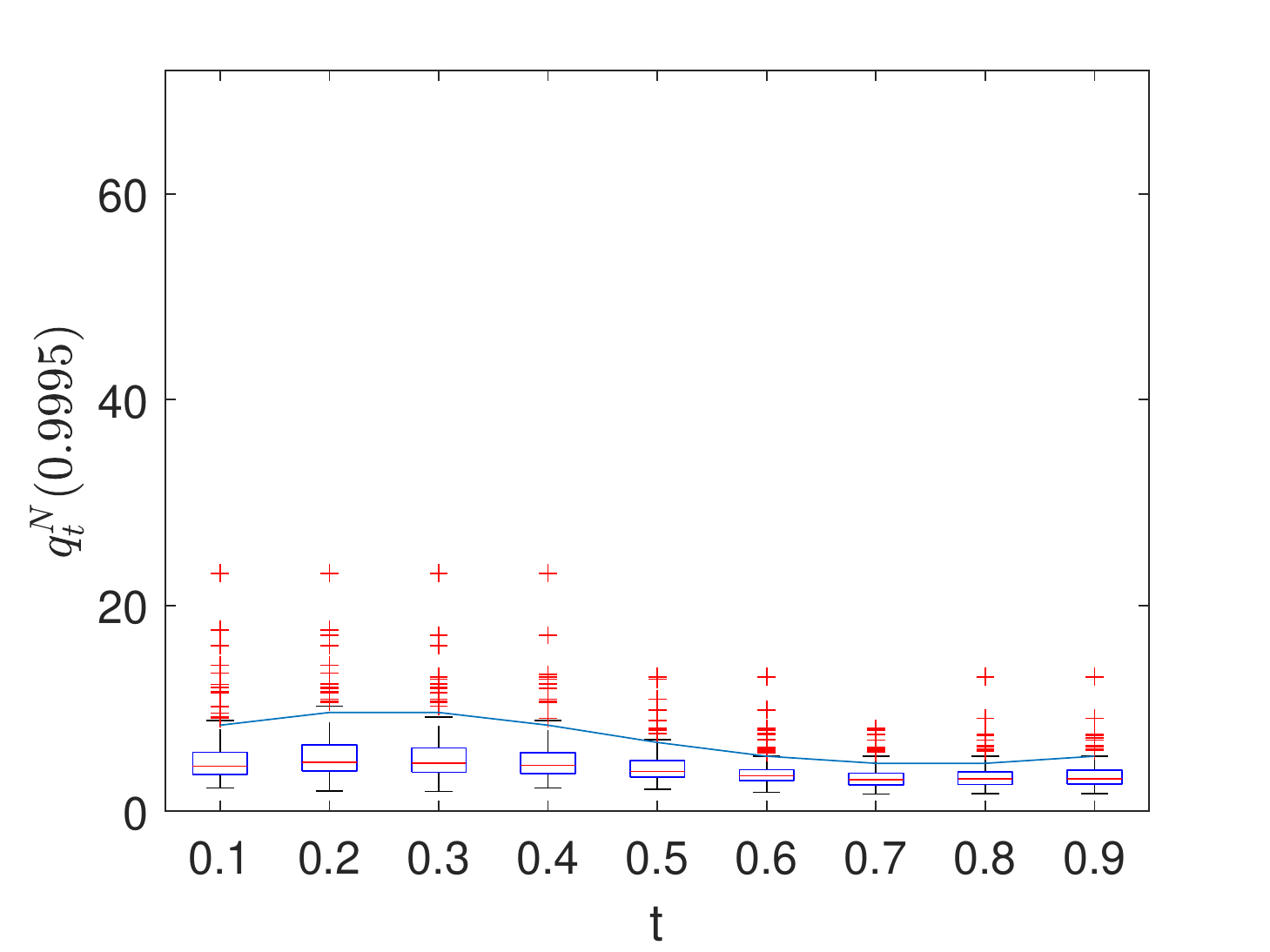}
		
		\includegraphics[width = .45\textwidth]{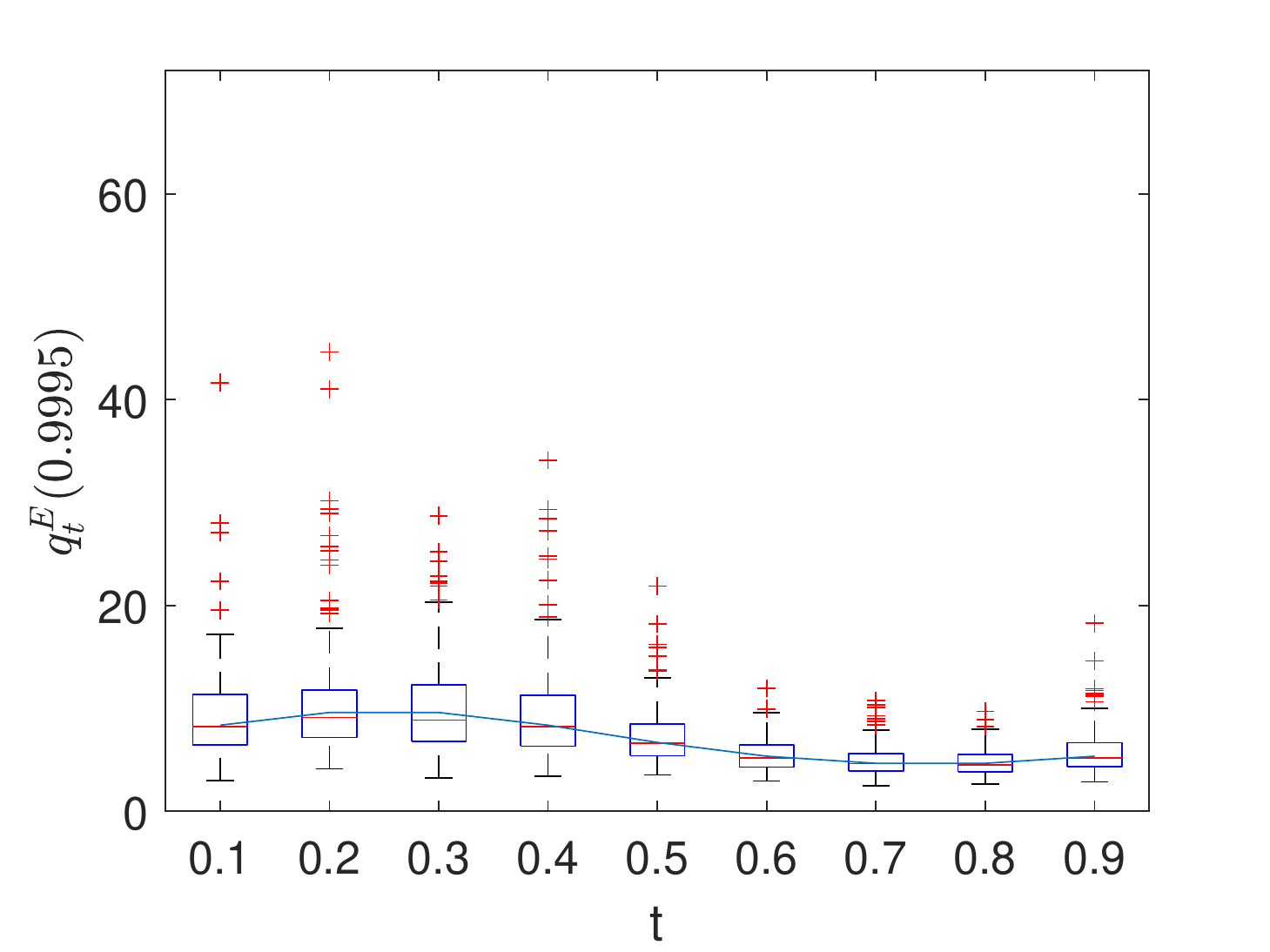}
		\includegraphics[width = .45\textwidth]{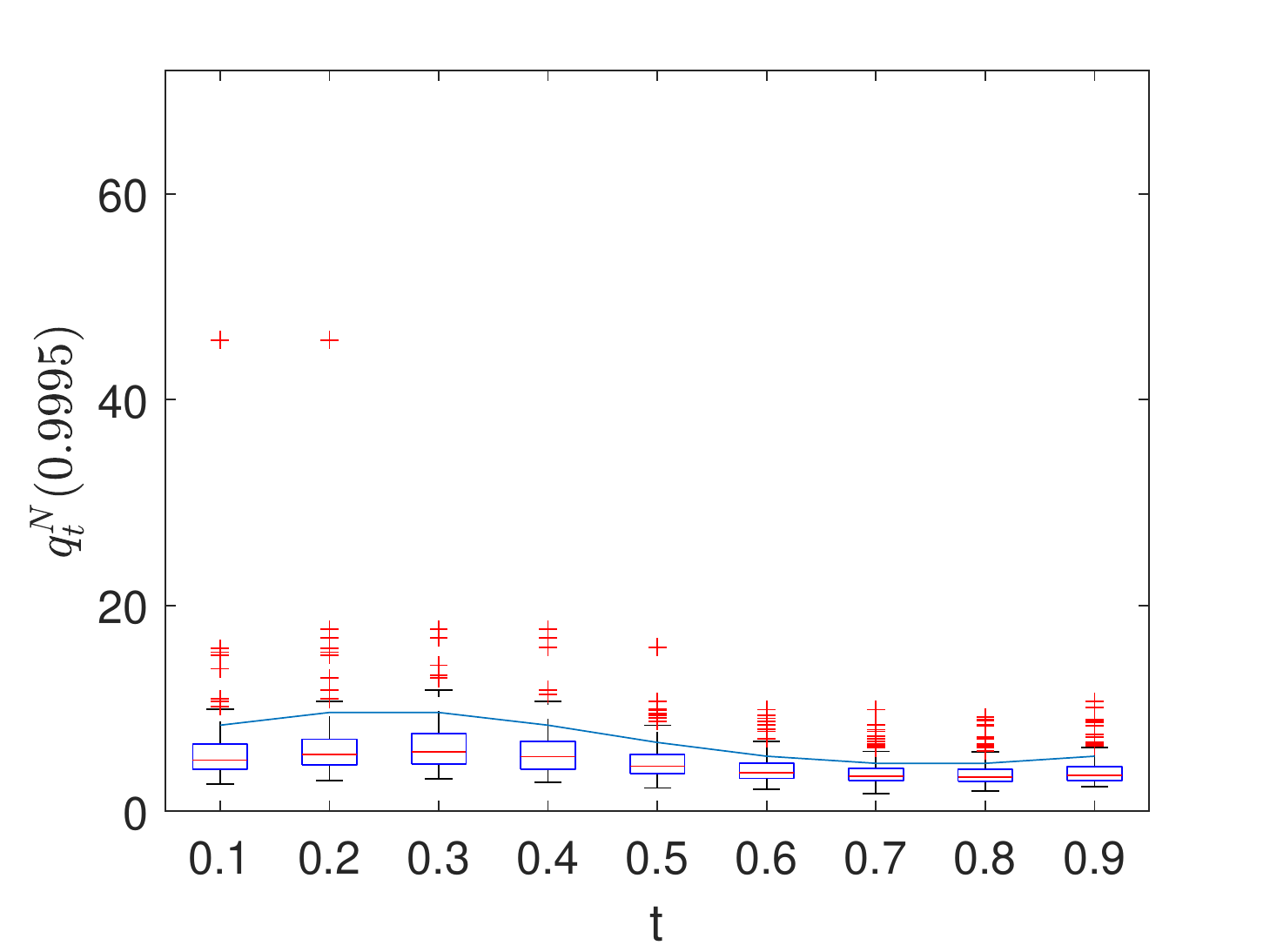}
		
		\includegraphics[width = .45\textwidth]{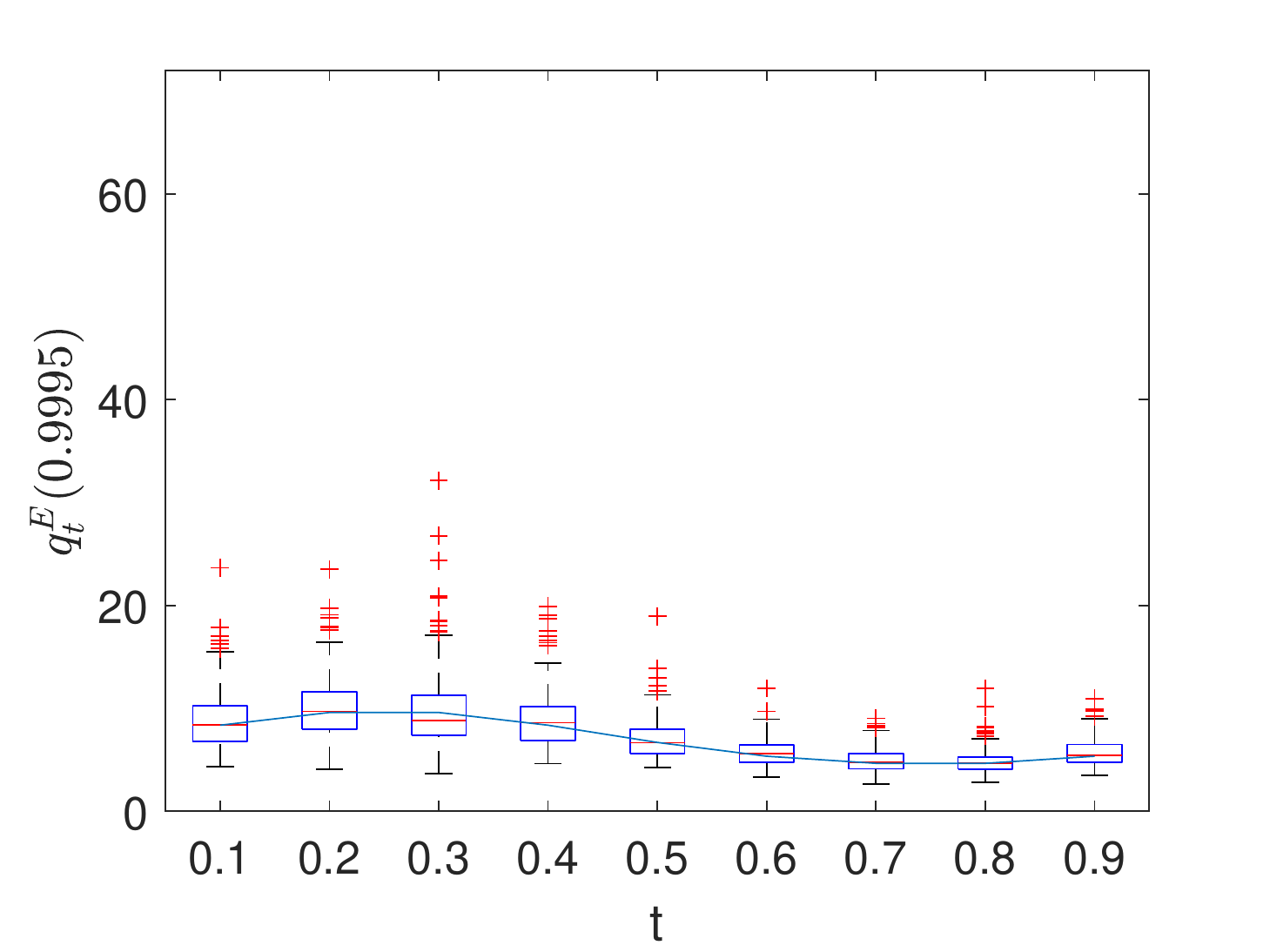}
		\includegraphics[width = .45\textwidth]{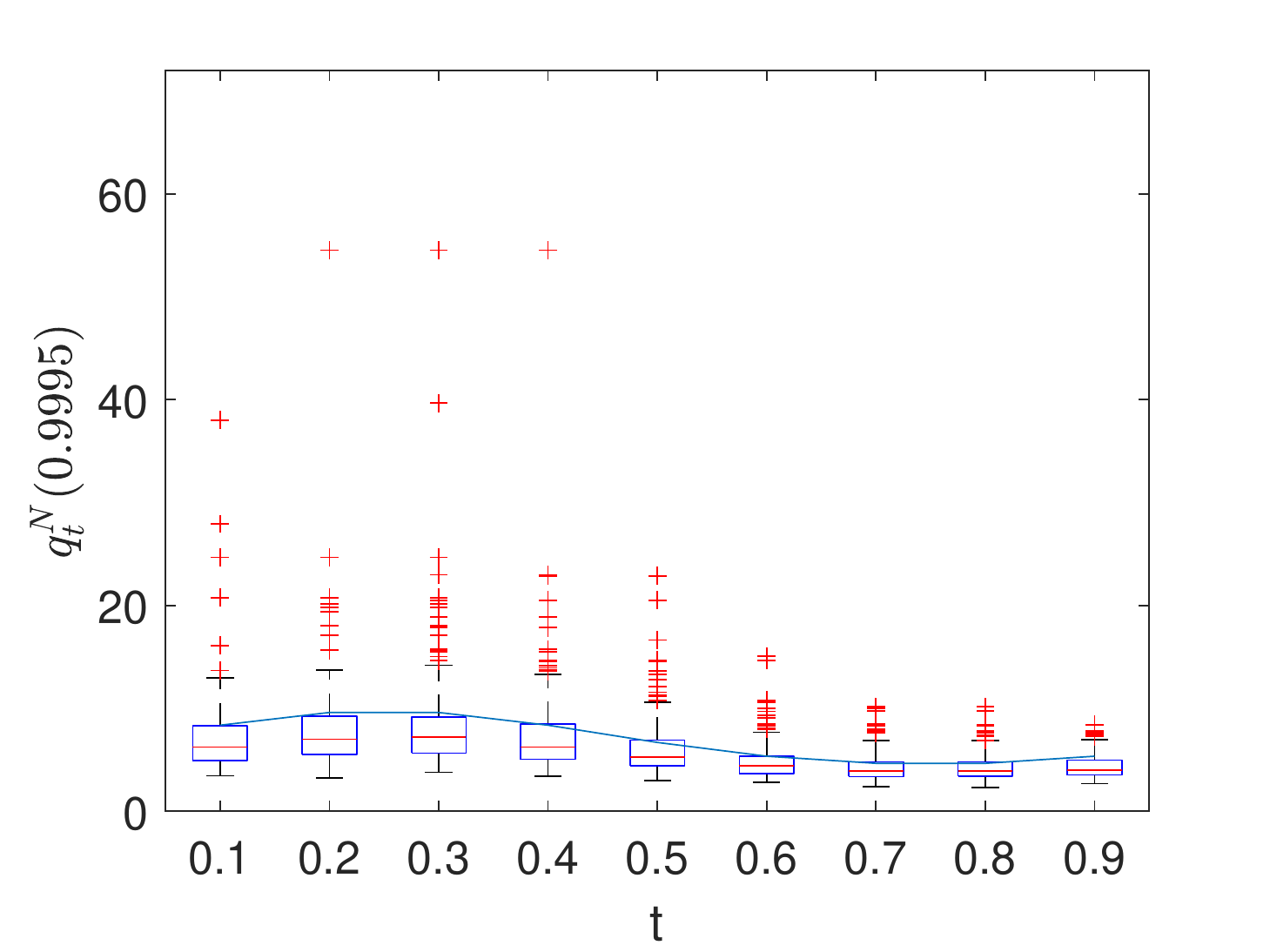}
		
		\caption{The boxplots of the extreme quantiles of $\alpha_N=0.9995$ calculated using our method $\widehat{q}_t^E$ (left) and the naive method (right) from DGP1 with $N=500$ (row 1), 1000 (row 2), 2000 (row 3). }\label{fig:ExtremeQuantileDGP12}
	\end{figure}

\begin{figure}
		\centering
		\includegraphics[width = .45\textwidth]{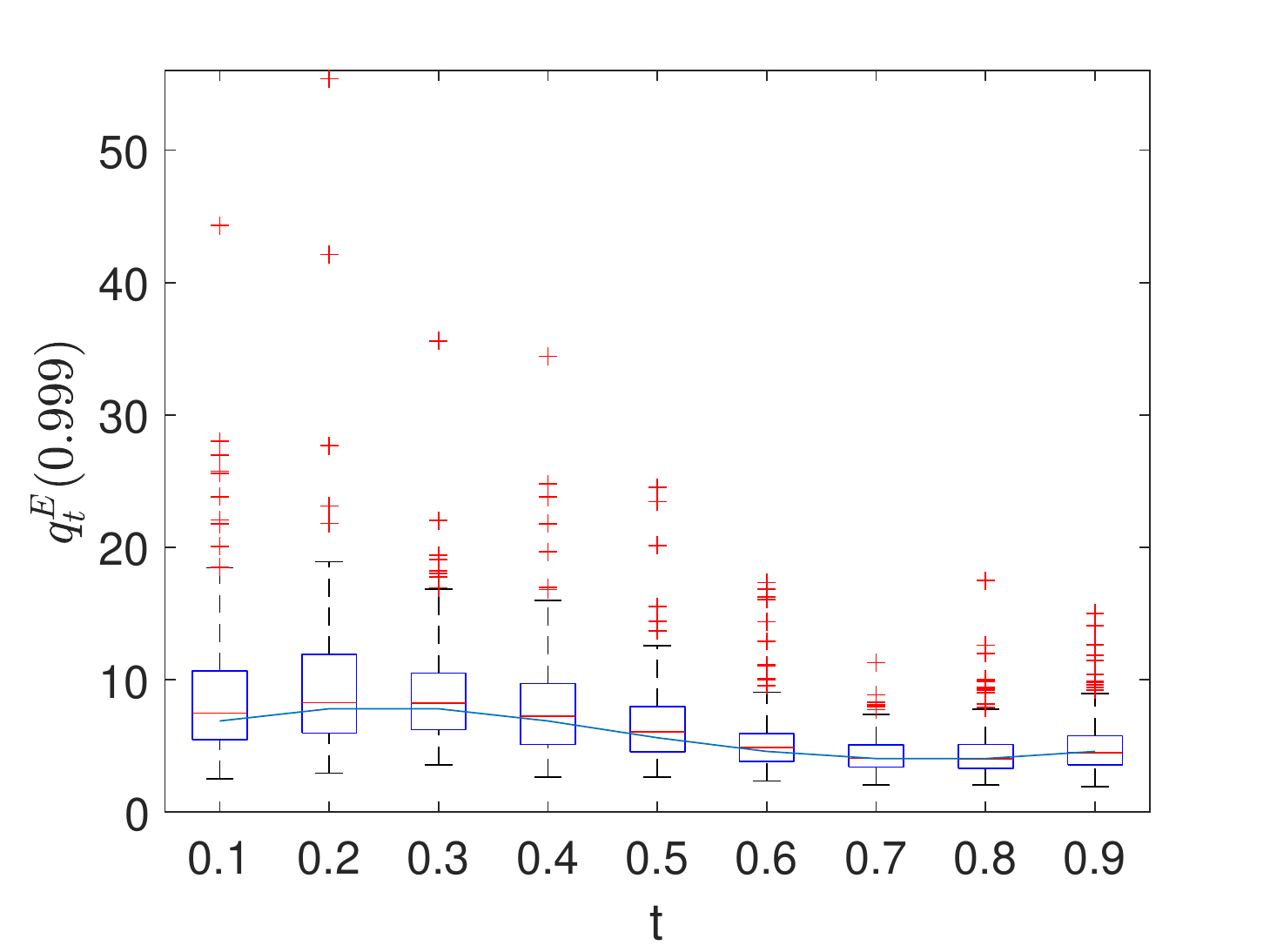}
		\includegraphics[width = .45\textwidth]{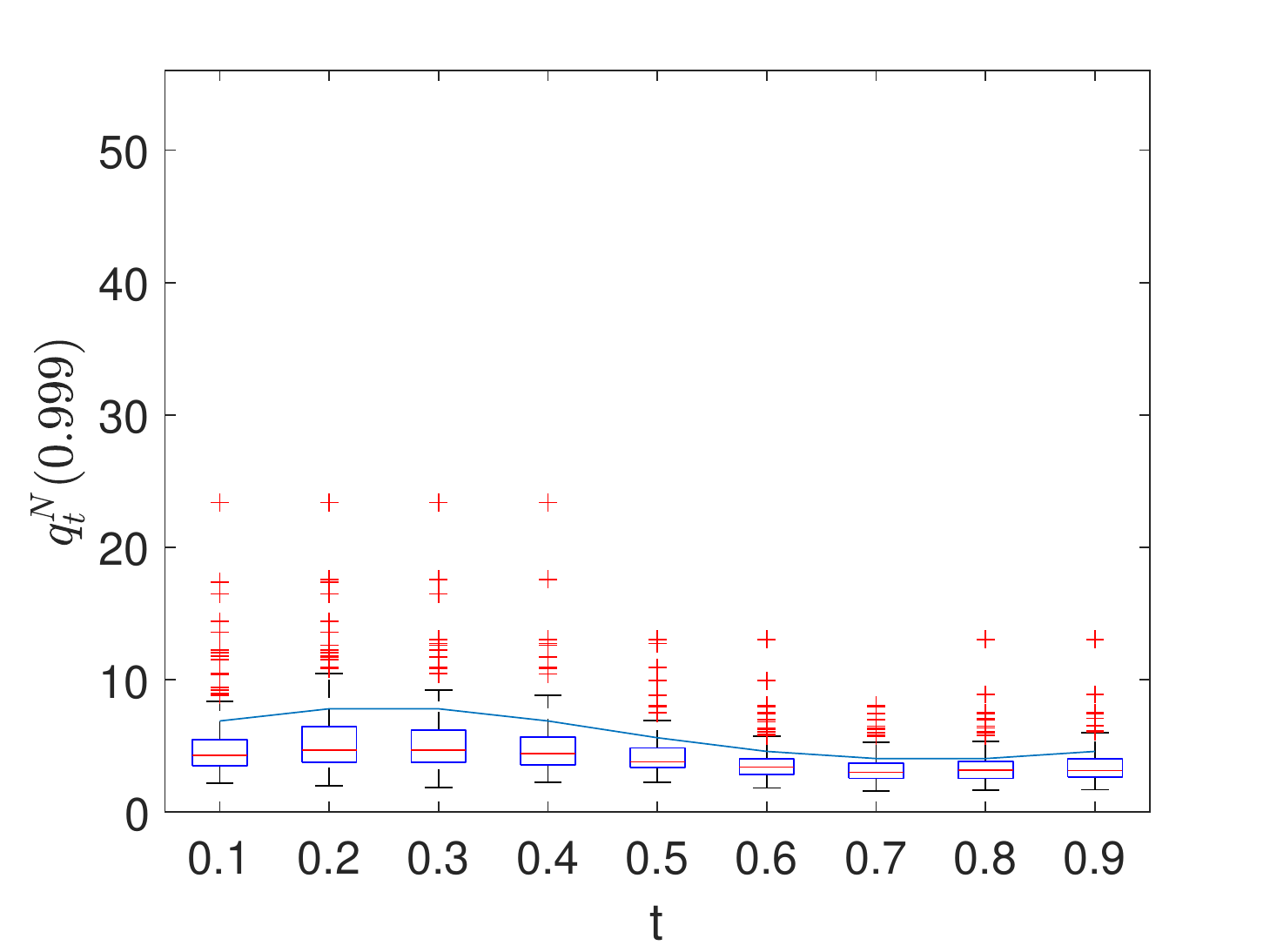}
		
		\includegraphics[width = .45\textwidth]{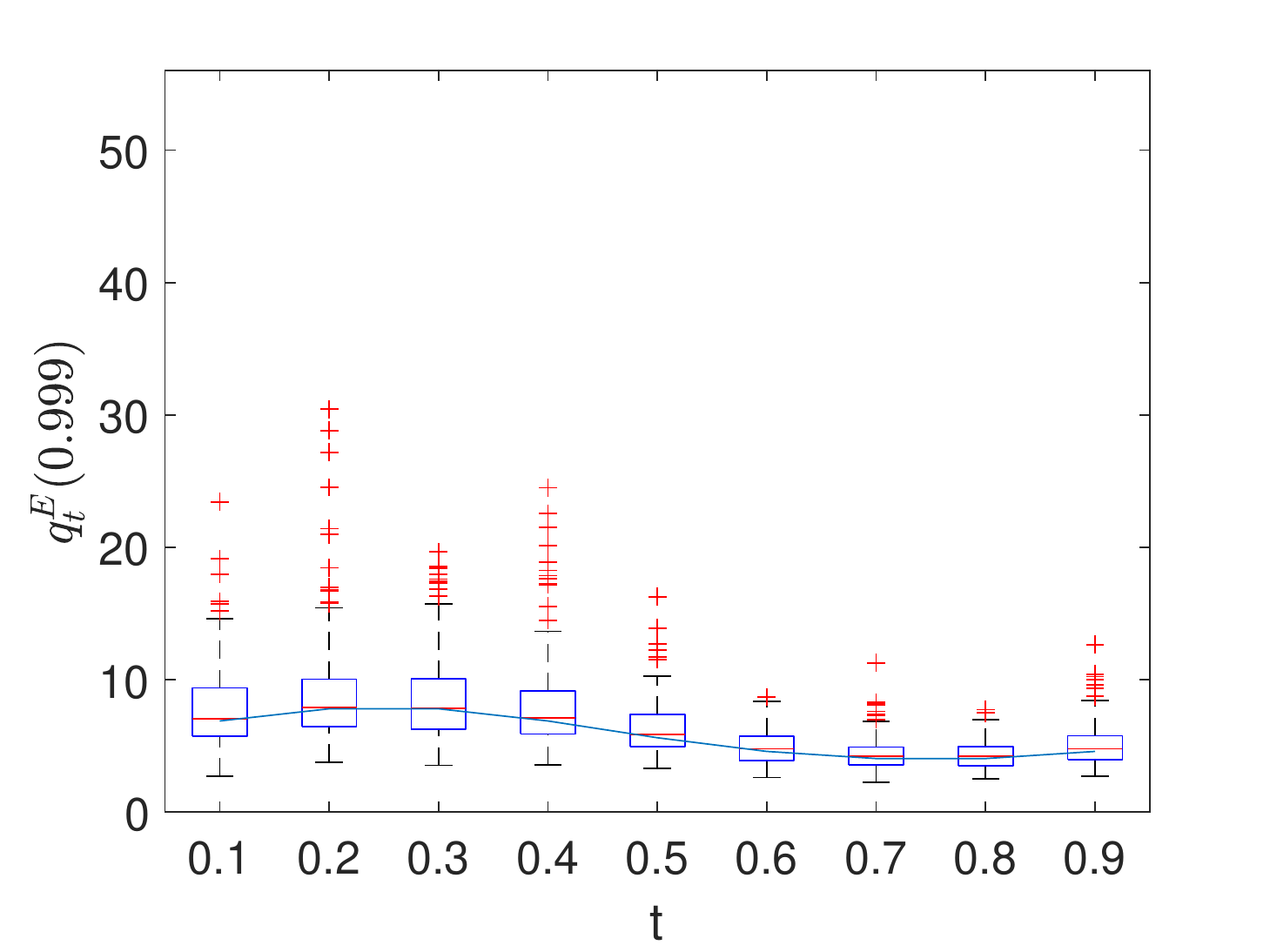}
		\includegraphics[width = .45\textwidth]{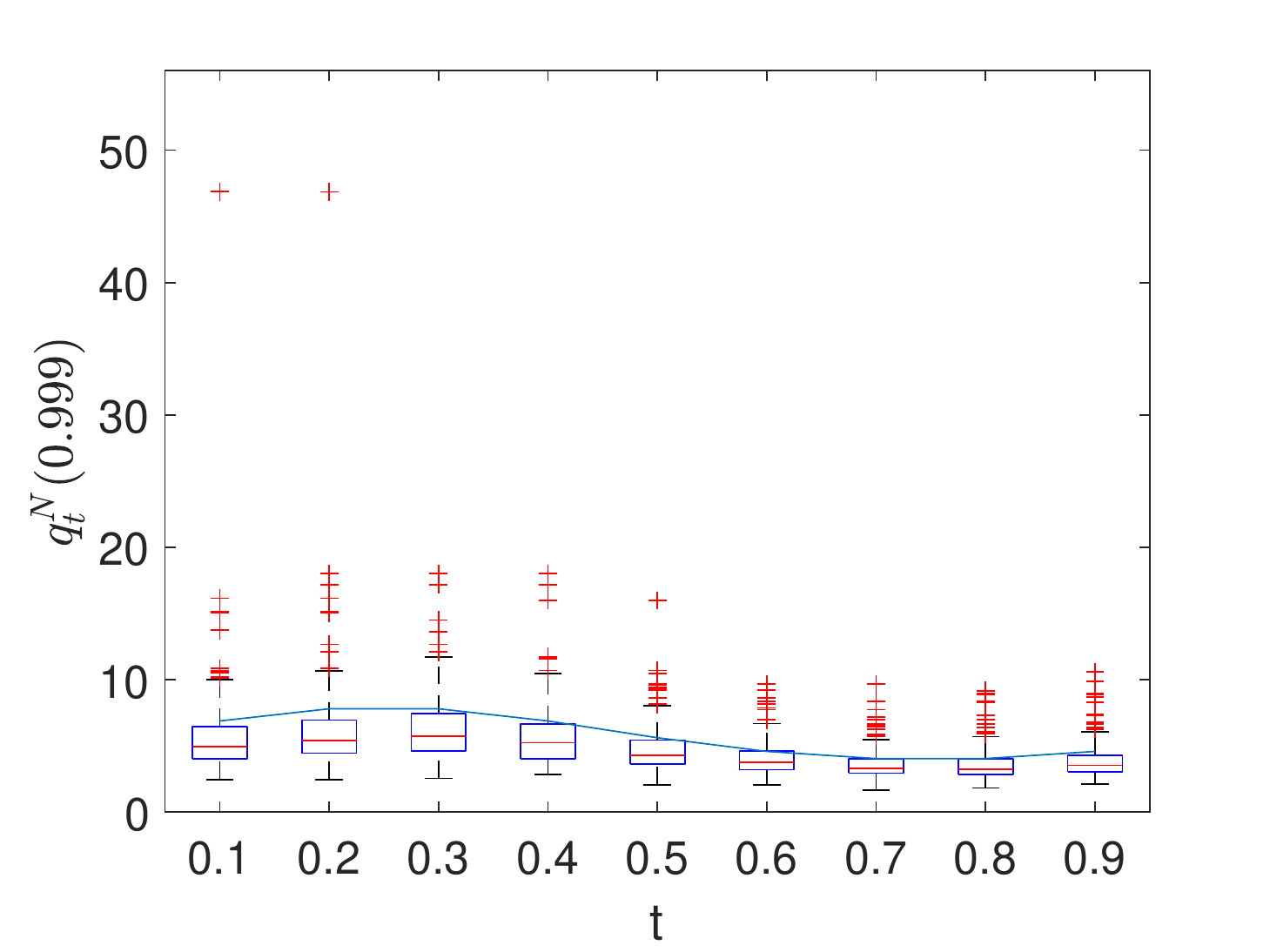}
		
		\includegraphics[width = .45\textwidth]{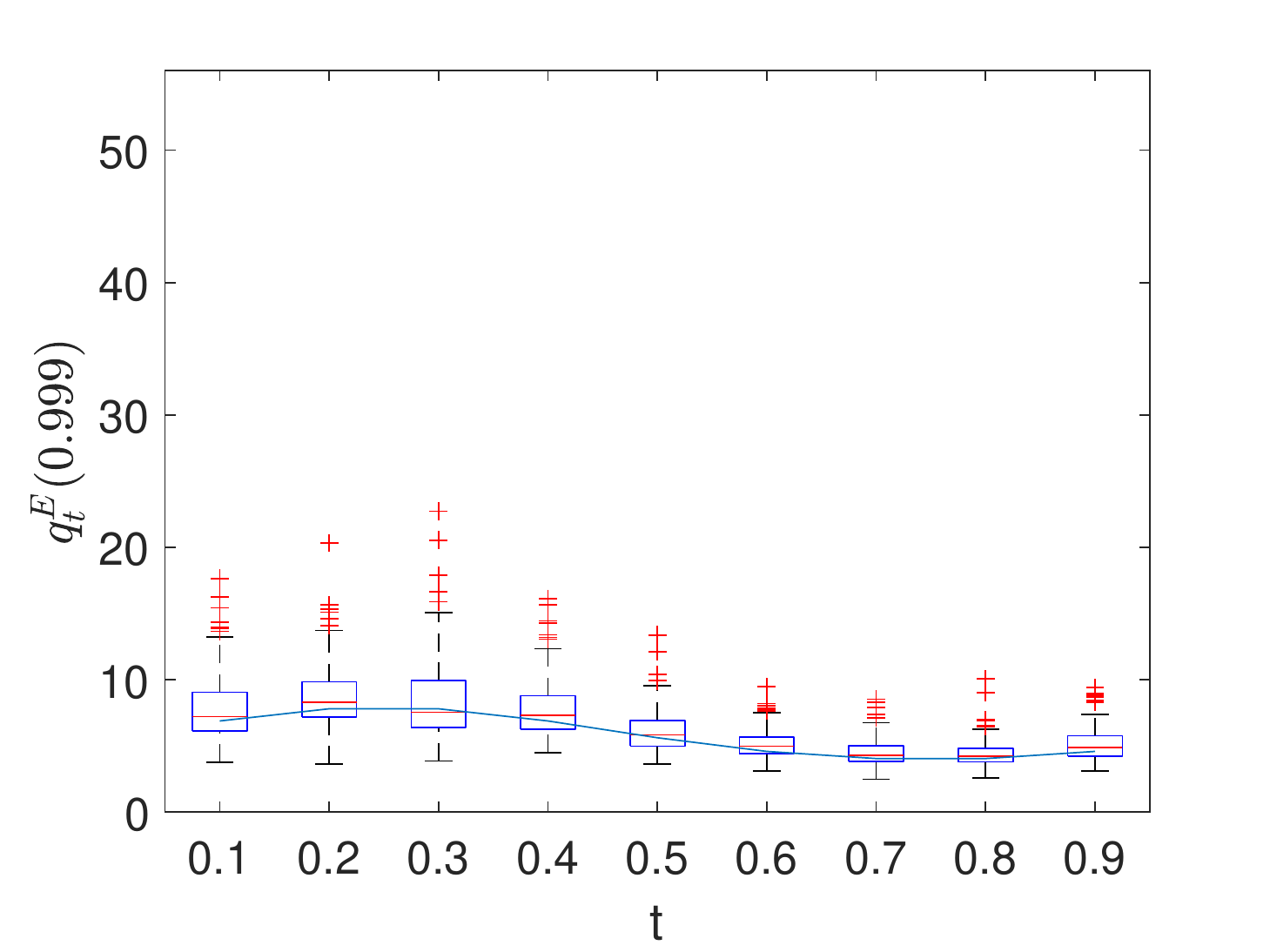}
		\includegraphics[width = .45\textwidth]{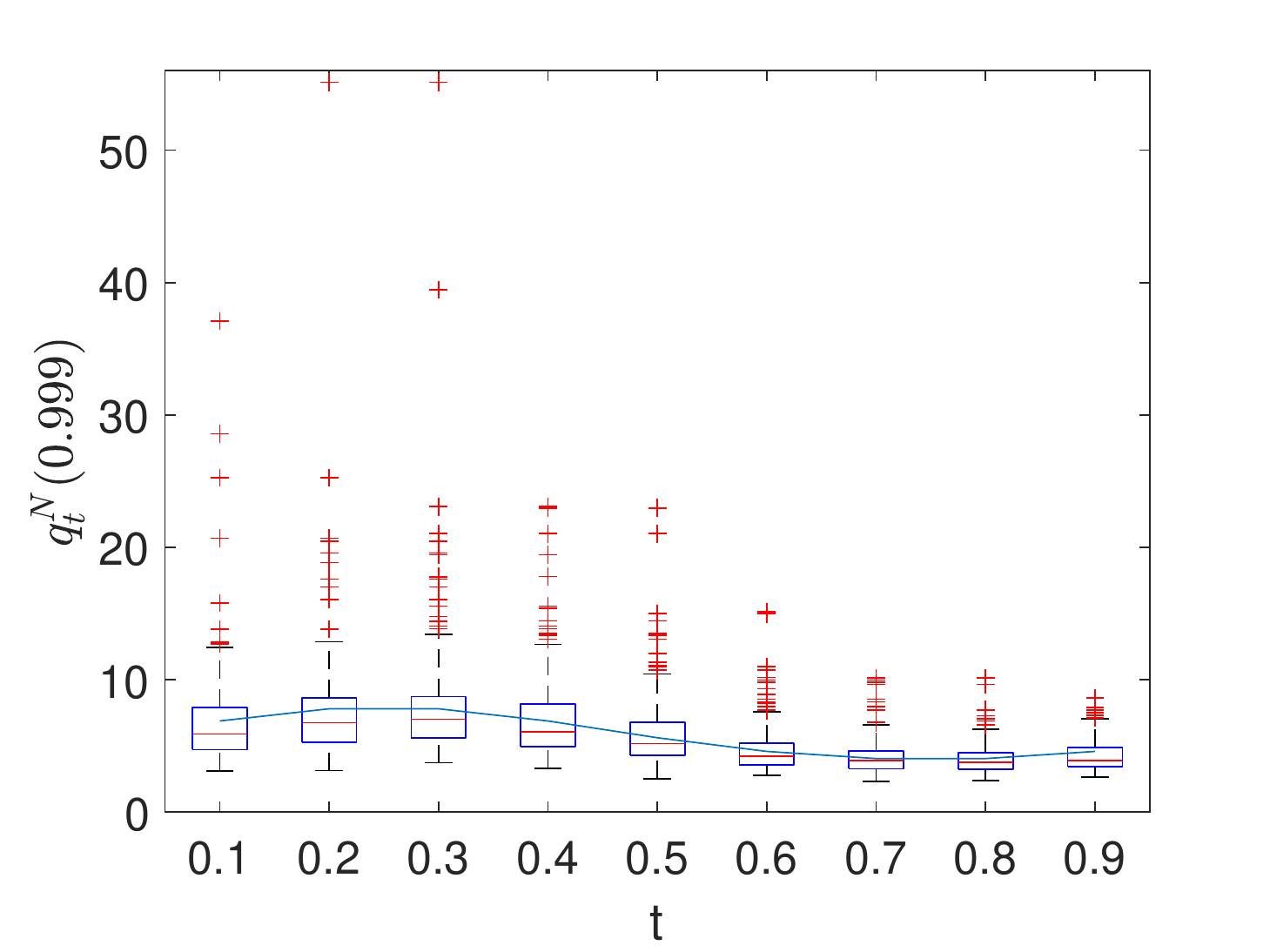}
		
		\caption{The boxplots of the extreme quantiles of $\alpha_N=0.999$ calculated using our method $\widehat{q}_t^E$ (left) and the naive method (right) from DGP2 with $N=500$ (row 1), 1000 (row 2), 2000 (row 3).}\label{fig:ExtremeQuantileDGP2}
	\end{figure}

\begin{figure}
		\centering
		\includegraphics[width = .45\textwidth]{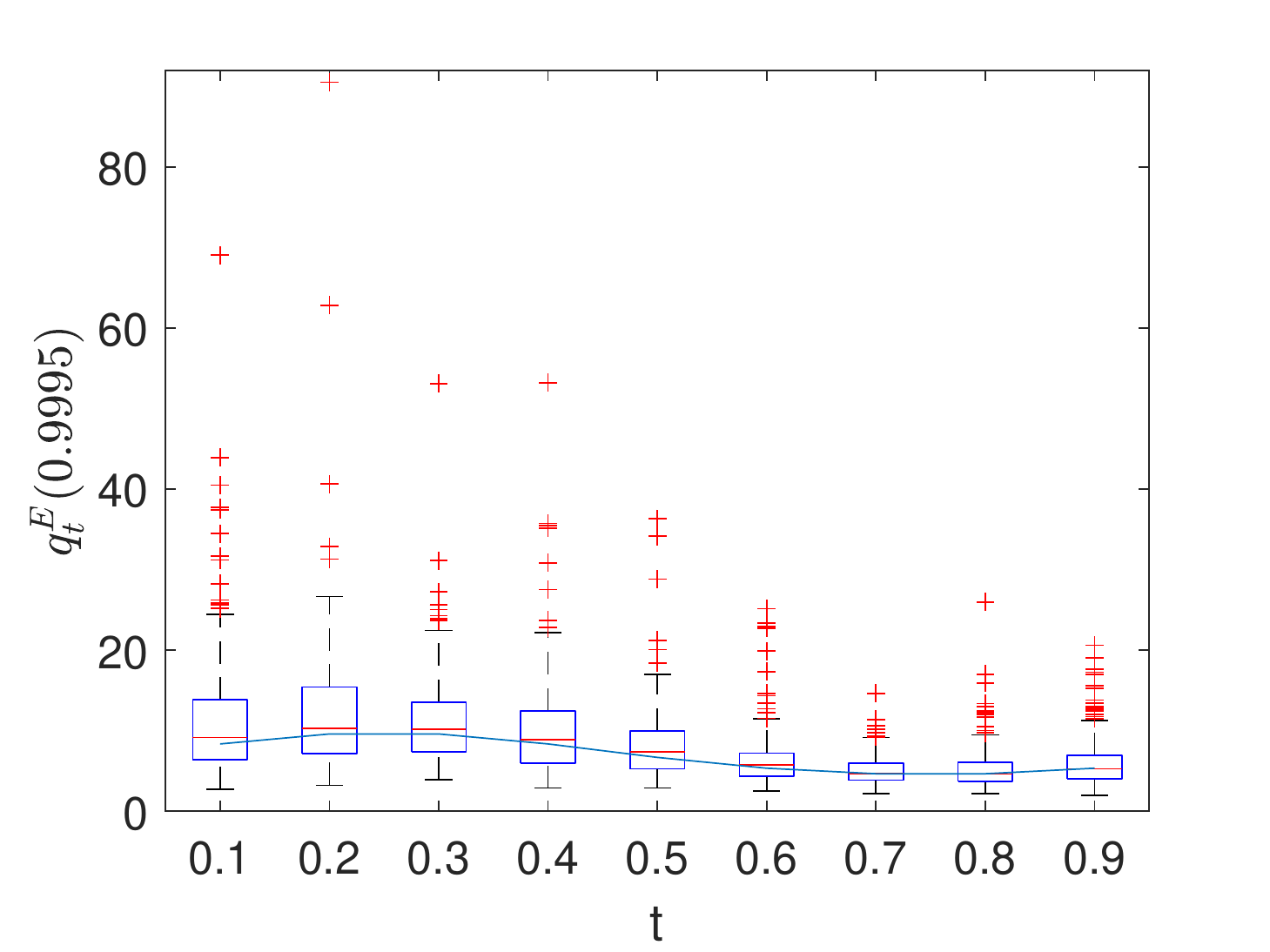}
		\includegraphics[width = .45\textwidth]{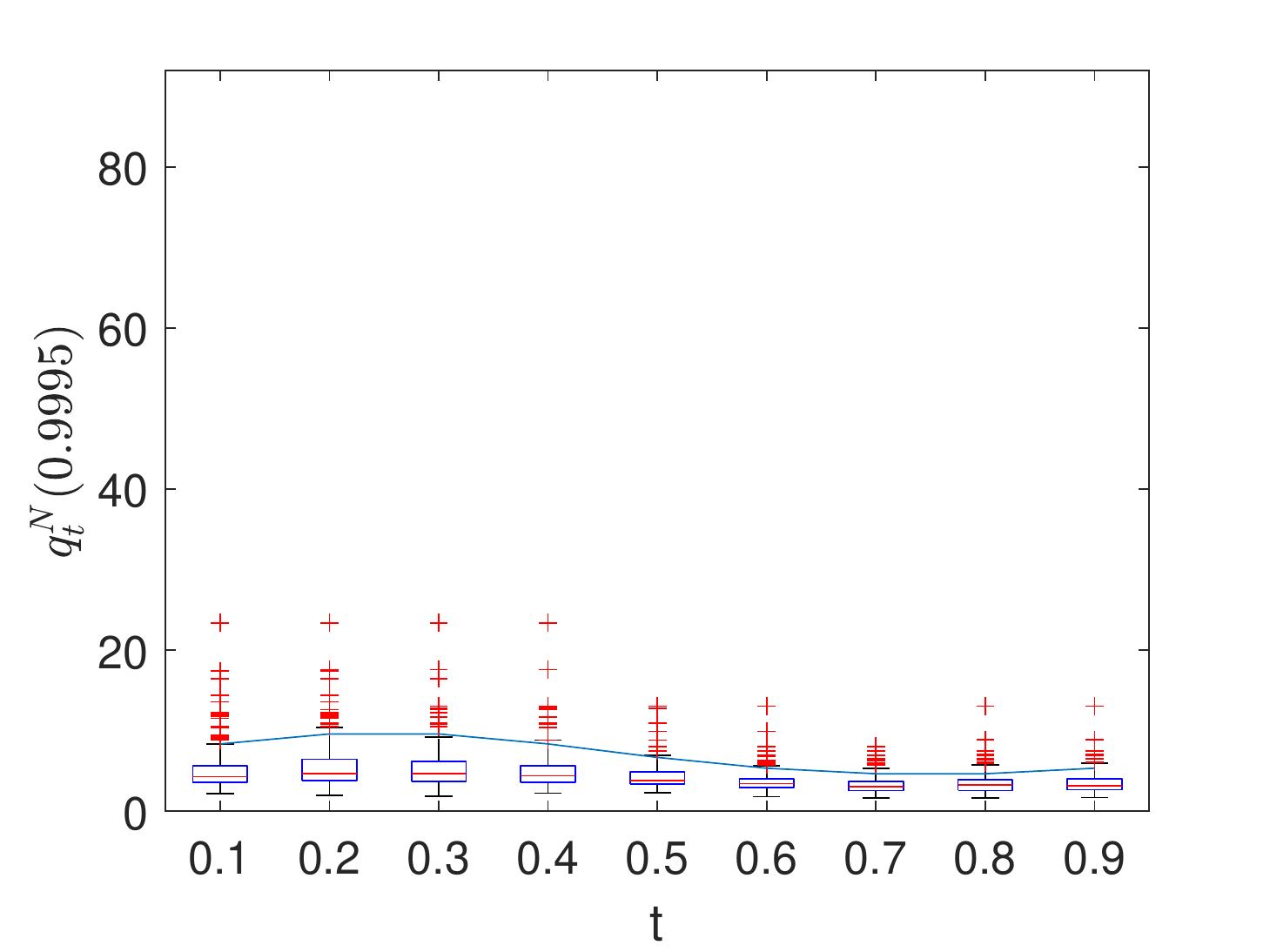}
		
		\includegraphics[width = .45\textwidth]{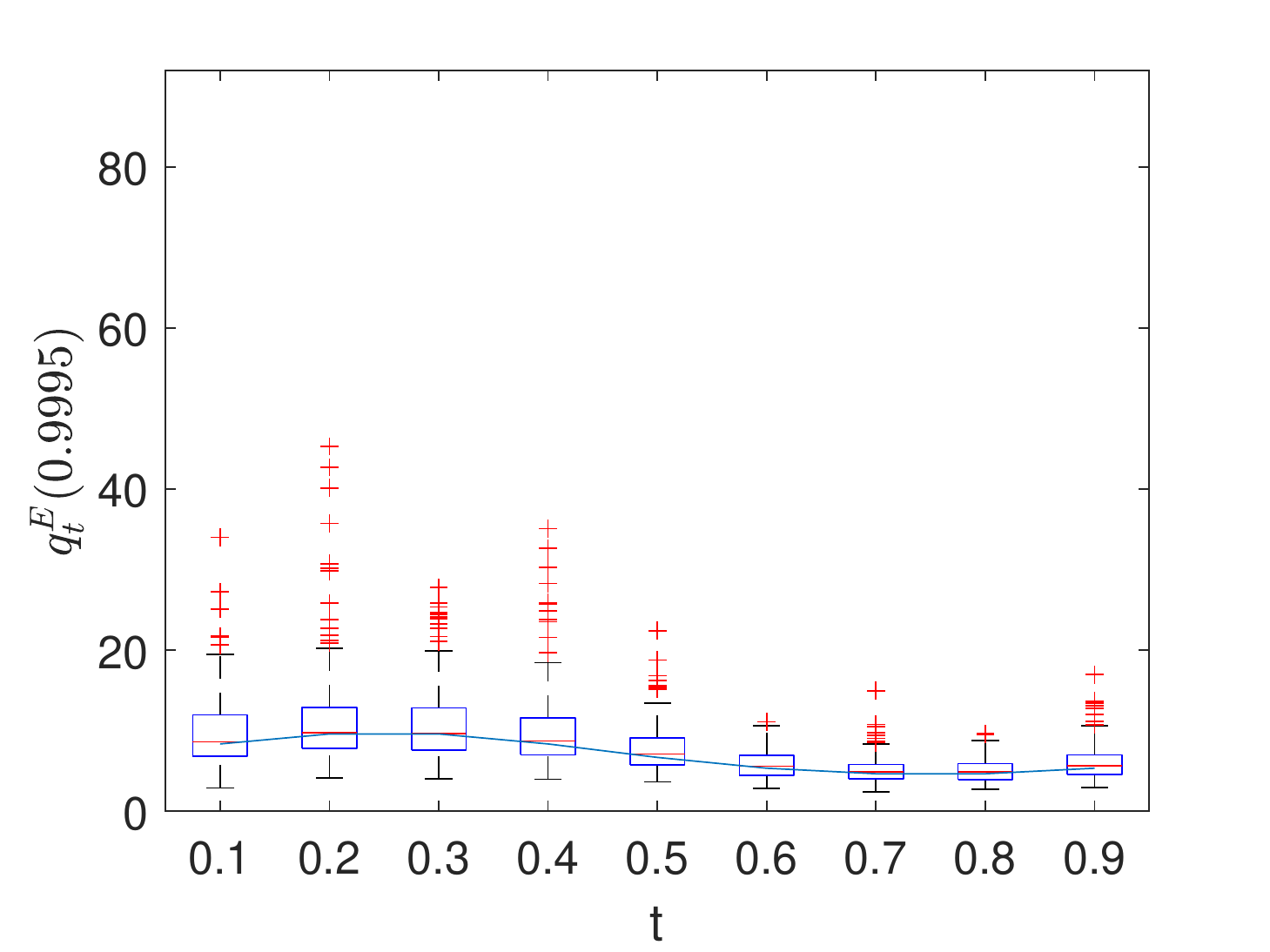}
		\includegraphics[width = .45\textwidth]{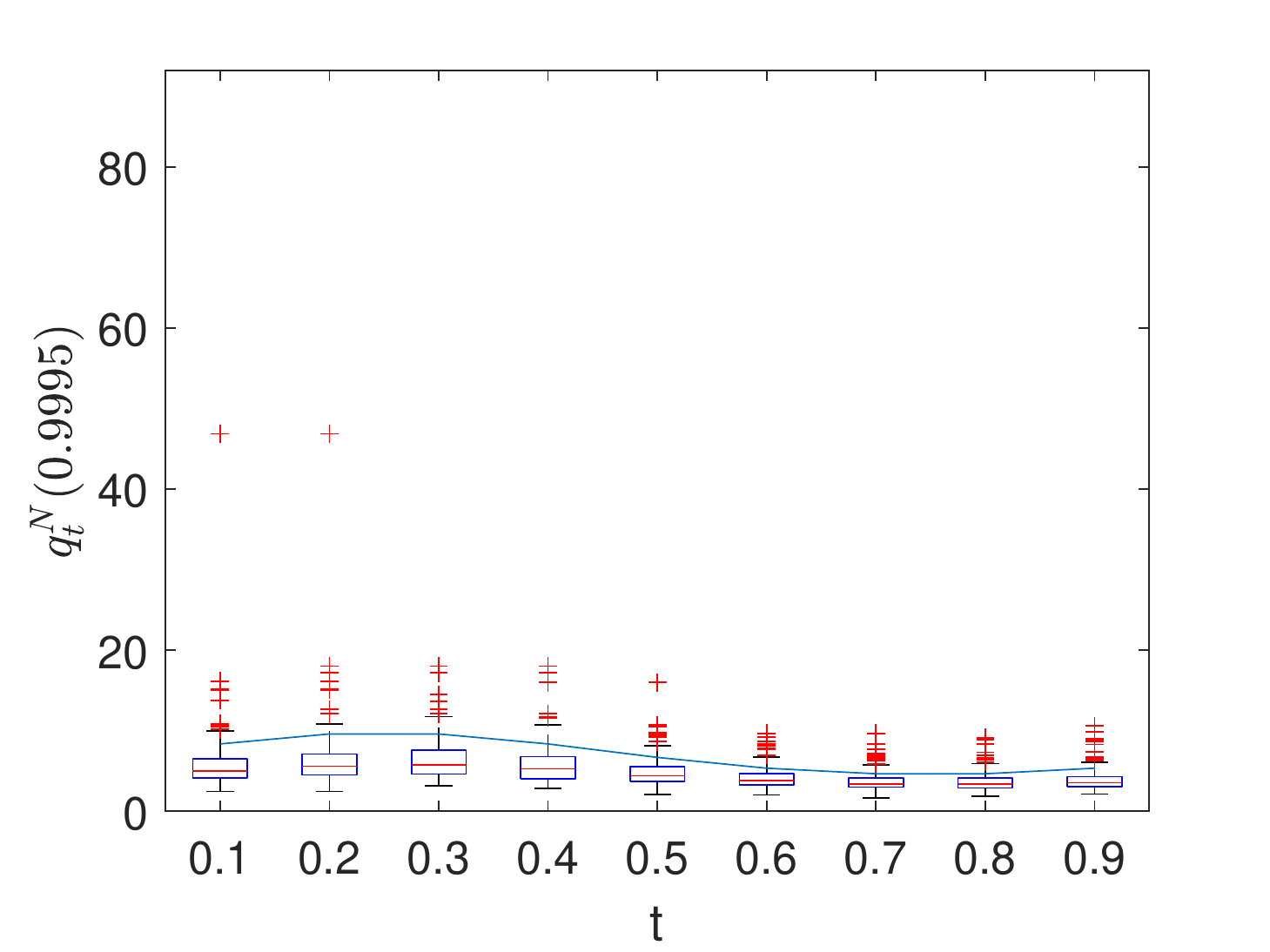}
		
		\includegraphics[width = .45\textwidth]{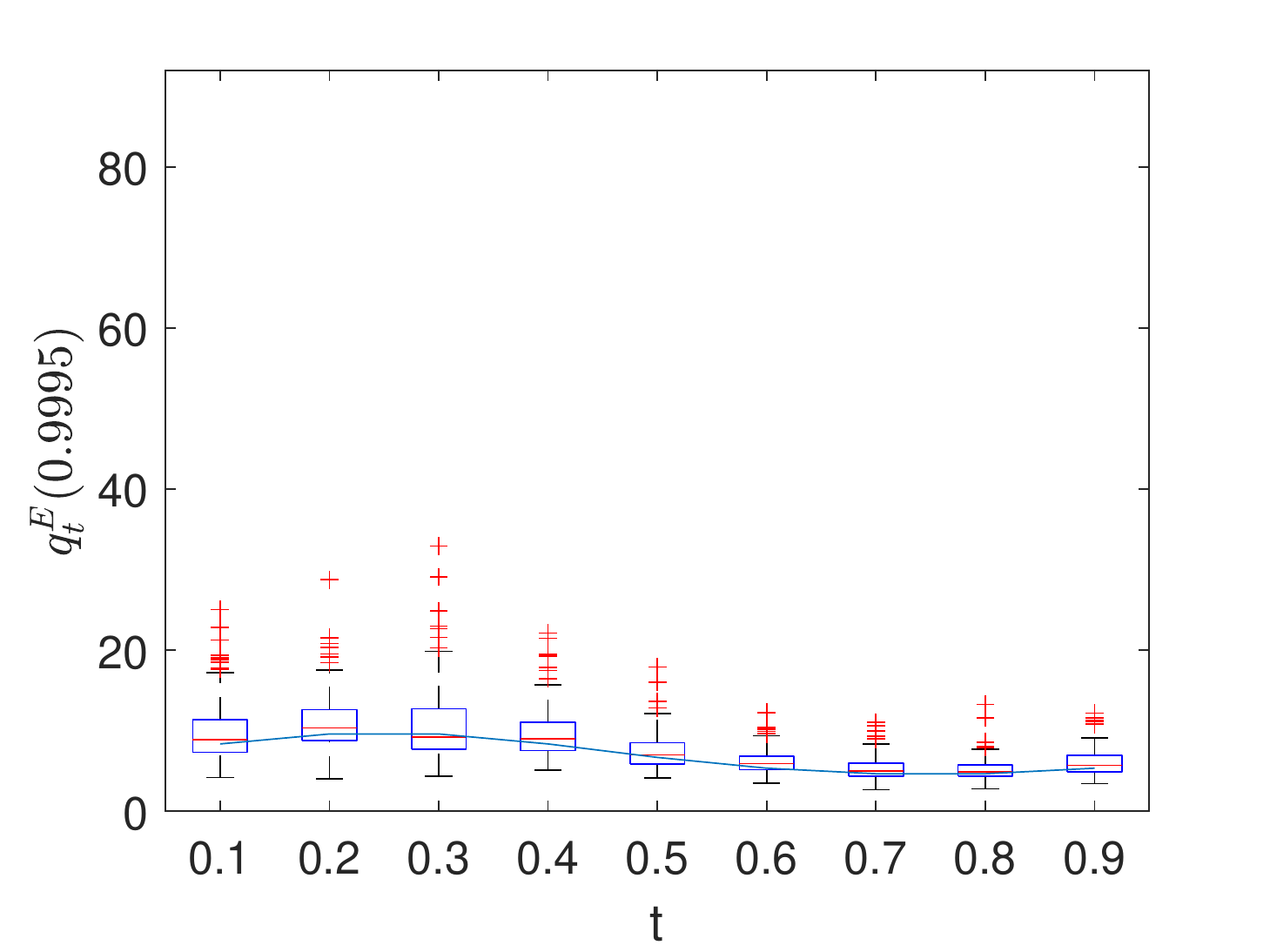}
		\includegraphics[width = .45\textwidth]{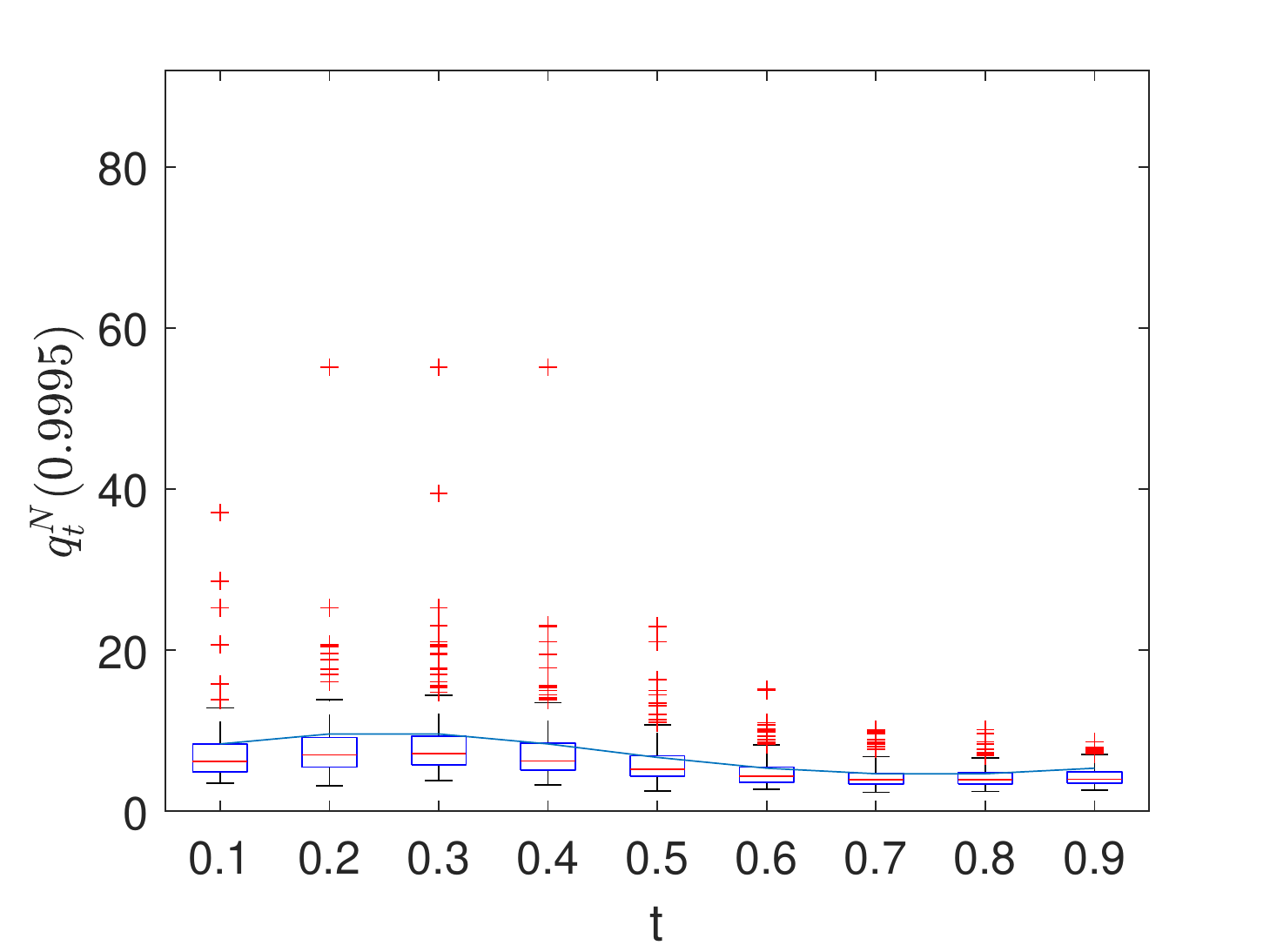}
		
		\caption{The boxplots of the extreme quantiles of $\alpha_N=0.9995$ calculated using our method $\widehat{q}_t^E$ (left) and the naive method (right) from DGP2 with $N=500$ (row 1), 1000 (row 2), 2000 (row 3).}\label{fig:ExtremeQuantileDGP22}
	\end{figure}

	From our simulation studies, we found that $6\leq J \leq 10$ gives similarly good results. Figures~\ref{fig:ExtremeQuantileDGP1} to \ref{fig:ExtremeQuantileDGP22} show the boxplots of our $\widehat{q}_t^E$, for $t=0.1,0.2,\ldots,0.9$, $\alpha_N=0.999$ and $0.9995$ with $J=8$ and the naive extrapolated quantiles from the 200 samples of DGP1 and DGP2. We can see that the naive estimator is biased and does not show any convergence, while our proposed estimator remains unbiased and converges to the true extreme quantile as $N$ increases.
	
\begin{figure}
		\centering
		\includegraphics[width = .43\textwidth]{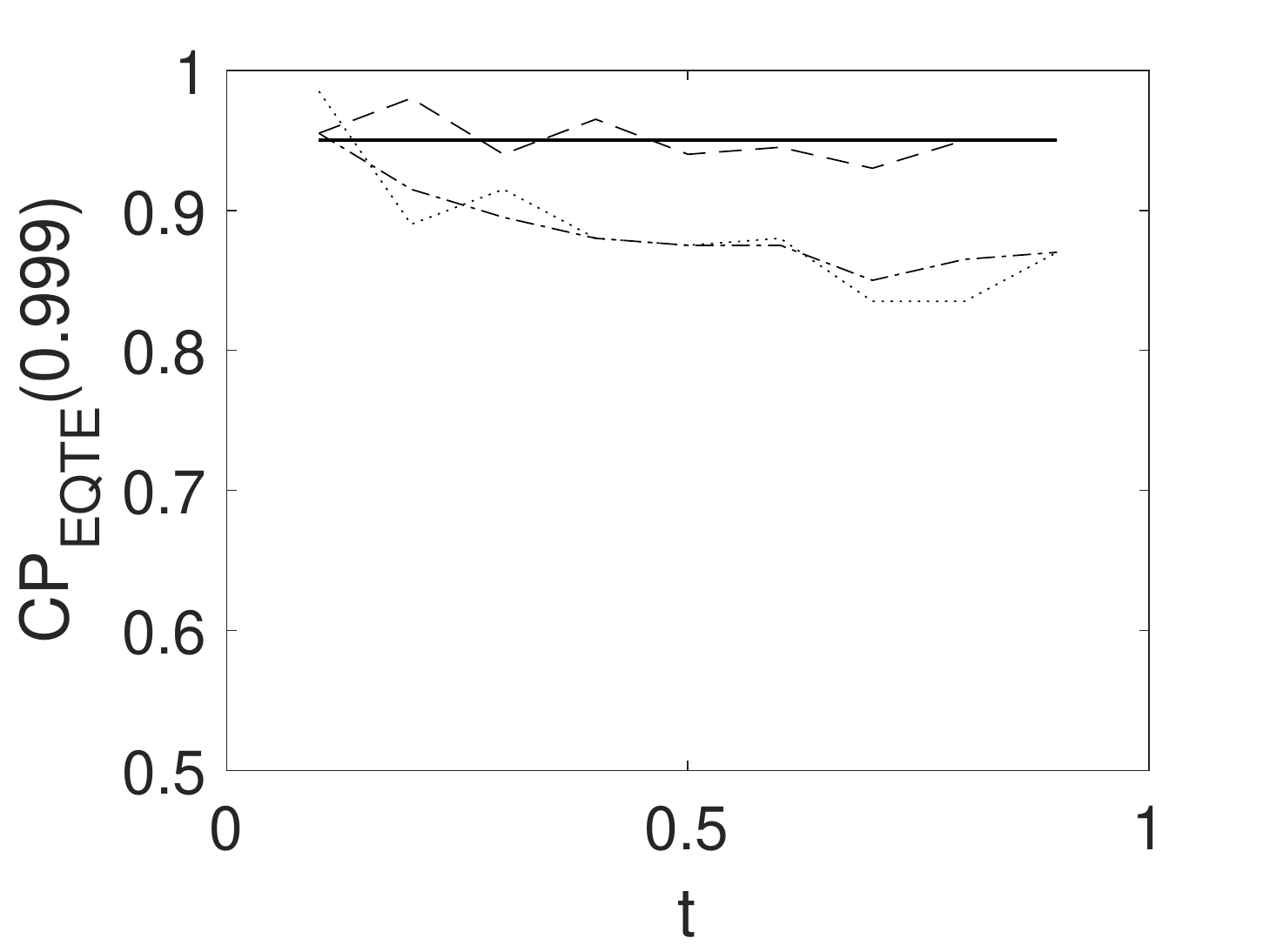}
		\includegraphics[width = .43\textwidth]{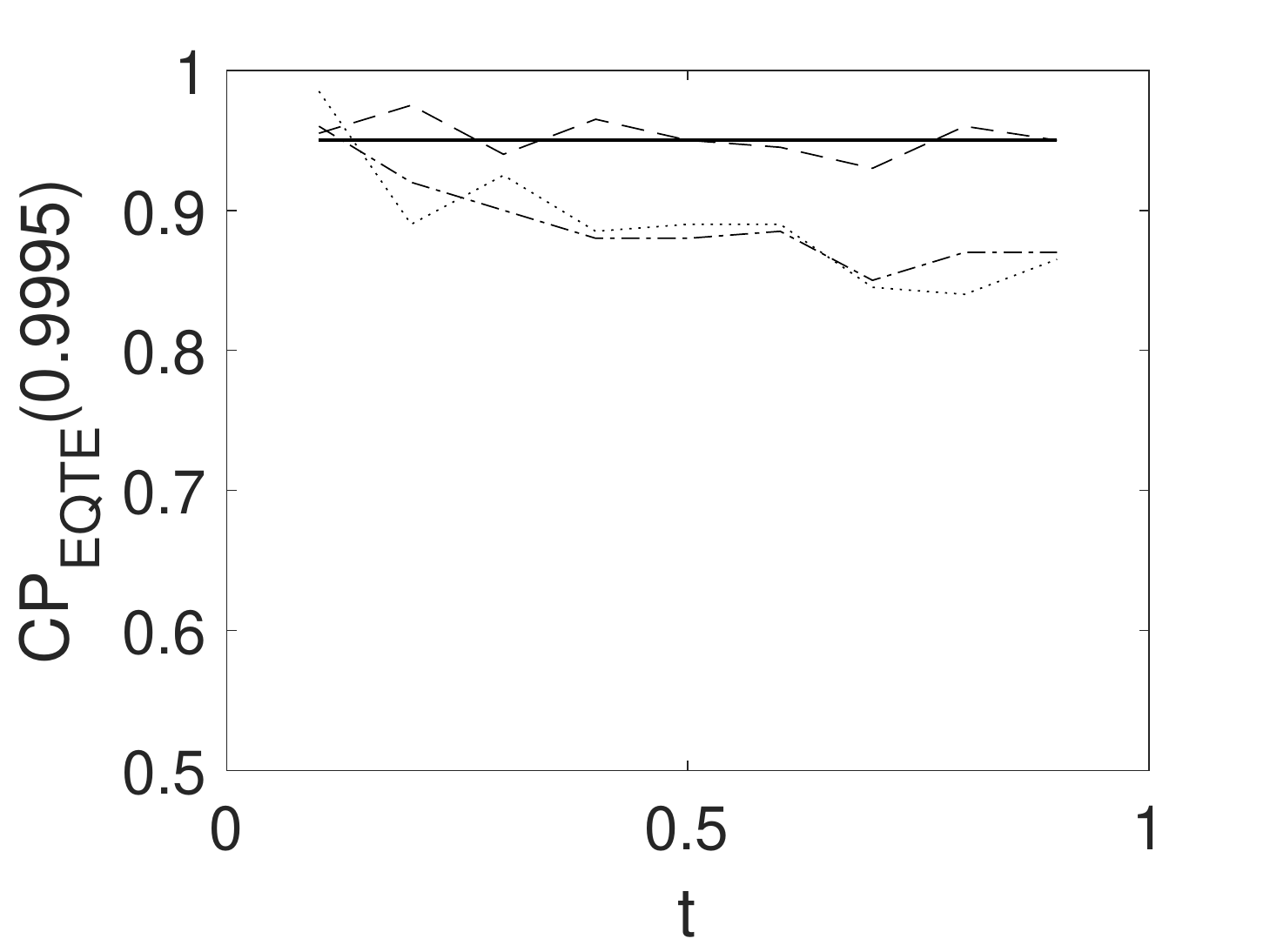}
		
		\includegraphics[width = .43\textwidth]{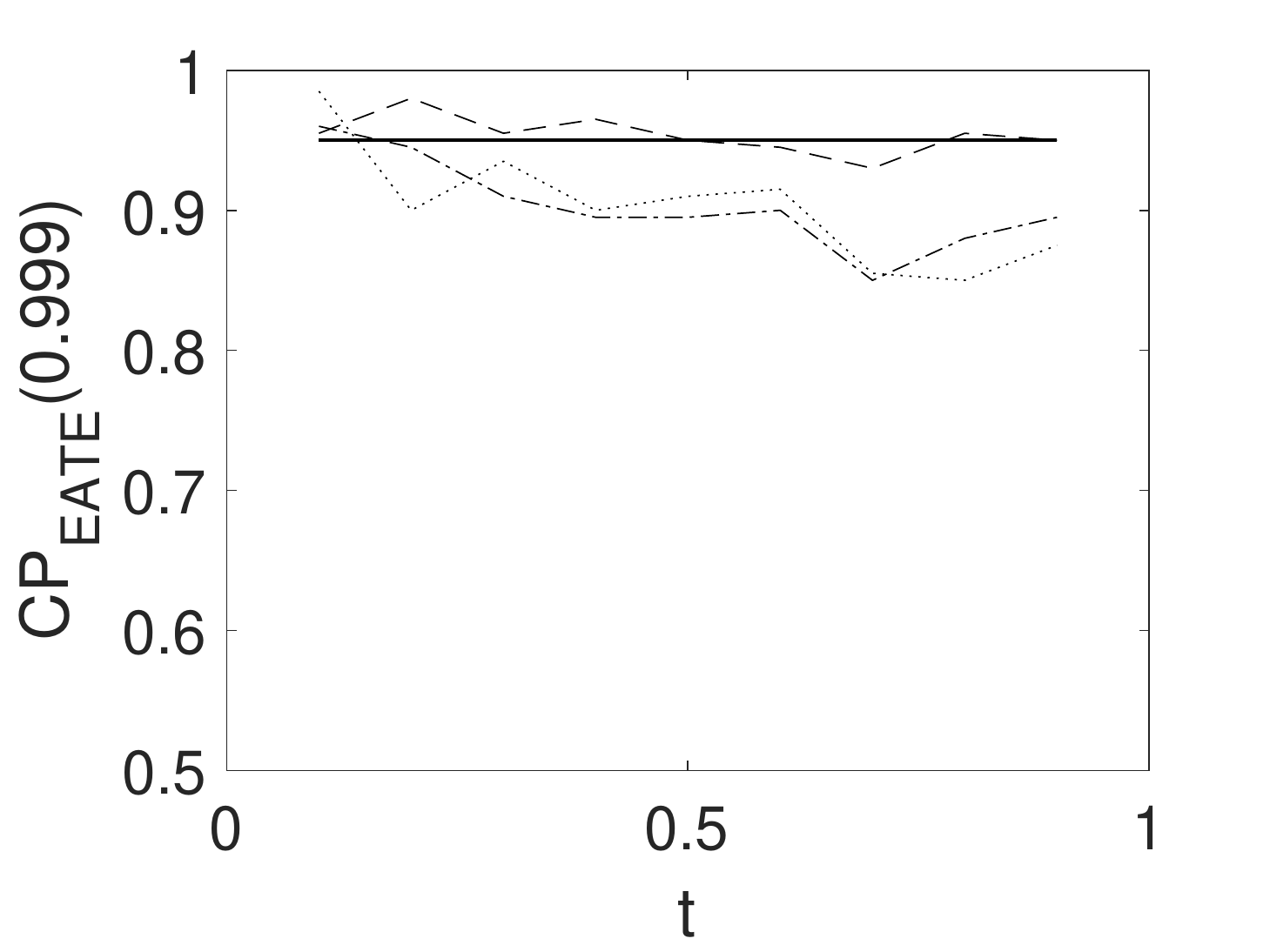}
		\includegraphics[width = .43\textwidth]{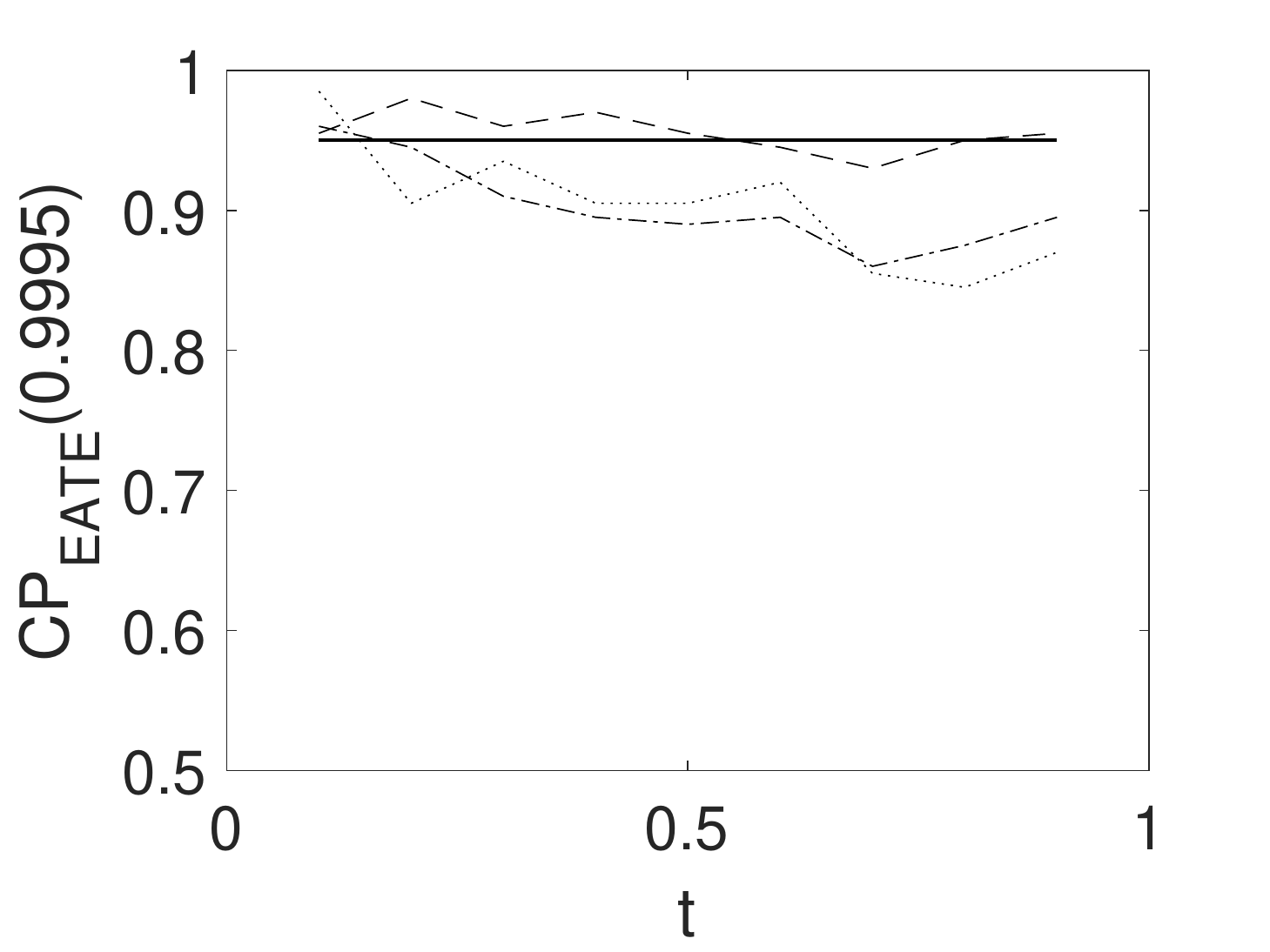}
		
		\caption{The empirical coverage probability of the 95\% confidence intervals for EQTE (rows 1) and EATE (rows 2) of DGP1 with $\alpha_N=0.999$ (left) and $\alpha_N=0.9995$ (right), sample size $N=500$ (dotted line), $N=1000$ (dash-dotted line) and $N=2000$ (dashed line). }\label{fig:CoverageProbDGP1}
	\end{figure}

\begin{figure}
		\centering
		\includegraphics[width = .45\textwidth]{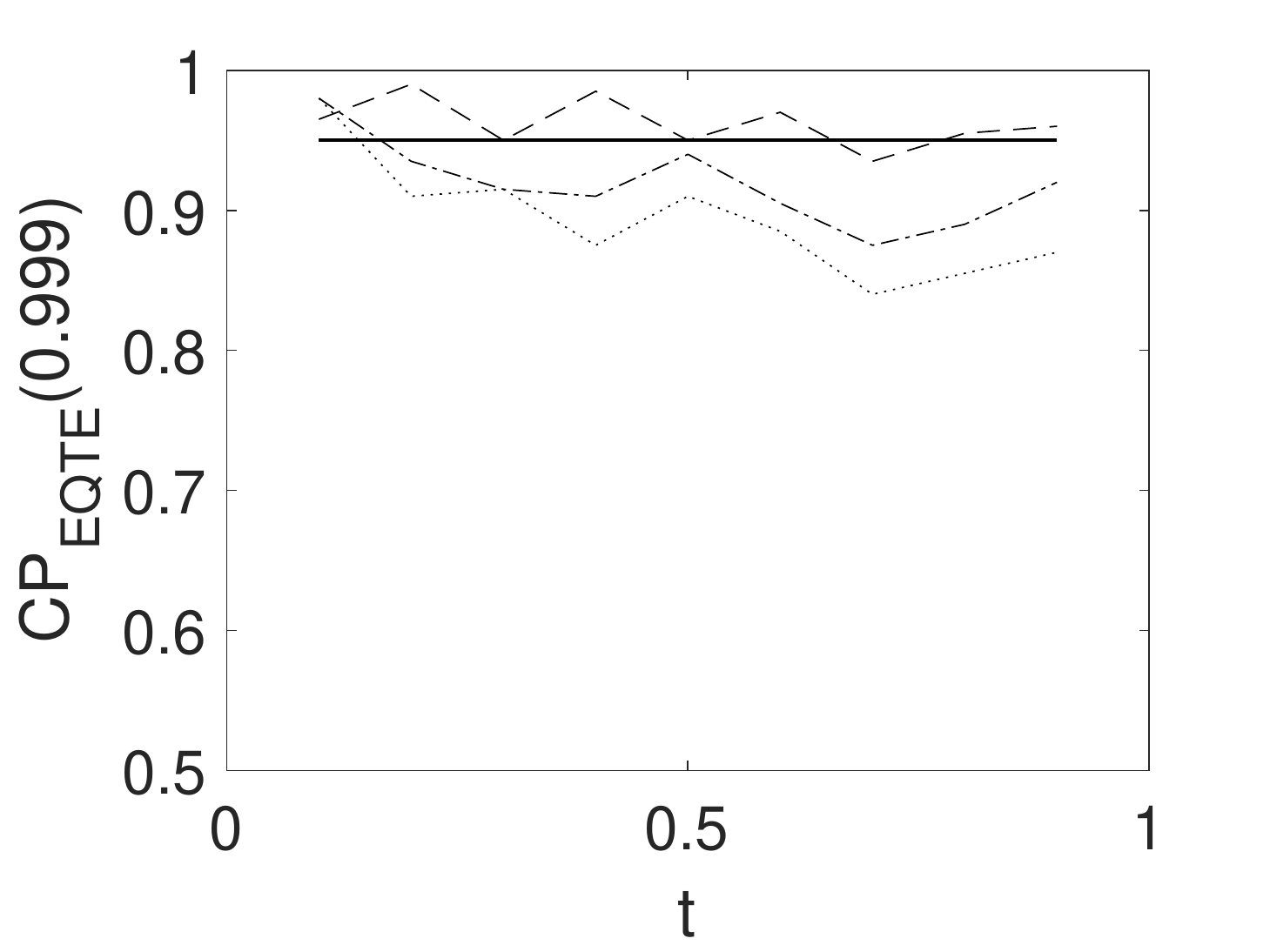}
		\includegraphics[width = .45\textwidth]{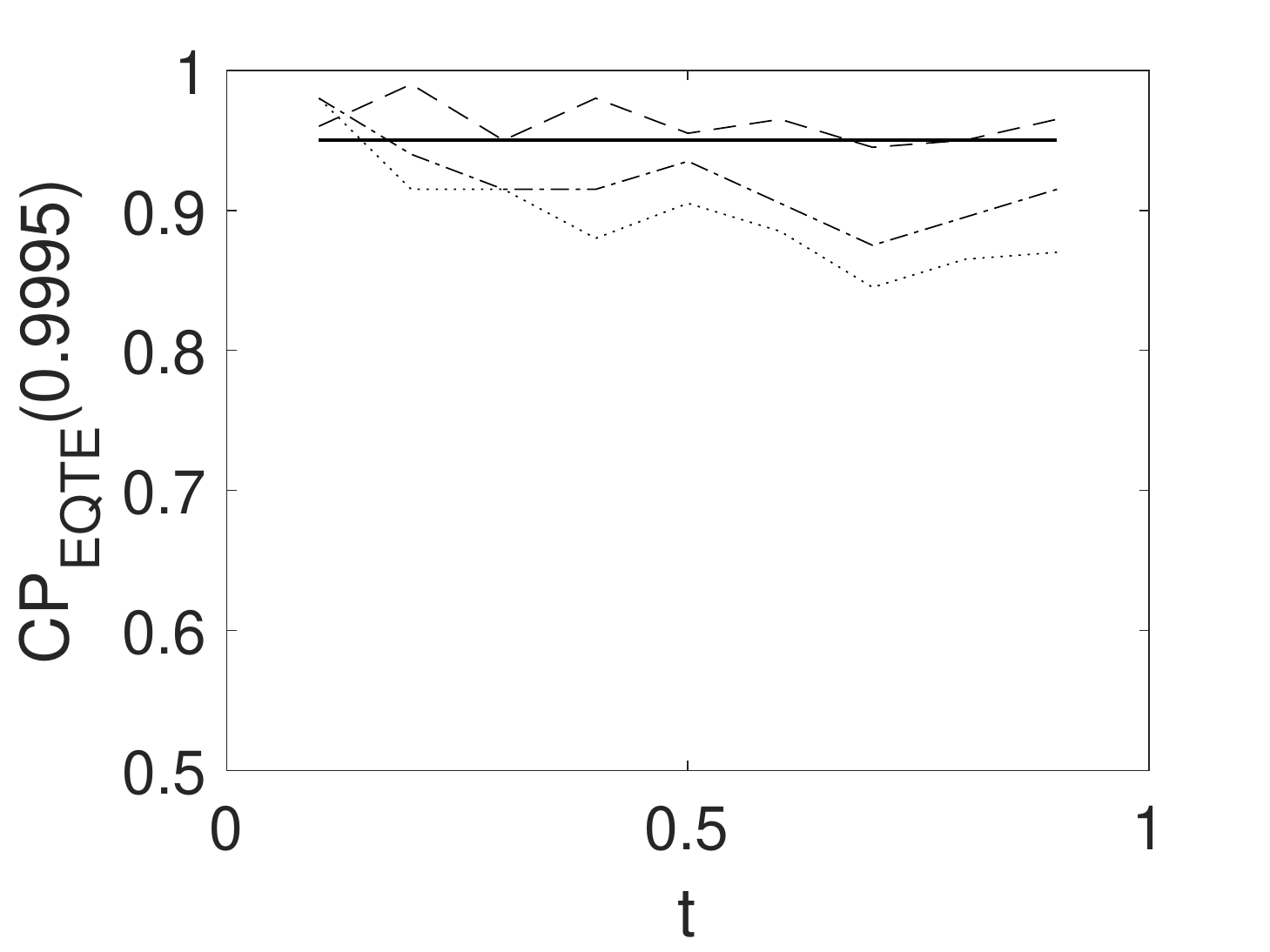}
		
		\includegraphics[width = .45\textwidth]{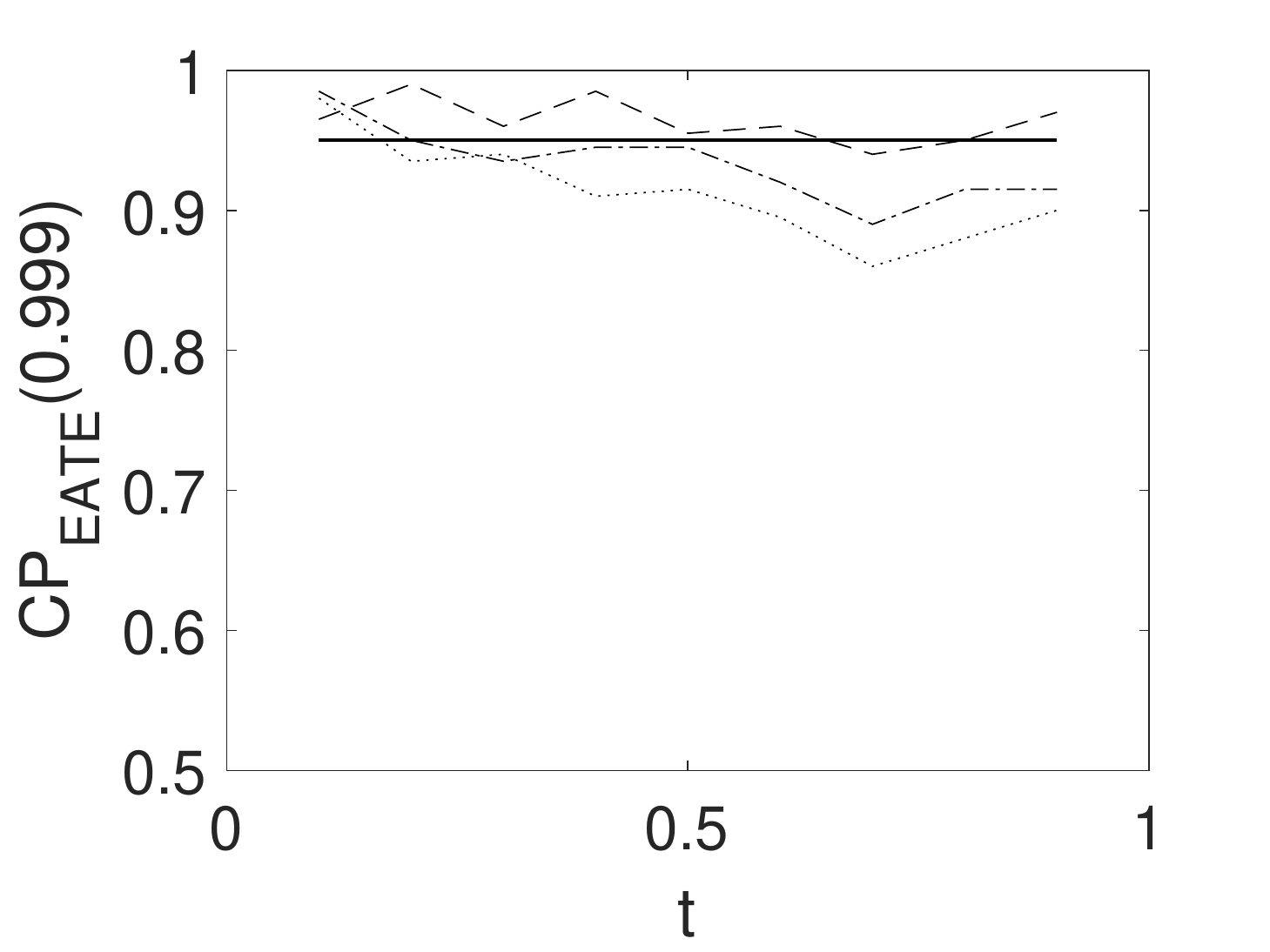}
		\includegraphics[width = .45\textwidth]{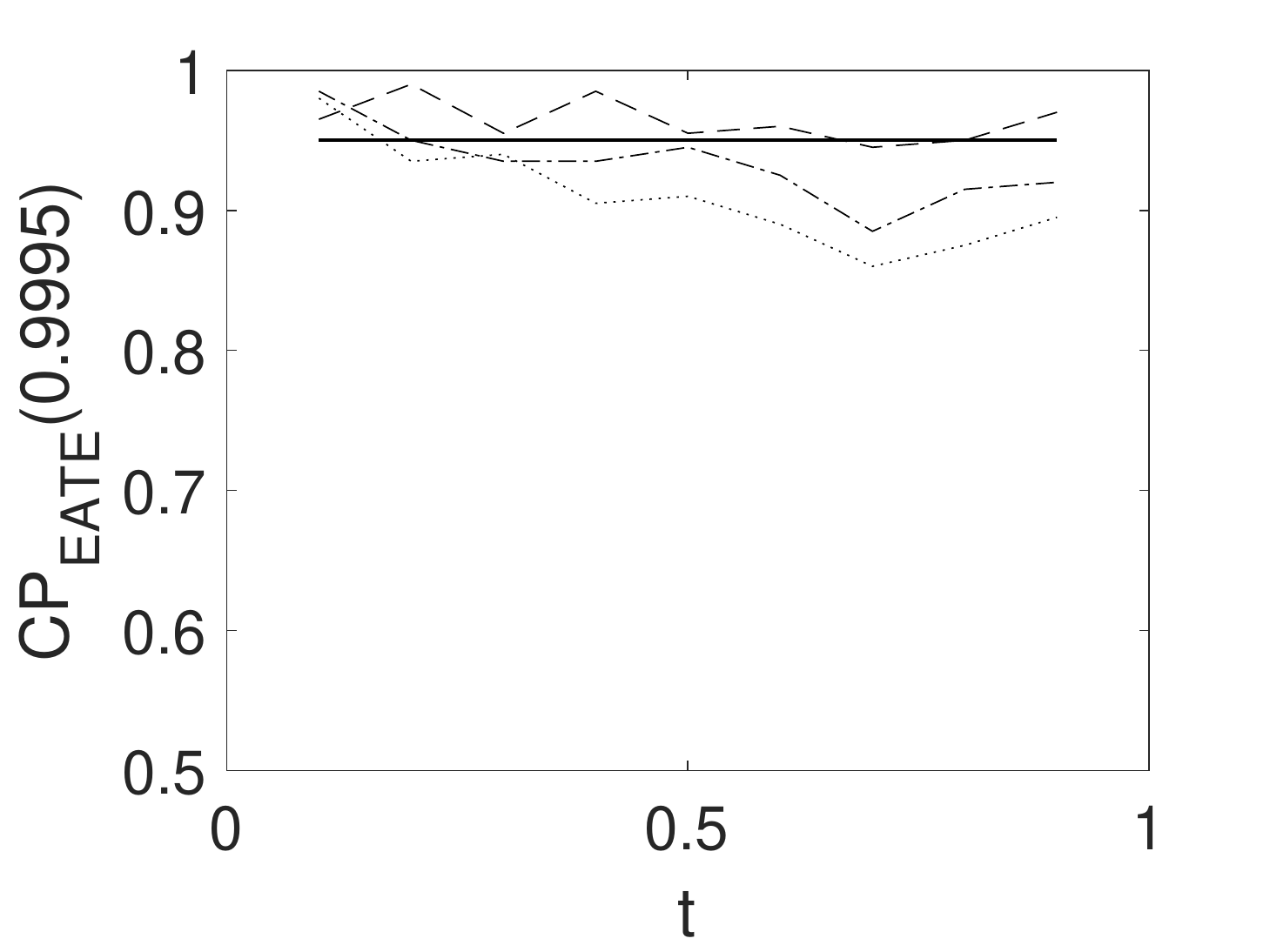}
		
		\caption{The empirical coverage probability of the 95\% confidence intervals for EQTE (row 1) and EATE (row 2) of DGP2 with $\alpha_N=0.999$ (left) and $\alpha_N=0.9995$ (right), sample size $N=500$ (dotted line), $N=1000$ (dash-dotted line) and $N=2000$ (dashed line). }\label{fig:CoverageProbDGP2}
	\end{figure}

	Figures~\ref{fig:CoverageProbDGP1} and \ref{fig:CoverageProbDGP2} depict the the empirical coverage probability of our proposed 95\% confidence intervals for $\text{EQTE}_{0,t}(\alpha_N)$ and $\text{EATE}_{0,t}(\alpha_N)$, for $t=0.1,\ldots,0.9$ and $\alpha_N=0.999$ and $0.9995$. For all the circumstances, our method gives reasonable coverage and converges to the nominal confidence level as sample size $N$ increases.
	
	\subsection{Real Data Example}\label{sec:real_data}
	
	In this section, we apply our method to the U.S. presidential campaign data used in \cite{huang2021unified}. The data have been analyzed several times in the treatment effect literature \citep{Urban_Niebler_2014,Fong_Hazlett_Imai_2018, Ai_Linton_Motegi_Zhang_cts_treat}, where the interest was to explore the casual relationship between advertising and campaign contributions. The treatment of interest is the number of political advertisements aired in each zip code from non-competitive states, which ranges from 0 to 22379 across $N=16265$ zip codes.
	
	\begin{figure}[t]
		\centering
		\includegraphics[width = .45\textwidth]{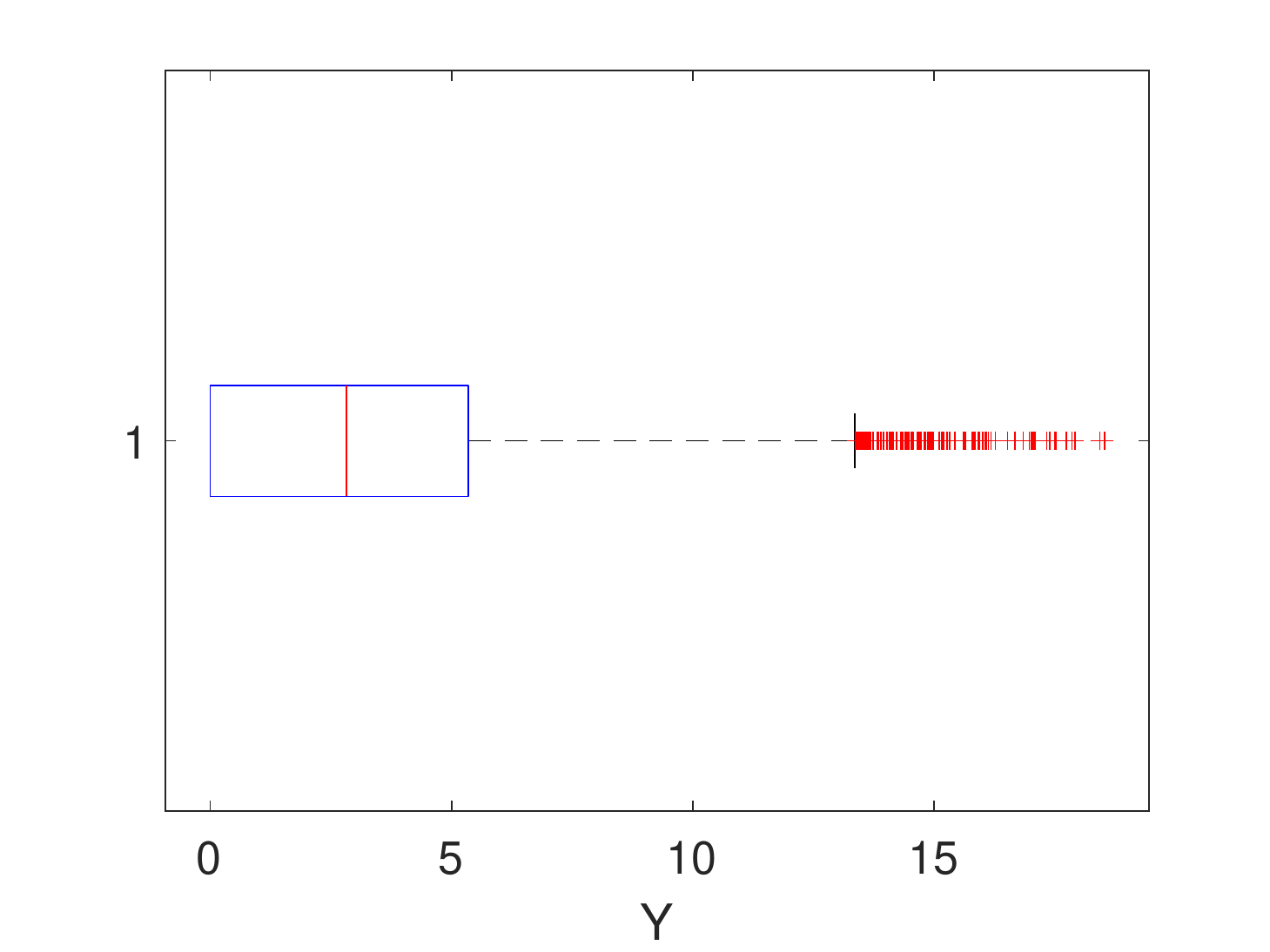}
		\includegraphics[width = .45\textwidth]{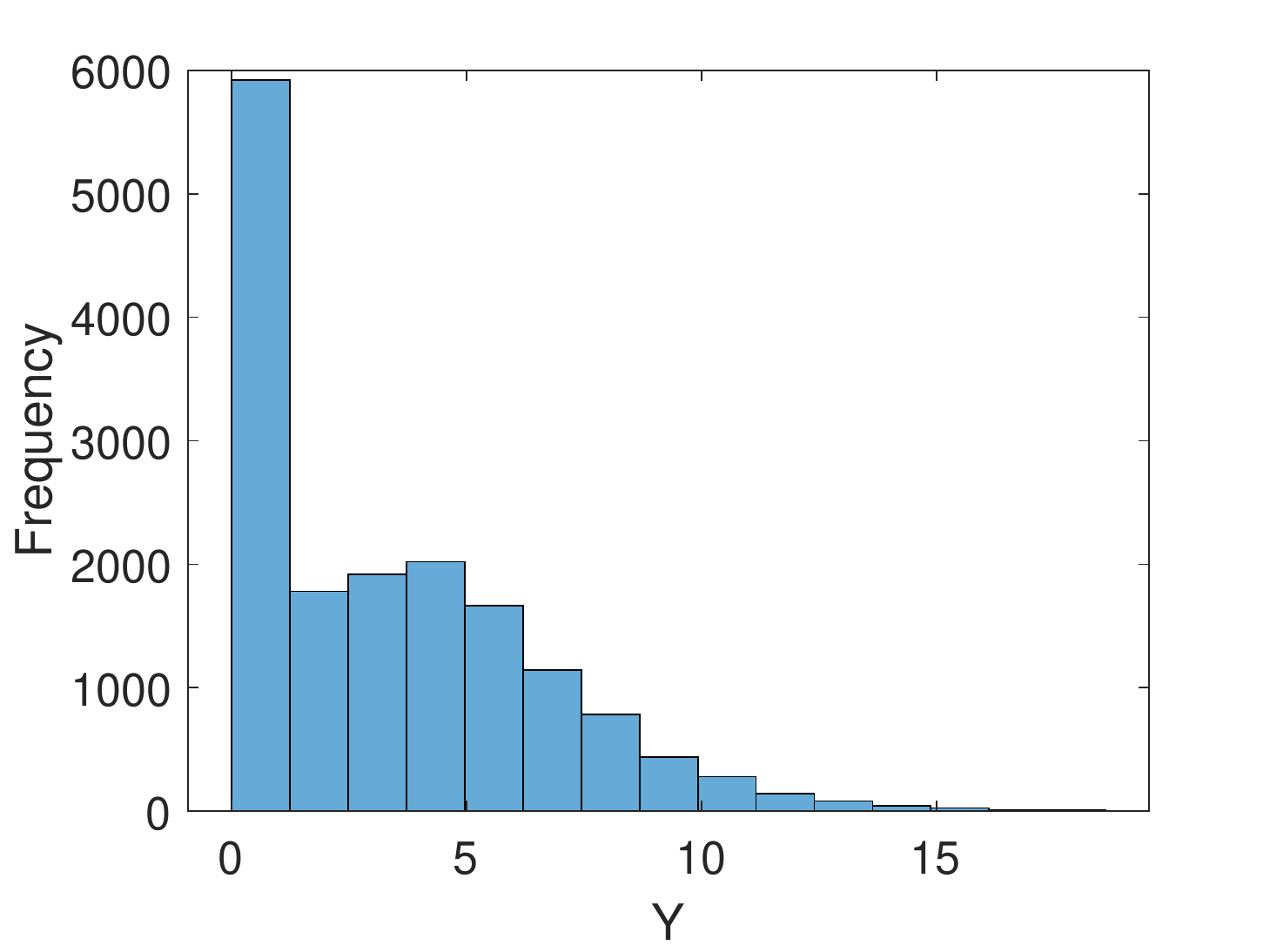}
		
		\centering
		\includegraphics[width=.45\textwidth]{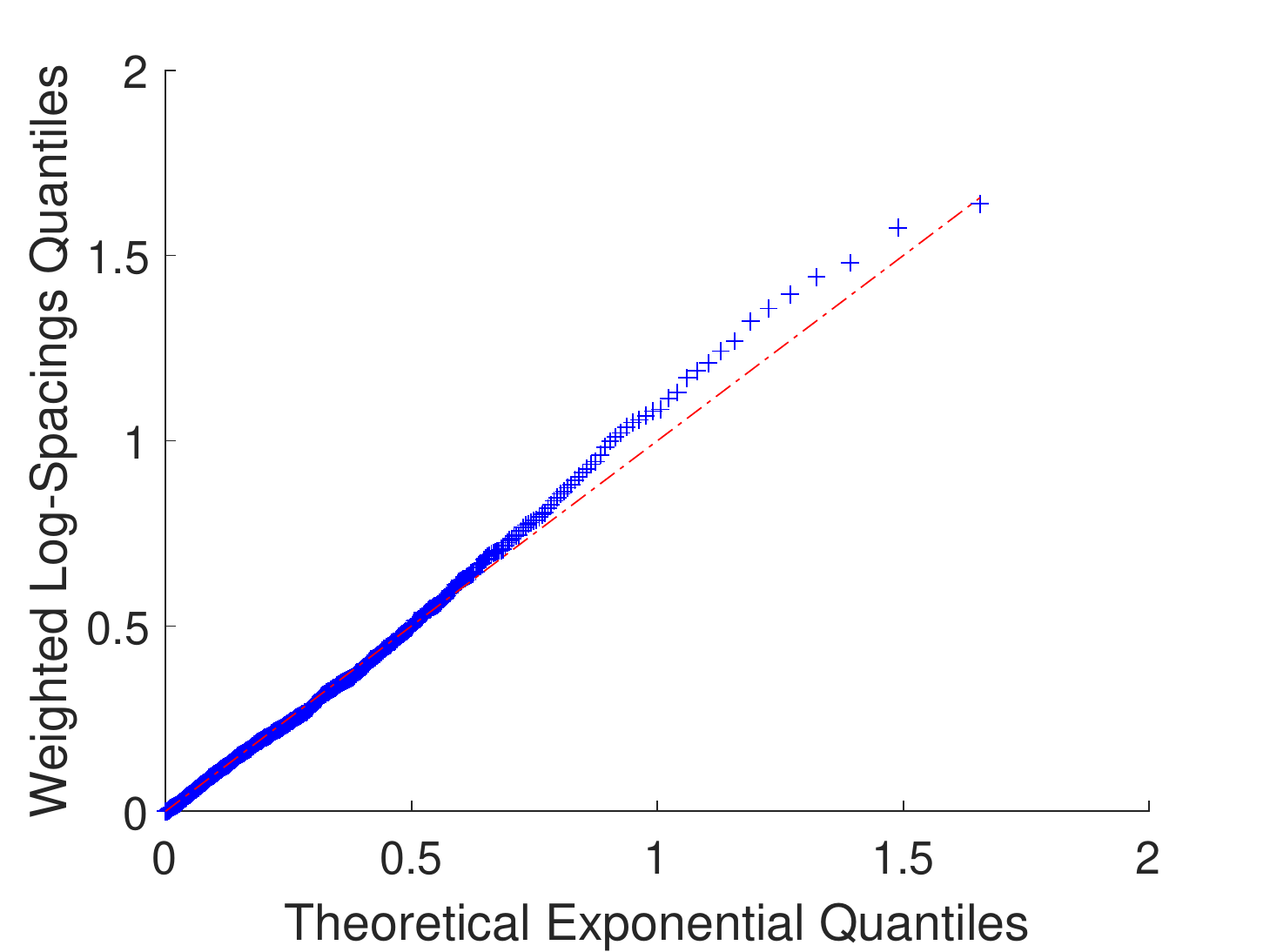}
		\caption{The boxplot (top left), histogram (top right) of the Box-Cox transformed campaign contribution and the Q-Q plot of $\{Z^{(i)}\}_{i=1}^{\lfloor N/8 \rfloor}$ versus exponential distribution (bottom).} \label{fig:USCampaign:Distribution_Y}
	\end{figure}

	The raw data of the campaign contributions are highly right-skewed. \cite{huang2021unified} searched across Box-Cox transformations of the campaign contribution data, $\text{BoxCox}(\text{Contribution},\lambda_1,\lambda_2):=\{(\text{Contribution}+\lambda_2)^{\lambda_1}-1\}/\lambda_1$, with respect to $\lambda_1, \lambda_2$, whose sample quantiles have the largest correlation with those of a standard normal distribution. This yielded $(\tilde{\lambda}_1,\tilde{\lambda}_2) = (0.1397,0.0176)$. They then took 
	\begin{equation}
		Y :=\text{BoxCox}(\text{Contribution},\tilde{\lambda}_1,\tilde{\lambda}_2)- \min\big\{\text{BoxCox}(\text{Contribution},\tilde{\lambda}_1,\tilde{\lambda}_2)\big\}\,,\label{TransContributions}
	\end{equation}
	so that the minimum response data is 0 and fitted a Tobit model to the data to estimate the average dose-response function $\mathbb{E}\{Y(t)\}$.

	In this paper, we focus on the quantiles of the potential outcome $Y(t)$ where the distribution of $Y(t)$ has a heavy tail. We first exam the heavy-tail assumption of the data. From the boxplot and histogram in the top row of Figure~\ref{fig:USCampaign:Distribution_Y}, we see that even the response data have been transformed as closer to a normal distribution as possible, the corresponding empirical distribution still has a heavy right tail. Moreover, let $Y^{(1)}\geq Y^{(2)}\geq \ldots Y^{(N)}$ be the order statistics of the sample $Y_1,\ldots, Y_N$, and $Z^{(i)} = i\log(Y^{(i)}/Y^{(i+1)})$ for $1\leq i \leq N-1$. It is known that if $Y$ is heavy-tailed, then, for small $i$, the $Z^{(i)}$'s are approximately independent copies of an exponential random variable (see e.g.~\citealp{beirlant2004statistics},~page~109--110). The bottom panel of Figure~\ref{fig:USCampaign:Distribution_Y} gives the quantile-quantile plot of the $Z^{(i)}$ for $i=1,\ldots,\lfloor N/8 \rfloor$. The approximately linear relationship in the plot further confirmed the heavy-tail assumption of our response data. This was ignored in the literature where linear models or normal residuals are assumed (e.g.~\citealp{Fong_Hazlett_Imai_2018, Ai_Linton_Motegi_Zhang_cts_treat, huang2021unified}).

	We then apply our method to estimate the extreme quantiles of the potential campaign contributions $Y(t)$ given an amount of political advertisements aired in the area. The treatment variable is a transform of the raw number of advertisements: Letting $\tilde{T} = \log(\# \text{advertisements}+1)$, $\tilde{T}_{\max} = \max(\tilde{T}_1,\ldots, \tilde{T}_N)$ and $\tilde{T}_{\min} = \min (\tilde{T}_1,\ldots, \tilde{T}_N)$, we take $T := (\tilde{T} - \tilde{T}_{\min}) / (\tilde{T}_{\max} - \tilde{T}_{\min})$ as the treatment variable. Following the literature, the covariates considered were 
	$$
	\boldsymbol{X} =\begin{bmatrix}
		\log(\text{Population})\\
		\%\text{Age over 65}\\
		\log(\text{Median Income})\\
		\%\text{Hispanic}\\
		\%\text{Black}\\
		\log(\text{Population density}+1)\\
		\%\text{College graduates}\\
		\mathbbm{1}(\text{Can commute to a competitive state})
	\end{bmatrix} \,.
	$$
	The definition of each covariate is almost self-explanatory, and one can refer to \cite{Fong_Hazlett_Imai_2018} for more details.
	
	\begin{figure}[t]
		\centering
		\includegraphics[width = .45\textwidth]{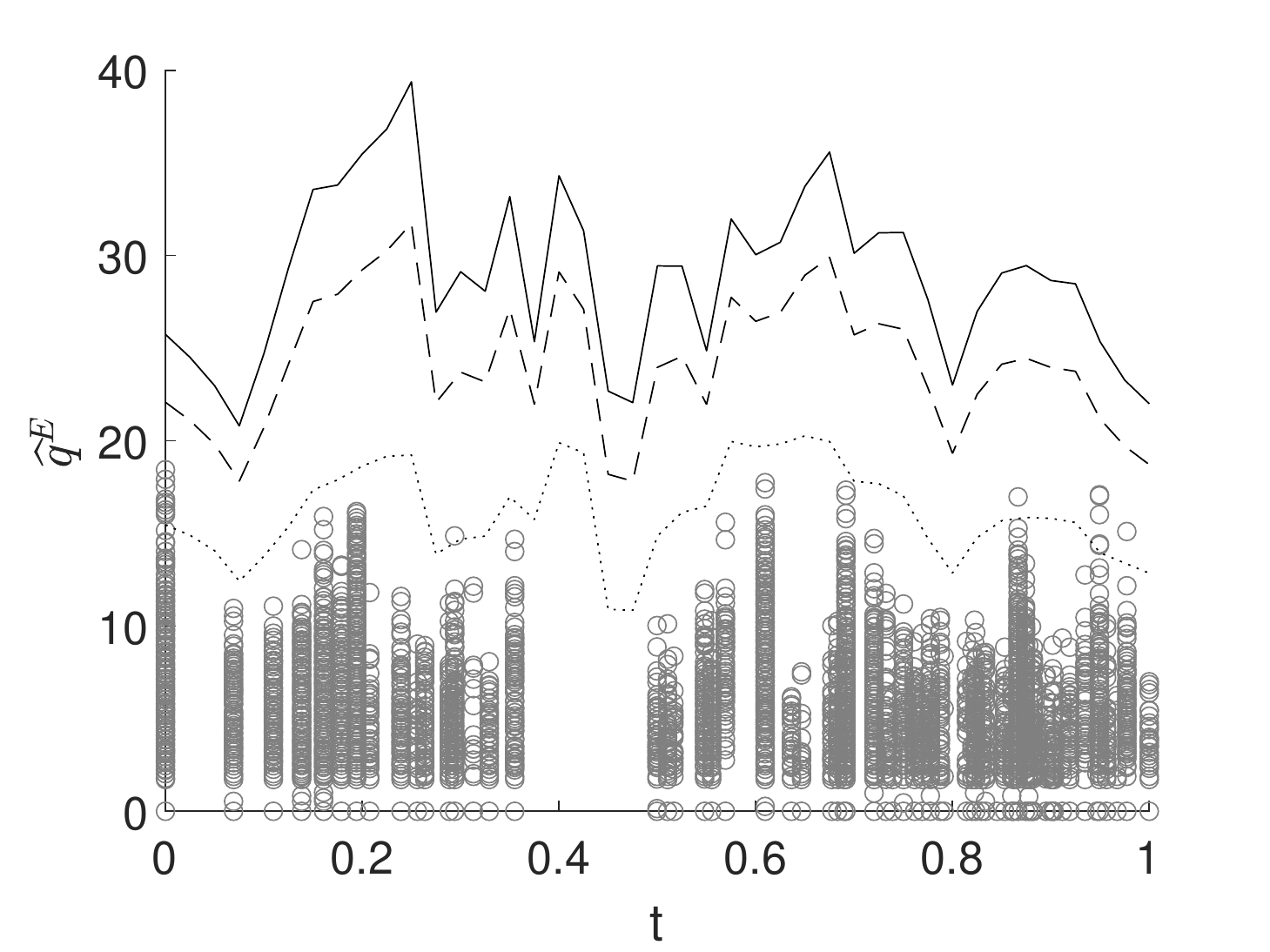}
		\includegraphics[width = .45\textwidth]{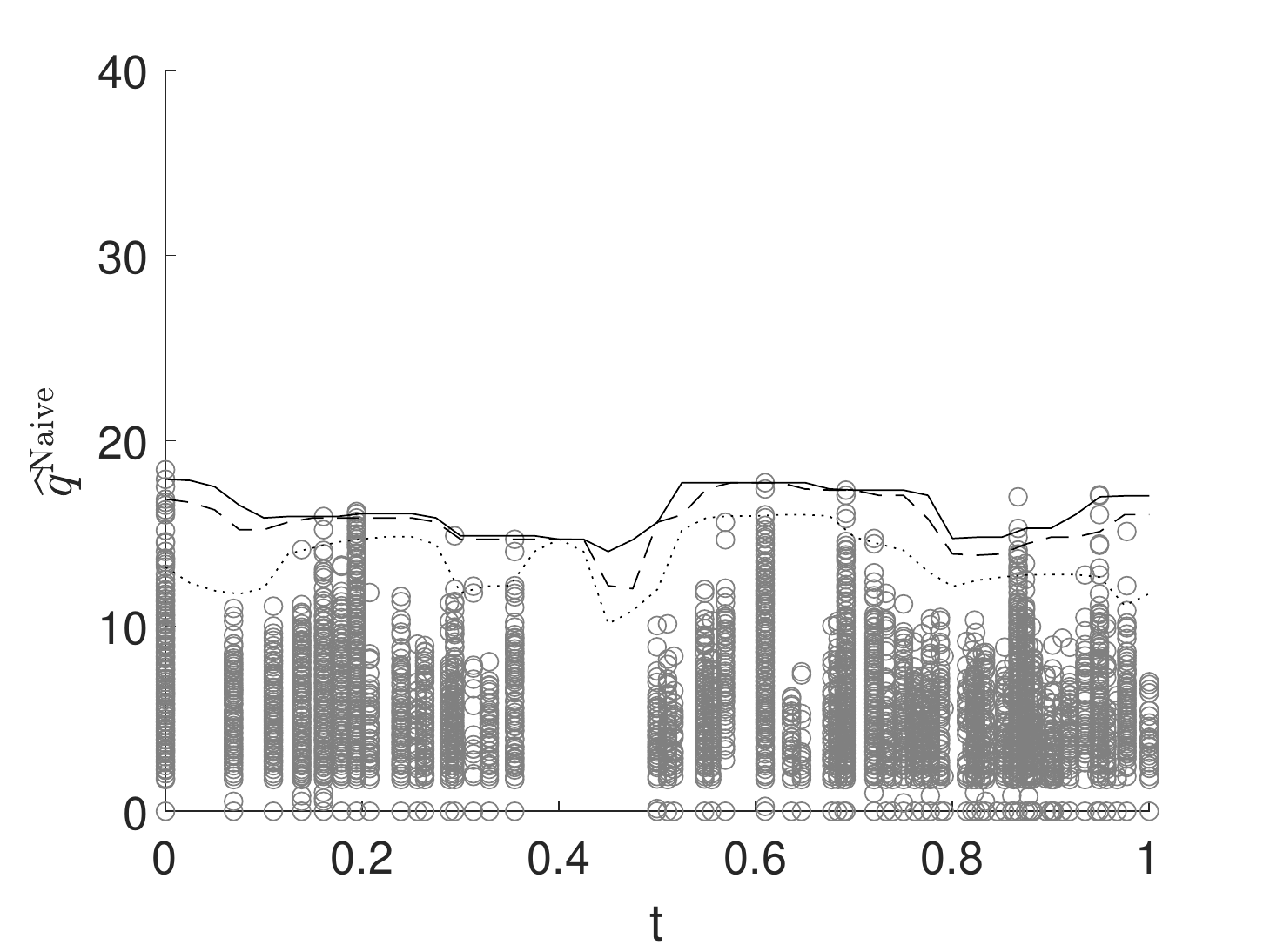}
		
		%\centering
		%\includegraphics[width=.45\textwidth]{Gammahat.pdf}
		\caption{Our estimated extreme quantiles (top left) and naive extreme quantiles (top right) of the potential Box-Cox transformed campaign contributions at $\alpha_N = 0.995$ (orange), $0.999$ (purple) and $0.9995$ (yellow)%; and the estimated index function (bottom).
		} \label{fig:USCampaign:Quantiles}
	\end{figure}
	
	\begin{figure}[t]
		\centering
		\includegraphics[width = .32\textwidth]{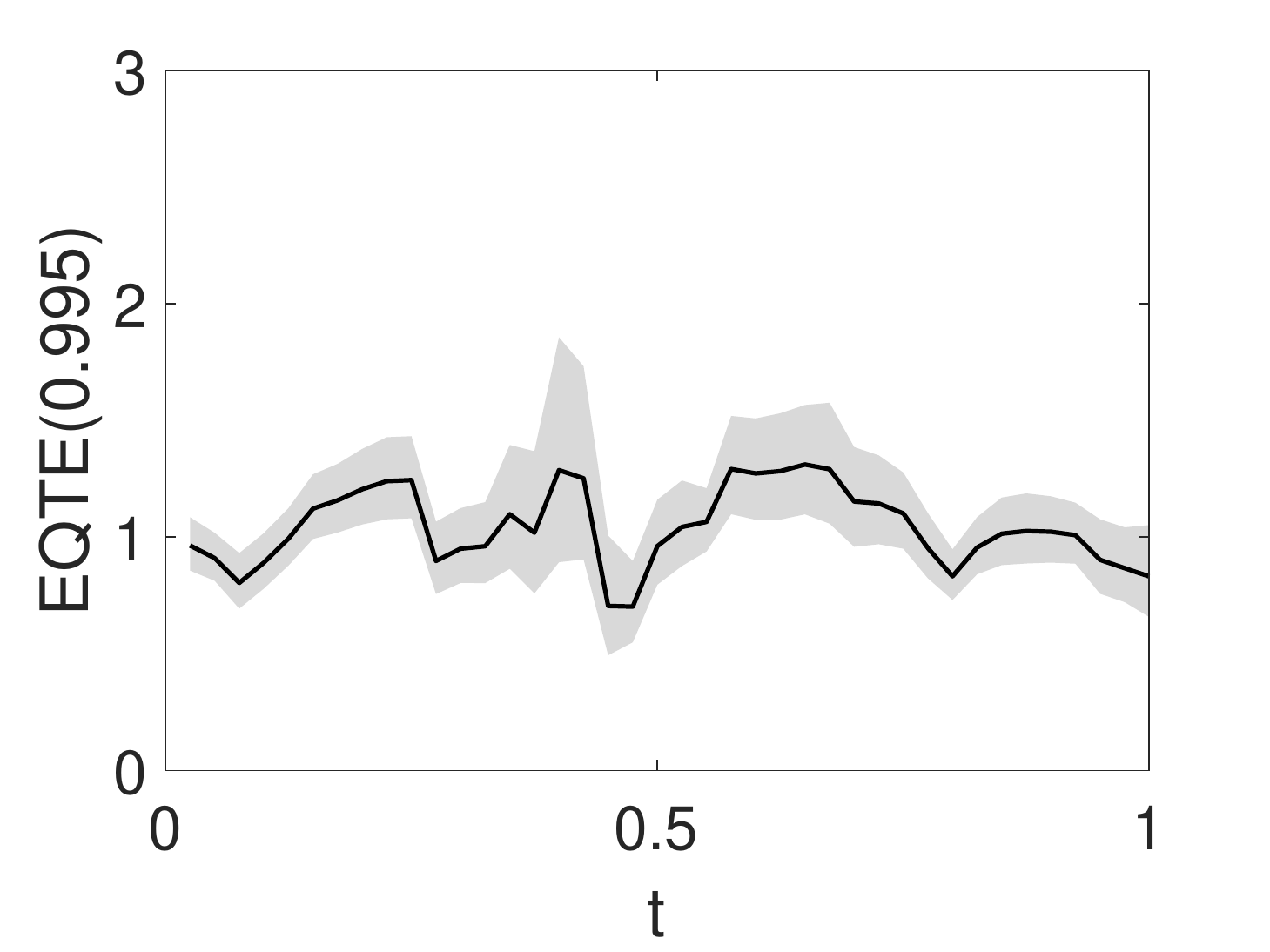}
		\includegraphics[width = .32\textwidth]{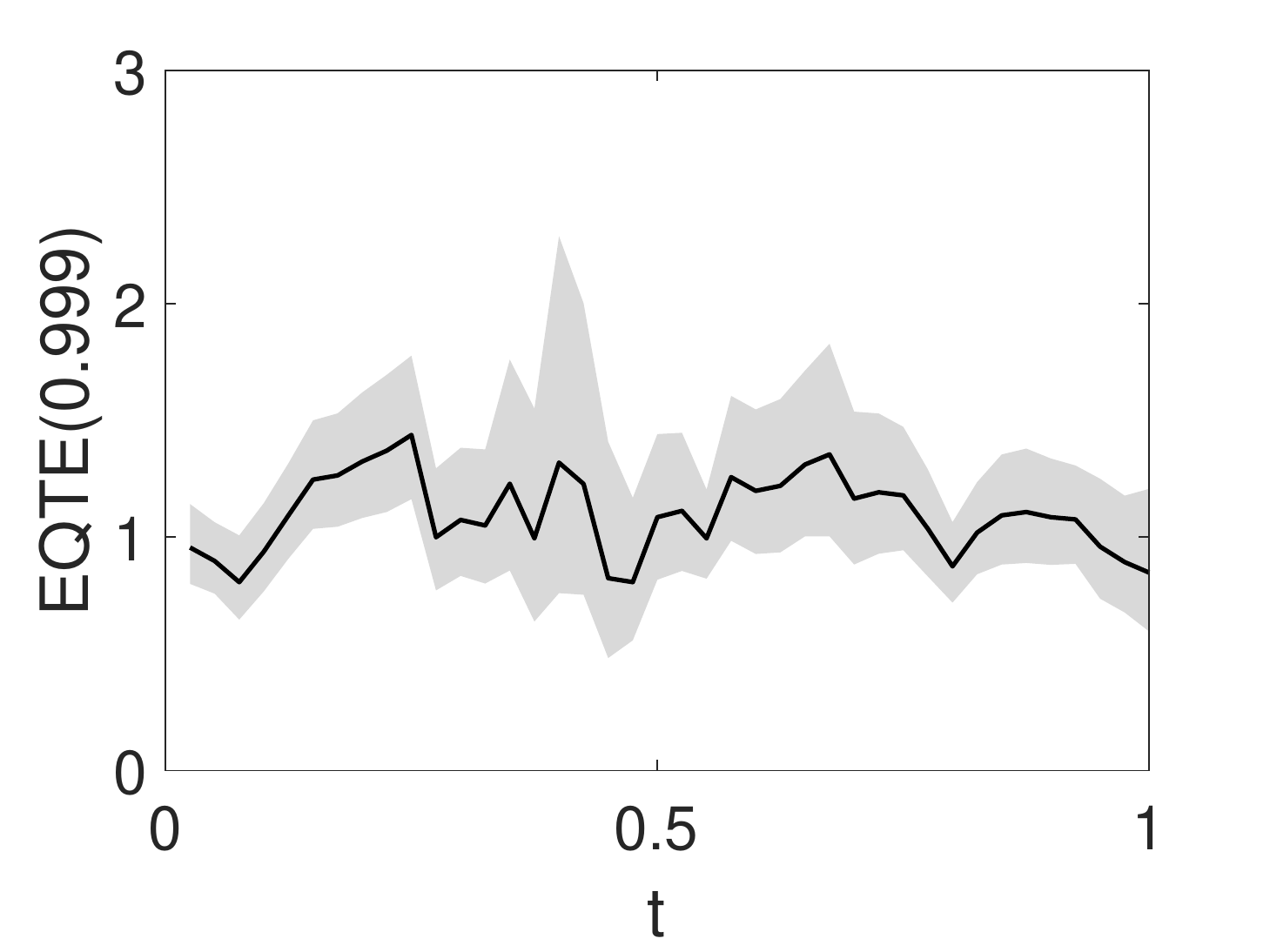}
		\includegraphics[width = .32\textwidth]{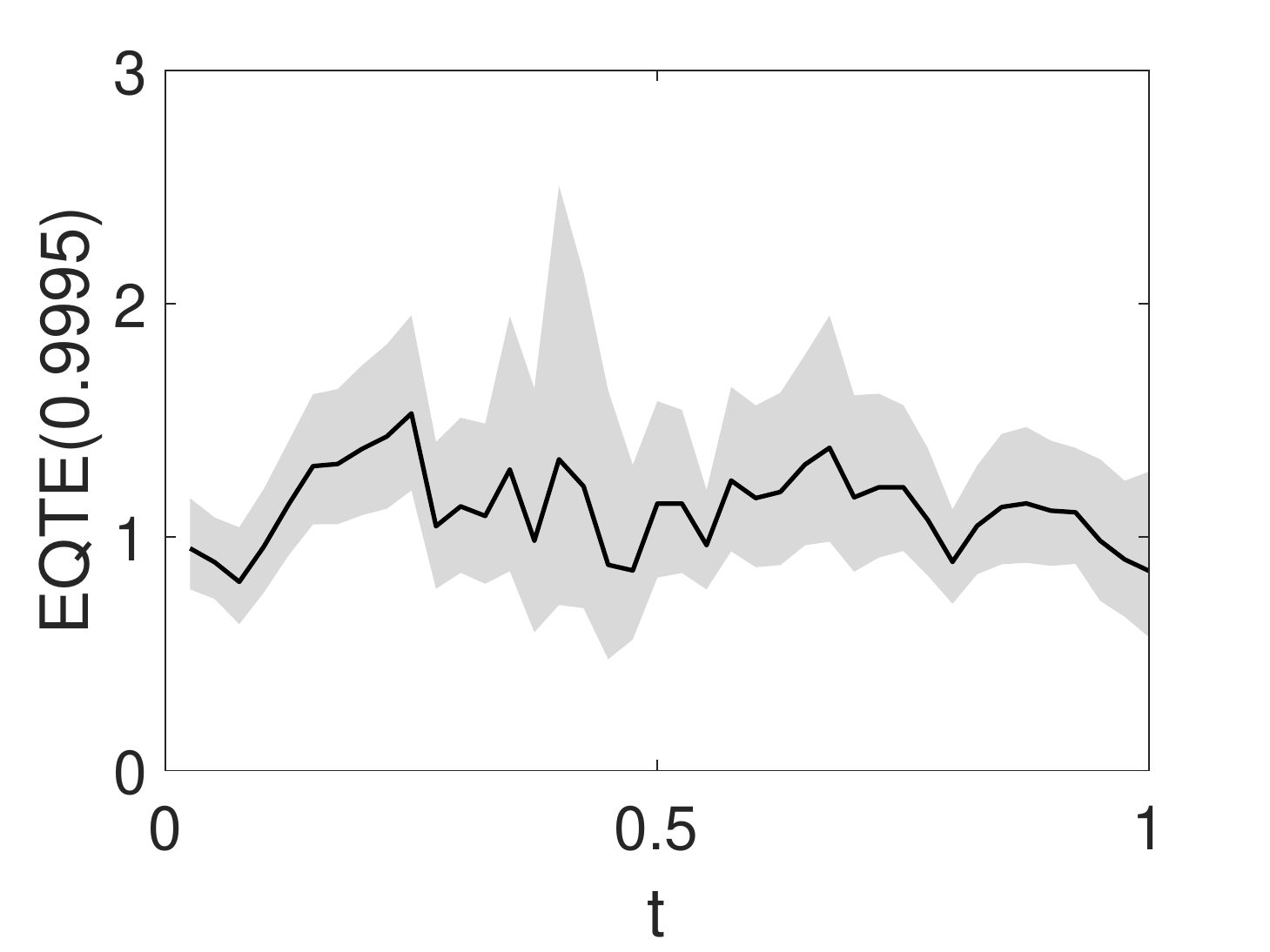}
		
		\includegraphics[width = .32\textwidth]{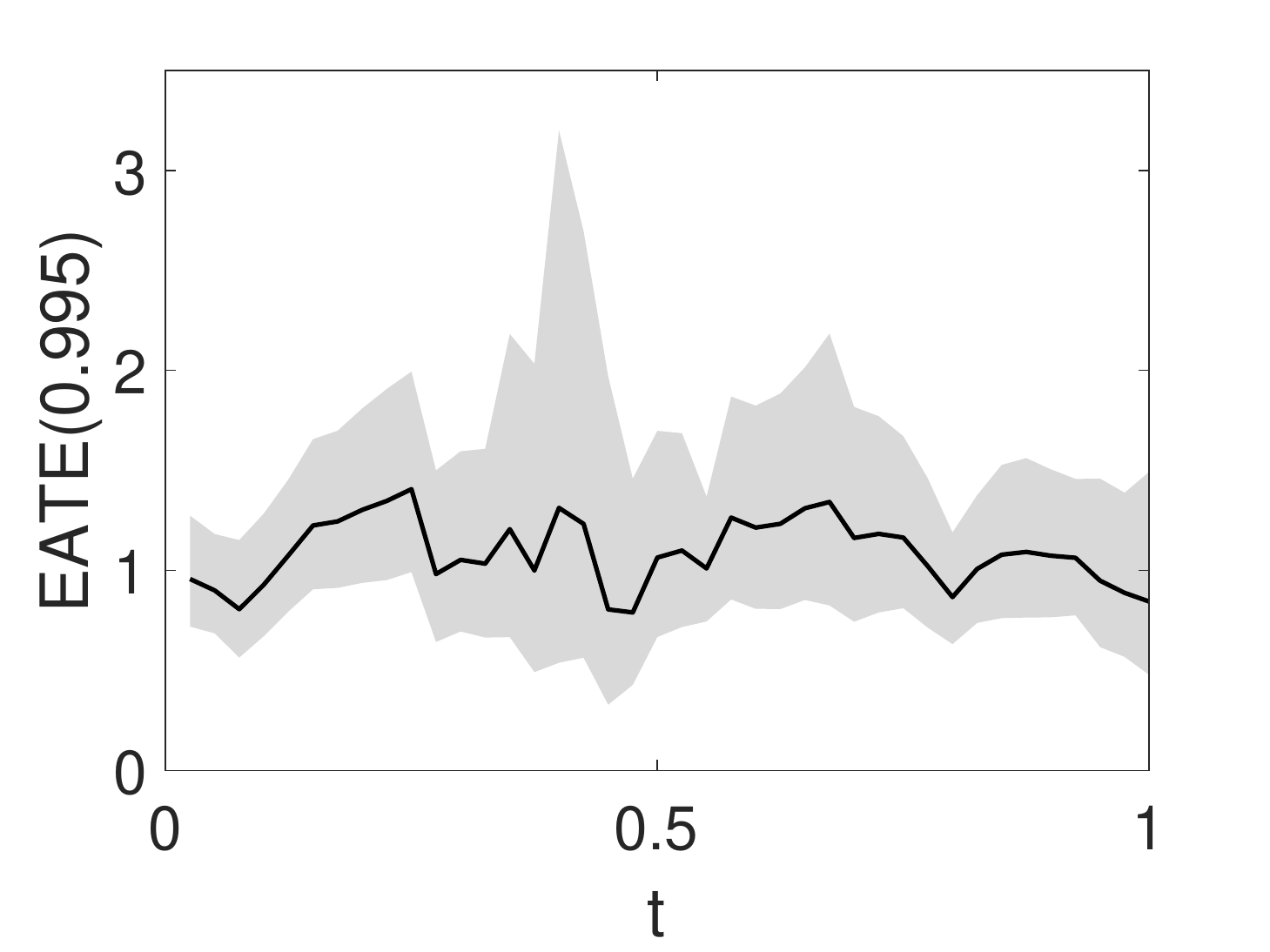}
		\includegraphics[width = .32\textwidth]{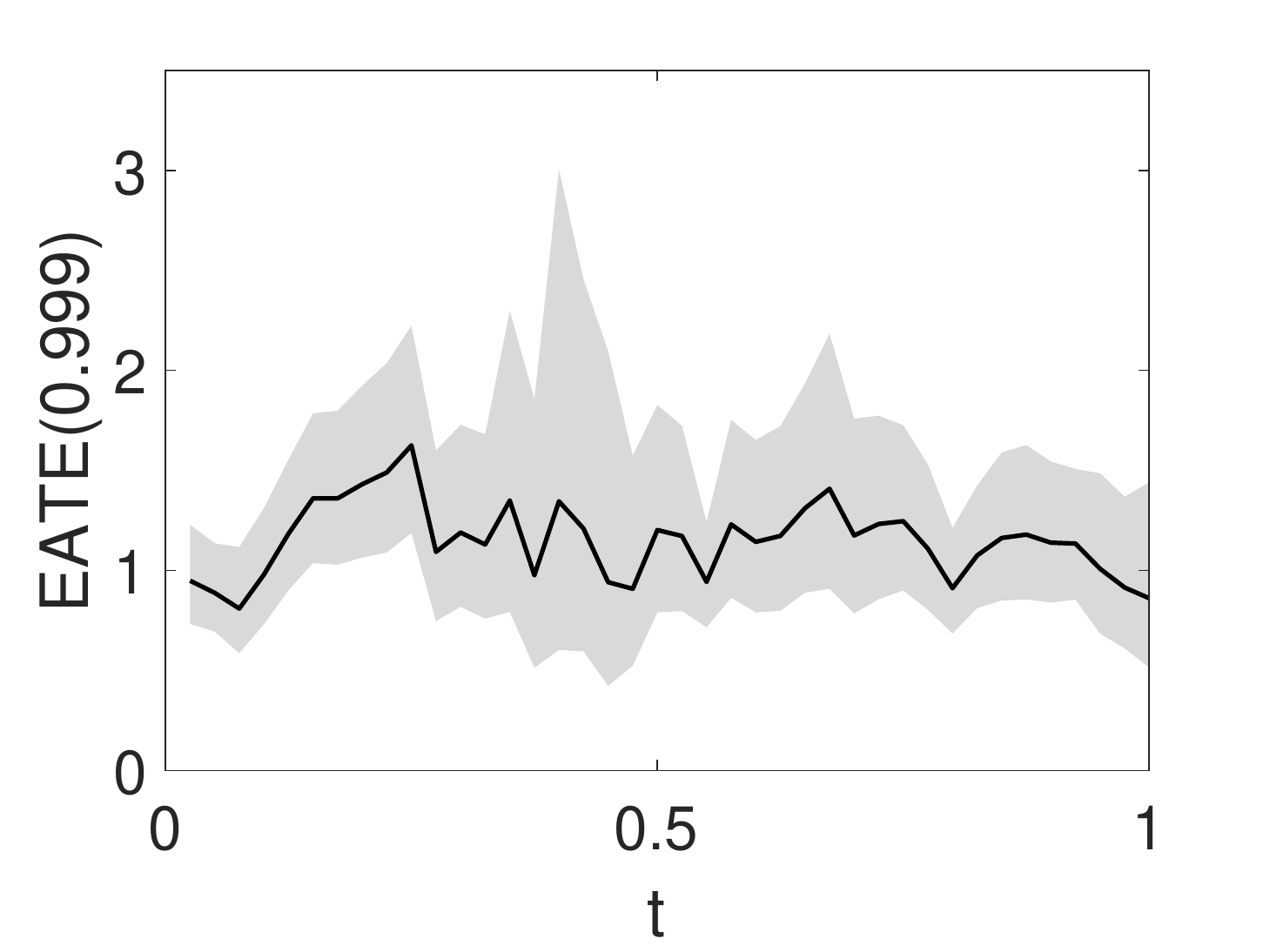}
		\includegraphics[width = .32\textwidth]{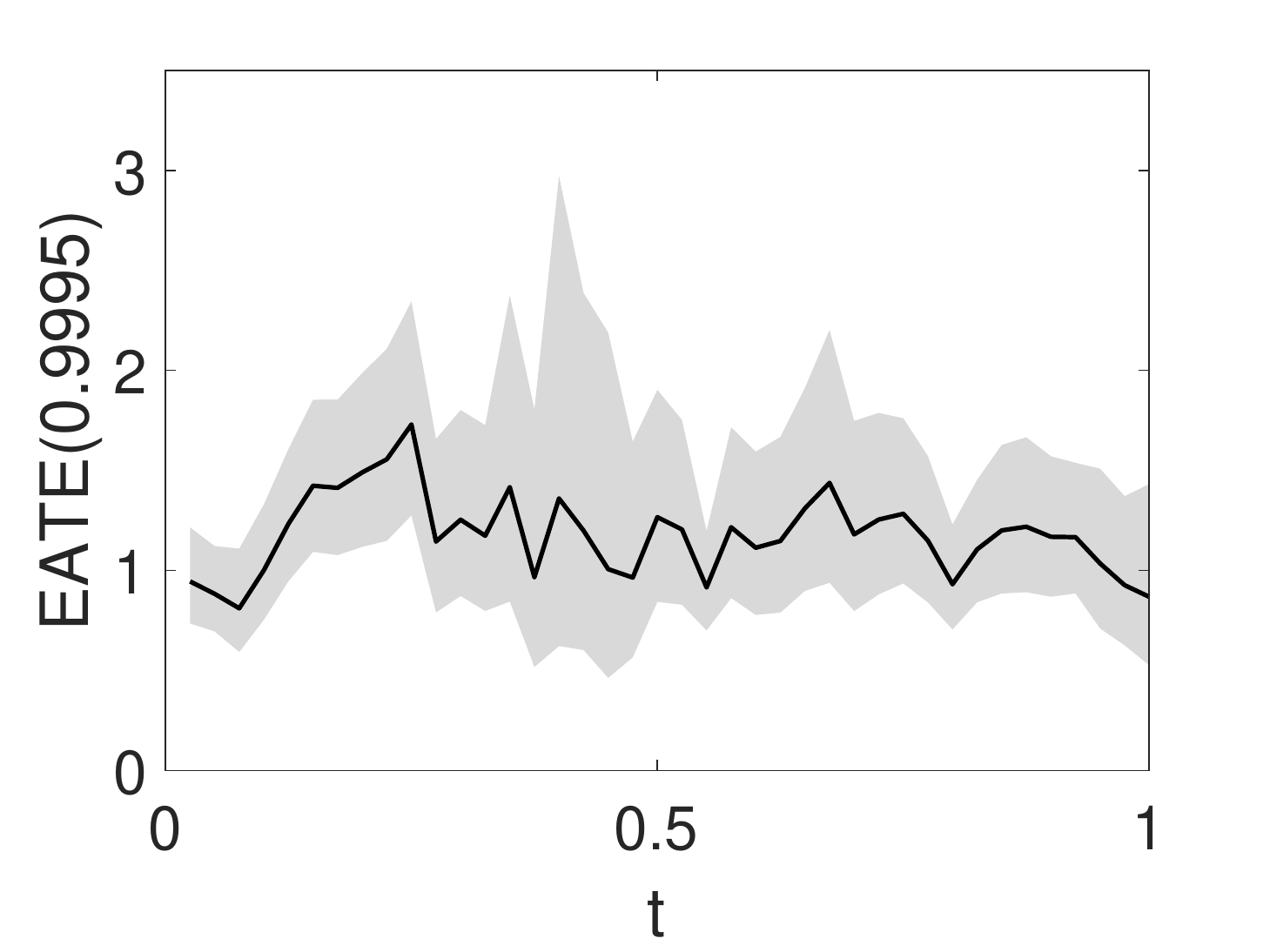}
		
		\caption{Estimated $\text{EQTE}_{0,t}(\alpha_N)$ (top) and $\text{EQTE}_{0,t}(\alpha_N)$ (bottom), for $\alpha_N=0.995$ (left), 0.999 (centre) and 0.9995(right) with 95\% confidence intervals.
		} \label{fig:USCampaign:TEs}
	\end{figure}
	
	Figure~\ref{fig:USCampaign:Quantiles} gives the plots of
	our estimated extreme quantiles $\widehat{q}^E_t(\alpha_N)$ and the naive extreme quantiles $\widehat{q}^{\text{Naive}}_{t,h}(\alpha_N)$ defined in \eqref{Def:NaiveQuantiles} at $\alpha_N = 0.995, 0.999$ and $0.9995$. Figure~\ref{fig:USCampaign:TEs} further depicts the estimated $\text{EQTE}_{0,t}(\alpha_N)$ and $\text{EATE}_{0,t}(\alpha_N)$ for $0<t\leq 1$. Note that in the literature (see e.g.~\citealp{huang2021unified}), it is found that with a relatively small increase in the number of political advertisements from 0, the average potential campaign contribution increases rapidly. Then the improvement gradually becomes marginal as the number of advertisements gets large. Now, our results show that this conclusion also applies to extreme cases. That is, even in the area where the amount of campaign contributions is extremely large, enough political advertisements still play an essential role in affecting the potential campaign contributions. After a sufficient amount of advertisements has been aired, more advertisements do not seem to help. On the contrast, the naive method is not capable of extrapolating and gives nearly flat extreme quantile curves due to insufficient empirical extreme quantile data. Thus, its extreme result is unreliable.
	
		\section{Conclusion and Discussion}
	    In this paper, we consider estimation and inference for extreme treatment effects in the continuous treatment effect model. We define the EQTE and EATE through the quantile and tail mean function. We derive the limiting theory by employing established results in EVT. Depending on whether extrapolation is needed, we distinguish our limiting results between the intermediate scenario and extreme scenario. In addition, our limiting results are for EQTE and EATE processes indexed by a range of tail levels and hence facilitates uniform inference. Simulations confirm that our estimators have favorable accuracy and the confidence intervals enjoy satisfactory coverage. An empirical application demonstrates the practical usefulness of our approach. For future research, it is of interest to consider other tail measures such as the expectile \citep[e.g.,][]{Daouia2017}.
	%\pagebreak
	
	\bigskip
	\section*{Supplementary Material}
	%Supplementary materials are only for online publication. 
	The supplemental material contains preliminary lemmas and proofs of main results in the paper.
	
	\subsection{Preliminary lemmas}
	
	\begin{lemma}\label{lemma:02}
		Let $y_N\to \infty$ and $h\to 0$ as $N\to\infty$. In addition, for any fixed $t\in\mathcal{T}$, $\omega_t(y_{N},h)\log(y_N)\to0$.
		Then %uniformly in $y^{\prime}\in[y_N, \Delta y_N]$ and $t^{\prime}\in B(t, h)$,
		\begin{equation}
			\sup_{y\in[y_N,\Delta y_N], t^{\prime}\in B(t, h)}\left|\frac{\bar{F}_{t^{\prime}}(y)}{\bar{F}_{t}(y)} - 1\right| = O\left(\omega_{t}(y_{N},h)\log(y_N)\right).
		\end{equation}
		%where $\Delta\geq 1$ is a constant.
	\end{lemma}
	
	\begin{proof}
		The definition of $\omega_{t}(y_N,h)$ in \eqref{eq:oscillation} indicates that 
		\begin{eqnarray}
			\notag & & \frac{1}{\log(\Delta y_N)}\sup_{y\in[y_N,\Delta y_N], t^{\prime}\in B(t, h)}\left|\log\left\{\frac{\bar{F}_{t^{\prime}}(y)}{\bar{F}_{t}(y)}\right\}\right| \\
			\notag & \leq & \sup_{y\in[y_N,\Delta y_N], t^{\prime}\in B(t, h)}\frac{1}{\log(y)}\left|\log\left\{\frac{\bar{F}_{t^{\prime}}(y)}{\bar{F}_{t}(y)}\right\}\right|  \leq  \omega_{t}(y_N,h),
		\end{eqnarray}
		from which we know that
		\begin{equation}
			\notag \sup_{y\in[y_N,\Delta y_N], t^{\prime}\in B(t, h)}\left|\log\left\{\frac{\bar{F}_{t^{\prime}}(y)}{\bar{F}_{t}(y)}\right\}\right| = O(\omega_{t}(y_N, h)\log(\Delta y_N)) = O\left(\omega_{t}(y_{N},h)\log(y_N)\right).
		\end{equation}
		Thus, we finish the proof of this lemma by the fact that $|e^x-1|\leq 2x$ uniformly for $x\in[0,1]$.
	\end{proof}

	\begin{lemma}\label{lemma:03}
		Suppose Assumptions \ref{UnconfoundAssump}, \ref{as:suppX}, \ref{assump:f_T} and \ref{asaump:K}  hold. Let $y_N\to \infty$ and $h\to 0$ as $N\to\infty$. In addition, for any fixed $t\in\mathcal{T}$, $\omega_t(y_{N},h)\log(y_N)\to0$. 
		Then, for any fixed $t\in\mathcal{T}$, \\
		(i) $\mathbb{E}\left\{\pi_{0}(T,\mathbf{X})\frac{1}{h}K\left(\frac{T-t}{h}\right)\right\} = \mathbb{E}\left\{\frac{1}{h}K\left(\frac{T-t}{h}\right)\right\}=f_{T}(t)+\frac{\kappa_{21}}{2}f^{\prime\prime}_{T}(t)h^2+O(h^3)$. \\
		(ii)
		$\frac{\mathbb{E}\left\{\pi_0(T,\bs{X})\mathbbm{1}(Y> y)\frac{1}{h}K\left(\frac{T-t}{h}\right)\right\}}{\bar{F}_{t}(y)f_{T}(t)} - 1 = O\left(\omega_{t}(y_{N},h)\log(y_N)+h\right)
		$ uniformly in $y\in[y_N,\Delta y_N]$. \\
		(iii) Denote $\bar{F}_{t,h}(y)=\frac{\mathbb{E}\left\{\pi_0(T,\bs{X})\mathbbm{1}(Y> y)\frac{1}{h}K\left(\frac{T-t}{h}\right)\right\}}{\mathbb{E}\left\{\pi_{0}(T,\mathbf{X})\frac{1}{h}K\left(\frac{T-t}{h}\right)\right\}}$, then uniformly in $y\in[y_N,\Delta y_N]$,
		\begin{eqnarray}
			\notag & & \frac{\bar{F}_{t,h}(y)}{\bar{F}_{t}(y)} - 1 = O\left(\omega_{t}(y_{N},h)\log(y_N)+h\right).
		\end{eqnarray}
	\end{lemma}
	
	\begin{proof}
		(i) First, by simple algebra,
		\begin{eqnarray}
			\notag & &  \mathbb{E}\left\{\pi_{0}(T,\mathbf{X})\frac{1}{h}K\left(\frac{T-t}{h}\right)\right\} = \int\pi_{0}(t^{\prime}, \mathbf{x})\frac{1}{h}K\left(\frac{t^{\prime}-t}{h}\right)f_{T,\mathbf{X}}(t^{\prime},\mathbf{x})d\mathbf{x}dt^{\prime} \\
			\notag & = & \int\frac{f_{T}(t^{\prime})}{f_{T|\mathbf{X}}(t^{\prime}|\mathbf{x})}\frac{1}{h}K\left(\frac{t^{\prime}-t}{h}\right)f_{T|\mathbf{X}}(t^{\prime}|\mathbf{x})f_{\mathbf{X}}(\mathbf{x})d\mathbf{x}dt^{\prime} \\
			\notag & = & \int f_{\mathbf{X}}(\mathbf{x})d\mathbf{x}\int f_{T}(t^{\prime})\frac{1}{h}K\left(\frac{t^{\prime}-t}{h}\right)dt^{\prime} = \mathbb{E}\left\{\frac{1}{h}K\left(\frac{T-t}{h}\right)\right\}.
		\end{eqnarray}
		Then, by change of variable and Taylor's expansion, we have
		\begin{eqnarray}
			\notag & & \mathbb{E}\left\{\pi_{0}(T,\mathbf{X})\frac{1}{h}K\left(\frac{T-t}{h}\right)\right\}  = \int f_{T}(t^{\prime})\frac{1}{h}K\left(\frac{t^{\prime}-t}{h}\right)dt^{\prime} = \int_{\mathcal{K}} f_{T}(t+hu)K\left(u\right)du \\
			\notag & = & f_{T}(T)\int_{\mathcal{K}} K(u)du + hf_{T}^{\prime}(t)\int_{\mathcal{K}} uK(u)du + \frac{h^2}{2}f_{T}^{\prime\prime}(t)\int_{\mathcal{K}} u^2K(u)du + O(h^3) \\
			\notag & = & f_{T}(t) + \frac{\kappa_{21}}{2}f_{T}^{\prime\prime}(t)h^2 + O(h^3). 
		\end{eqnarray}%\textcolor{red}{(may not need to expand to the third order)}

		(ii) First, we can rewrite
		\begin{eqnarray}
			\notag & & \mathbb{E}\left\{\pi_0(T,\bs{X})\mathbbm{1}(Y> y)\frac{1}{h}K\left(\frac{T-t}{h}\right)\right\} \\
			\notag & = & \int\frac{1}{h}K\left(\frac{t^{\prime}-t}{h}\right)\frac{f_{T}(t^{\prime})}{f_{T|\mathbf{X}}(t^{\prime}|\mathbf{x})}\mathbb{P}(Y(t^{\prime})> y|T=t^{\prime},\mathbf{X}=\mathbf{x})f_{T|\mathbf{X}}(t^{\prime}|\mathbf{x})f_{\mathbf{X}}(\mathbf{x})dt^{\prime}d\mathbf{x} \\
			\notag & = & \int\frac{1}{h}K\left(\frac{t^{\prime}-t}{h}\right)f_{T}(t^{\prime})\mathbb{P}(Y(t^{\prime})> y|\mathbf{X}=\mathbf{x})f_{\mathbf{X}}(\mathbf{x})dt^{\prime}d\mathbf{x} \\
			\notag & = & \int\frac{1}{h} K\left(\frac{t^{\prime}-t}{h}\right)f_{T}(t^{\prime})\bar{F}_{t^{\prime}}(y)dt^{\prime},
		\end{eqnarray}
		where the second equation utilizes Assumption \ref{UnconfoundAssump}.
		By change of variable,
		\begin{eqnarray}
			\notag & & \int\frac{1}{h} K\left(\frac{t^{\prime}-t}{h}\right)f_{T}(t^{\prime})\bar{F}_{t^{\prime}}(y)dt^{\prime}=  \int_{\mathcal{K}} K\left(u\right)f_{T}(t+hu)\bar{F}_{t+hu}(y)du.
		\end{eqnarray}
		Now consider
		\begin{eqnarray}
			\notag & & \left|\frac{\int_{\mathcal{K}} K\left(u\right)f_{T}(t+hu)\bar{F}_{t+hu}(y)du}{\bar{F}_{t}(y)f_{T}(t)} -1 \right| \\
			\notag & \leq & \frac{1}{\bar{F}_{t}(y)f_{T}(t)}\left|\int_{\mathcal{K}} K\left(u\right)f_{T}(t+hu)\bar{F}_{t+hu}(y)du - \int_{\mathcal{K}} K\left(u\right)f_{T}(t+hu)\bar{F}_{t}(y)du\right| \\
			\notag & & + \frac{1}{\bar{F}_{t}(y)f_{T}(t)}\left|\int_{\mathcal{K}} K\left(u\right)f_{T}(t+hu)\bar{F}_{t}(y)du - \bar{F}_{t}(y)f_{T}(t)\right| \\
			\notag & \leq & \frac{1}{f_{T}(t)}\int_{\mathcal{K}} K\left(u\right)f_{T}(t+hu)\left|\frac{\bar{F}_{t+hu}(y)}{\bar{F}_{t}(y)}-1\right|du + \frac{1}{f_{T}(t)}\int_{\mathcal{K}} K\left(u\right)\left|f_{T}(t+hu)-f_{T}(t)\right|du.
		\end{eqnarray}
		According to Lemma \ref{lemma:02}, as $0<\eta_{3}<f_{T}(t)<\eta_{4}<\infty$ for all $t\in\mathcal{T}$ and $\mathcal{K}\subseteq[-1,1]$ implies $\{t+hu: u\in\mathcal{K}\}\subset B(t, h)$, we have uniformly in $y\in[y_N, \Delta y_N]$,
		\begin{eqnarray}
			\notag & & \frac{1}{f_{T}(t)}\int_{\mathcal{K}} K\left(u\right)f_{T}(t+hu)\left|\frac{\bar{F}_{t+hu}(y)}{\bar{F}_{t}(y)}-1\right|du \\
			\notag & \leq & \frac{1}{f_{T}(t)}\int_{\mathcal{K}} K\left(u\right)f_{T}(t+hu)O\left(\omega_{t}(y_{N},h)\log(y_N)\right)du \\
			\notag & = & \frac{1}{f_{T}(t)}\int_{\mathcal{K}} K\left(u\right)f_{T}(t+hu)du\times O\left(\omega_{t}(y_{N},h)\log(y_N)\right) \\
			& = & O\left(\omega_{t}(y_{N},h)\log(y_N)\right). \label{eq:lemma_02_01}
		\end{eqnarray}
		Furthermore, 
		\begin{eqnarray}
			\notag & & \frac{1}{f_{T}(t)}\int_{\mathcal{K}} K\left(u\right)\left|f_{T}(t+hu)-f_{T}(t)\right|du \\
			& \leq & \frac{1}{f_{T}(t)}\left\{f^{\prime}_{T}(t)h\int |u|K\left(u\right)du+O(h^2)\right\} = O(h) \label{eq:lemma_02_02}.
		\end{eqnarray}
		Finally, \eqref{eq:lemma_02_01} and \eqref{eq:lemma_02_02} imply that
		\begin{eqnarray}
			\notag & & \sup_{y\in[y_N,\Delta y_N]}\left|\frac{\mathbb{E}\left\{\pi_0(T,\bs{X})\mathbbm{1}(Y> y)\frac{1}{h}K\left(\frac{T-t}{h}\right)\right\}}{\bar{F}_{t}(y)\}f_{T}(t)} - 1\right| = O\left(\omega_{t}(y_{N},h)\log(y_N)+h\right).
		\end{eqnarray} 
		
		(iii) In view of (i) and (ii), uniformly in $y\in[y_N, \Delta y_N]$,
		\begin{eqnarray}
			\notag & & \frac{\bar{F}_{t,h}(y)}{\bar{F}_{t}(y)} - 1 = \frac{1}{\bar{F}_{t}(y)}\left[\frac{\mathbb{E}\left\{\pi_0(T,\bs{X})\mathbbm{1}(Y> y)\frac{1}{h}K\left(\frac{T-t}{h}\right)\right\}}{\mathbb{E}\left\{\pi_{0}(T,\mathbf{X})\frac{1}{h}K\left(\frac{T-t}{h}\right)\right\}}-\bar{F}_{t}(y)\right] \\
			\notag & = & \frac{1}{\bar{F}_{t}(y)}\left[\frac{\bar{F}_{t}(y)f_{T}(t)\left\{1+O\left(\omega_{t}(y_{N},h)\log(y_N)+h\right)\right\}}{f_{T}(t)+\frac{\kappa_{21}}{2}f^{\prime\prime}_{T}(t)h^2+O(h^3)}-\bar{F}_{t}(y)\right] \\
			\notag & = & \frac{f_{T}(t)\left\{1+O\left(\omega_{t}(y_{N},h)\log(y_N)+h\right)\right\}-f_{T}(t)-\frac{\kappa_{21}}{2}f^{\prime\prime}_{T}(t)h^2+O(h^3)}{f_{T}(t)+\frac{\kappa_{21}}{2}f^{\prime\prime}_{T}(t)h^2+O(h^3)} \\
			\notag & = & O\left(\omega_{t}(y_{N},h)\log(y_N)+h\right).
		\end{eqnarray}
		Thus, we finish the proof of this lemma.
	\end{proof}

	\begin{lemma}\label{lemma:04}
		Suppose Assumptions \ref{assump:extreme_domain}, \ref{UnconfoundAssump}, \ref{as:suppX}, \ref{assump:f_T} and \ref{asaump:K}  hold. Let $y_N\to \infty$ and $h\to 0$ as $N\to\infty$. In addition, for any fixed $t\in\mathcal{T}$, $\omega_t(y_{N},h)\log(y_N)\to0$ and $Nh\bar{F}_{t}(y_N)\to\infty$.
		Denote $\Lambda_{N,1}(y,t,h) = N^{-1}\sum_{i=1}^{N}\pi_0(T_i,\mathbf{X}_i)\mathbbm{1}(Y_i>y)\frac{1}{h}K\left(\frac{T_i-t}{h}\right)$, $\Lambda_{N,2}(t,h) = N^{-1}\sum_{i=1}^{N}\pi_0(T_i,\mathbf{X}_i)\frac{1}{h}K\left(\frac{T_i-t}{h}\right)$. Then, for any fixed $t\in\mathcal{T}$, \\
		(i) $\mathbb{E}\left\{\Lambda_{N,1}(y,t,h)\right\} =  \bar{F}_{t}(y)f_{T}(t)\left\{1+O\left(\omega_{t}(y_{N},h)\log(y_N)+h\right)\right\}$ uniformly in $y\in[y_N,\Delta y_N]$, and $\mathbb{E}\left\{\Lambda_{N,2}(t,h)\right\} = f_{T}(t)+\frac{\kappa_{21}}{2}f^{\prime\prime}_{T}(t)h^2+O(h^3)$. \\
		(ii) $\frac{\mathbb{E}\left[\left\{\Lambda_{N,1}(y,t,h)\right\}^2\right]}{\left[\mathbb{E}\left\{\Lambda_{N,1}(y,t,h)\right\} \right]^2}-1= O\left(N^{-1}h^{-1}\left\{\bar{F}_{t}(y_N)\right\}^{-1}\right)$ uniformly in $y\in[y_N,\Delta y_N]$, and
		$\mathbb{E}\left[\Lambda_{N,2}(t,h)-\mathbb{E}\left\{\Lambda_{N,2}(t,h)\right\}\right]^2=O\left(N^{-1}h^{-1}\right)$. \\
		(iii) $\frac{\Lambda_{N,1}(y,t,h)}{\mathbb{E}\left\{\Lambda_{N,1}(y,t,h)\right\}}-1=O_{p}\left(N^{-1/2}h^{-1/2}\left\{\bar{F}_{t}(y_N)\right\}^{-1/2}\right)$ uniformly in $y\in[y_N,\Delta y_N]$, and $\frac{\Lambda_{N,2}(t,h)}{\mathbb{E}\left\{\Lambda_{N,2}(t,h)\right\}}-1=O_{p}\left(N^{-1/2}h^{-1/2}\right)$.
	\end{lemma}
	
	\begin{proof}
		(i) The results of this part  follow immediately from Lemma \ref{lemma:03} and the fact that $\{T_i,\boldsymbol{X}_i,Y_i\}_{i=1}^N$ are i.i.d..
		
		(ii) 
		Under Assumption \ref{assump:extreme_domain}, by Potter's inequality \citep[Theorem B.1.9]{deHaan2006}, for any $y\in[y_N,\Delta y_N]$,
		\begin{eqnarray}
			\notag & & O(1) = O\left(\Delta^{-1/\gamma(t)}\right) = \frac{\bar{F}_{t}(\Delta y_N)}{\bar{F}_{t}(y_N)} \leq \frac{\bar{F}_{t}(y)}{\bar{F}_{t}(y_N)} \leq 1.
		\end{eqnarray}
		Thus, $\bar{F}_{t}(y)\asymp\bar{F}_{t}(y_N)$ uniformly in $y\in[y_N,\Delta y_N]$, which implies $1/\bar{F}_{t}(y)\asymp 1/\bar{F}_{t}(y_N)$ uniformly in $y\in[y_N,\Delta y_N]$.

		We consider $\Lambda_{N,1}(y,t,h)$ first.  Simple algebra yields that
		\begin{eqnarray}
			\notag & & \mathbb{E}\left[\left\{\Lambda_{N,1}(y,t,h)\right\}^2\right] \\
			\notag & = & (1-N^{-1})\left[\mathbb{E}\left\{\Lambda_{N,1}(y,t,h)\right\}\right]^2 + N^{-1}\mathbb{E}\left[\left\{\pi_{0}(T,\mathbf{X})\right\}^2\mathbbm{1}(Y>y)\frac{1}{h^2}\left\{K\left(\frac{T-t}{h}\right)\right\}^2\right].
		\end{eqnarray}
		Under the assumption that $0<\eta_{1}\leq\pi_{0}(t,\boldsymbol{x})\leq\eta_{2}<\infty$ for all $(t,%
		\boldsymbol{x})\in\mathcal{T}\times\mathcal{X}$, straightforward calculations lead to
		\begin{eqnarray}
			\notag & & \mathbb{E}\left[\left\{\pi_{0}(T,\mathbf{X})\right\}^2\mathbbm{1}(Y>y)\frac{1}{h^2}\left\{K\left(\frac{T-t}{h}\right)\right\}^2\right] \\
			\notag & = & \frac{1}{h^2}\int\left\{K\left(\frac{t^{\prime}-t}{h}\right)\right\}^2\left\{\pi_{0}(t^{\prime},\mathbf{x})\right\}^2\mathbb{P}\left(Y(t^{\prime})>y|T=t^{\prime},\mathbf{X}=\mathbf{x}\right)f_{T|\mathbf{X}}(t^{\prime}|\mathbf{x})f_{\mathbf{X}}(\mathbf{x})d\mathbf{x}dt^{\prime} \\
			\notag & \leq & \frac{\eta_{2}}{h^2}\int\left\{K\left(\frac{t^{\prime}-t}{h}\right)\right\}^2\pi_{0}(t^{\prime},\mathbf{x})\mathbb{P}\left(Y(t^{\prime})>y|\mathbf{X}=\mathbf{x}\right)f_{T|\mathbf{X}}(t^{\prime}|\mathbf{x})f_{\mathbf{X}}(\mathbf{x})d\mathbf{x}dt^{\prime} \\
			\notag & = & \frac{\eta_{2}}{h^2}\int\left\{K\left(\frac{t^{\prime}-t}{h}\right)\right\}^2\bar{F}_{t^{\prime}}(y)f_{T}(t^{\prime})dt^{\prime} \\
			\notag & = & \frac{\eta_{2}}{h}\int_{\mathcal{K}}\left\{K\left(u\right)\right\}^2\bar{F}_{t+hu}(y)f_{T}(t+hu)du \\
			\notag & = & \frac{\eta_{2}}{h}\int_{\mathcal{K}}\left\{K\left(u\right)\right\}^2\bar{F}_{t}(y)f_{T}(t+hu)du + \frac{\eta_{2}}{h}\int_{\mathcal{K}}\left\{K\left(u\right)\right\}^2\left\{\bar{F}_{t+hu}(y)-\bar{F}_{t}(y)\right\}f_{T}(t+hu)du,
		\end{eqnarray}
		where the two terms in the last equation can be bounded by
		\begin{eqnarray}
			\notag \frac{\eta_{2}}{h}\int_{\mathcal{K}}\left\{K\left(u\right)\right\}^2\bar{F}_{t}(y)f_{T}(t+hu)du = O(h^{-1}\bar{F}_{t}(y_N))
		\end{eqnarray}
		and 
		\begin{eqnarray}
			\notag & & \left|\frac{\eta_{2}}{h}\int_{\mathcal{K}}\left\{K\left(u\right)\right\}^2\left\{\bar{F}_{t+hu}(y)-\bar{F}_{t}(y)\right\}f_{T}(t+hu)du\right| \\
			\notag & \leq & \frac{\eta_{2}}{h}\int_{\mathcal{K}}\left\{K\left(u\right)\right\}^2\left|\bar{F}_{t+hu}(y)-\bar{F}_{t}(y)\right|f_{T}(t+hu)du \\
			\notag & = & \frac{\eta_{2}\bar{F}_{t}(y)}{h}\int_{\mathcal{K}}\left\{K\left(u\right)\right\}^2\left|\frac{\bar{F}_{t+hu}(y)}{\bar{F}_{t}(y)}-1\right|f_{T}(t+hu)du \\
			\notag & = & O(h^{-1}\bar{F}_{t}(y_N)\omega_{t}(y_N, h)\log(y_N))
		\end{eqnarray}
		uniformly in $y\in[y_N,\Delta y_N]$. Therefore, 
		\begin{eqnarray}
			\notag & & \frac{\mathbb{E}\left[\left\{\Lambda_{N,1}(y,t,h)\right\}^2\right]}{\left[\mathbb{E}\left\{\Lambda_{N,1}(y,t,h)\right\} \right]^2}-1 \\
			\notag & = & -N^{-1} + N^{-1}\frac{O(h^{-1}\bar{F}_{t}(y_N)) + O(h^{-1}\bar{F}_{t}(y_N)\omega_{t}(y_N, h)\log(y_N)). }{\left[\bar{F}_{t}(y)f_{T}(t)\left\{1+O\left(\omega_{t}(y_{N},h)\log(y_N)+h\right)   \right\}\right]^2} \\
			\notag & = & O\left(N^{-1}h^{-1}\left\{\bar{F}_{t}(y_N)\right\}^{-1}\right).
		\end{eqnarray}
		
		Now we deal with $\Lambda_{N,2}(t,h)$. Similar approaches yield that
		\begin{eqnarray}
			\notag & & \mathbb{E}\left[\left\{\Lambda_{N,2}(t,h)\right\}^2\right] \\
			\notag & = & (1-N^{-1})\left[\mathbb{E}\left\{\Lambda_{N,2}(t,h)\right\}\right]^2 + N^{-1}\mathbb{E}\left[\left\{\pi_{0}(T,\mathbf{X})\right\}^2\frac{1}{h^2}\left\{K\left(\frac{T-t}{h}\right)\right\}^2\right],
		\end{eqnarray}
		where
		\begin{eqnarray}
			\notag & & \mathbb{E}\left[\left\{\pi_{0}(T,\mathbf{X})\right\}^2\frac{1}{h^2}\left\{K\left(\frac{T-t}{h}\right)\right\}^2\right] \leq \frac{\eta_{2}^2}{h^2}\mathbb{E}\left[\left\{K\left(\frac{T-t}{h}\right)\right\}^2\right] \\
			\notag & = & \frac{\eta_{2}^2}{h^2}\int\left\{K\left(\frac{t^{\prime}-t}{h}\right)\right\}^2f_{T}(t^{\prime})dt^{\prime} = \frac{\eta_{2}^2}{h}\int_{\mathcal{K}}\left\{K\left(u\right)\right\}^2f_{T}(t+hu)du \\
			\notag & = & O(h^{-1}).
		\end{eqnarray}
		Thus, we conclude that
		\begin{eqnarray}
			\notag & &  \frac{\mathbb{E}\left[\left\{\Lambda_{N,2}(t,h)\right\}^2\right]}{\left[\mathbb{E}\left\{\Lambda_{N,2}(t,h)\right\} \right]^2}-1 = -N^{-1} + O(N^{-1}h^{-1}) = O(N^{-1}h^{-1}).
		\end{eqnarray}
		
		(iii) The results of this part are direct consequences of part (i) and (ii).
	\end{proof}

	\begin{lemma}\label{lemma:05}
		Suppose Assumptions \ref{assump:extreme_domain}, \ref{UnconfoundAssump}, \ref{as:suppX}, \ref{assump:f_T}, \ref{asaump:K}, \ref{assump:pi} and \ref{asaump:y_N} hold.
		Then, for any fixed $t\in\mathcal{T}$, 
		\begin{eqnarray}
			\notag \sqrt{Nh\bar{F}_{t}(y)}\left\{\frac{\widehat{\bar{F}}_{t,h}(y)}{\bar{F}_{t}(y)}-1\right\} & = &  \sqrt{Nh\bar{F}_{t}(y)}\left\{\frac{\Lambda_{N,1}(y,t,h)}{\mathbb{E}\left\{\Lambda_{N,1}(y,t,h)\right\}}-1\right\} + o_{p}(1) \\
			\notag & = & \sqrt{Nh\bar{F}_{t}(y)}\left\{\frac{\Lambda_{N,1}(y,t,h)}{\bar{F}_{t}(y)f_{T}(t)}-1\right\} + o_{p}(1)
		\end{eqnarray}
		uniformly in $y\in[y_N, \Delta y_N]$, where $\Lambda_{N,1}(y,t,h) = N^{-1}\sum_{i=1}^{N}\pi_0(T_i,\mathbf{X}_i)\mathbbm{1}(Y_i>y)\frac{1}{h}K\left(\frac{T_i-t}{h}\right)$.
	\end{lemma}
	
	\begin{proof}
		Denote  $\Lambda_{N,2}(t,h) = N^{-1}\sum_{i=1}^{N}\pi_0(T_i,\mathbf{X}_i)\frac{1}{h}K\left(\frac{T_i-t}{h}\right)$ as in Lemma \ref{lemma:04}. In addition, define
		$$\Pi_{N,1}(y,t,h) = N^{-1}\sum_{i=1}^{N}\hat{\pi}(T_i,\mathbf{X}_i)\mathbbm{1}(Y_i>y)\frac{1}{h}K\left(\frac{T_i-t}{h}\right)$$ 
		and 
		$$\Pi_{N,2}(t,h) = N^{-1}\sum_{i=1}^{N}\hat{\pi}(T_i,\mathbf{X}_i)\frac{1}{h}K\left(\frac{T_i-t}{h}\right).$$ 
		Lemma \ref{lemma:04} implies that
		\begin{eqnarray}
			\notag & &  \Lambda_{N,1}(y,t,h) = \mathbb{E}\left\{\Lambda_{N,1}(y,t,h)\right\}\left\{1+O_{p}\left(N^{-1/2}h^{-1/2}\left\{\bar{F}_{t}(y_N)\right\}^{-1/2}\right)\right\} \\
			\notag & = & \bar{F}_{t}(y)f_{T}(t)\left\{1+O_{p}\left(\omega_{t}(y_{N},h)\log(y_N)\right) + O_{p}(h)+O_{p}\left(N^{-1/2}h^{-1/2}\left\{\bar{F}_{t}(y_N)\right\}^{-1/2}\right)\right\}
		\end{eqnarray}
		uniformly in $y\in[y_N, \Delta y_N]$, and
		\begin{eqnarray}
			\notag  \Lambda_{N,2}(t,h) & = & \mathbb{E}\left\{\Lambda_{N,2}(t,h)\right\}\left\{1+O_{p}\left(N^{-1/2}h^{-1/2}\right)\right\} \\
			\notag & = & \left\{f_{T}(t)+\frac{\kappa_{21}}{2}f^{\prime\prime}_{T}(t)h^2+O(h^3)\right\}\left\{1+O_{p}\left(N^{-1/2}h^{-1/2}\right)\right\} \\ 
			\notag & = & f_{T}(t)+O_{p}(h^2)+O_{p}\left(N^{-1/2}h^{-1/2}\right).
		\end{eqnarray}
		
		\textcolor{black}{Under Assumption \ref{assump:pi} that $\sup_{1\leq i\leq N}\left|\hat{\pi}(T_i,\mathbf{X}_i)-\pi_{0}(T_i,\mathbf{X}_i)\right|=O_{p}(\delta_N)$,} the difference between $\Pi_{N,1}(y,t,h)$ and $\Lambda_{N,1}(y,t,h)$ can be bounded by
		\begin{eqnarray}
			\notag & & \left|\frac{\Pi_{N,1}(y,t,h)}{\Lambda_{N,1}(y,t,h)}-1\right| \\
			\notag & \leq & \frac{1}{N\Lambda_{N,1}(y,t,h)}\sum_{i=1}^{N}\left|\hat{\pi}(T_i,\mathbf{X}_i)-\pi_{0}(T_i,\mathbf{X}_i)\right|\mathbbm{1}(Y_i>y)\frac{1}{h}K\left(\frac{T_i-t}{h}\right) \\
			\notag & \leq & \sup_{1\leq i\leq N}\left|\hat{\pi}(T_i,\mathbf{X}_i)-\pi_{0}(T_i,\mathbf{X}_i)\right|\frac{1}{N\Lambda_{N,1}(y,t,h)}\sum_{i=1}^{N}\mathbbm{1}(Y_i>y)\frac{1}{h}K\left(\frac{T_i-t}{h}\right) \\
			\notag & \leq & \frac{1}{\eta_{1}}\sup_{1\leq i\leq N}\left|\hat{\pi}(T_i,\mathbf{X}_i)-\pi_{0}(T_i,\mathbf{X}_i)\right|\frac{1}{N\Lambda_{N,1}(y,t,h)}\sum_{i=1}^{N}\pi_{0}(T_i,\mathbf{X}_i)\mathbbm{1}(Y_i>y)\frac{1}{h}K\left(\frac{T_i-t}{h}\right) \\
			\notag & = & O_{p}\left(\delta_N\right).
		\end{eqnarray}
		Similarly, we can bound $\Pi_{N,2}(t,h)/\Lambda_{N,2}(t,h)-1$ by
		\begin{eqnarray}
			\notag & & \left|\frac{\Pi_{N,2}(t,h)}{\Lambda_{N,2}(t,h)}-1\right| \\
			\notag & \leq & \frac{1}{\eta_{1}}\sup_{1\leq i\leq N}\left|\hat{\pi}(T_i,\mathbf{X}_i)-\pi_{0}(T_i,\mathbf{X}_i)\right|\frac{1}{N\Lambda_{N,2}(t,h)}\sum_{i=1}^{N}\pi_{0}(T_i,\mathbf{X}_i)\frac{1}{h}K\left(\frac{T_i-t}{h}\right) \\
			\notag & = & O_{p}\left(\delta_N\right).
		\end{eqnarray}
		
		Denote $\bar{F}_{t,h}(y)=\frac{\mathbb{E}\left\{\pi_0(T,\bs{X})\mathbbm{1}(Y> y)\frac{1}{h}K\left(\frac{T-t}{h}\right)\right\}}{\mathbb{E}\left\{\pi_{0}(T,\mathbf{X})\frac{1}{h}K\left(\frac{T-t}{h}\right)\right\}}=\frac{\mathbb{E}\left\{\Lambda_{N,1}(y,t,h)\right\}}{\mathbb{E}\left\{\Lambda_{N,2}(t,h)\right\}}$ as in Lemma \ref{lemma:03}, then, uniformly in $y\in[y_N,\Delta y_N]$,
		\begin{eqnarray}
			\notag & & \frac{\bar{F}_{t,h}(y)}{\bar{F}_{t}(y)} - 1 = O\left(\omega_{t}(y_{N},h)\log(y_N)+h\right).
		\end{eqnarray}
		Therefore, we decompose $\sqrt{Nh\bar{F}_{t}(y)}\left\{\frac{\widehat{\bar{F}}_{t,h}(y)}{\bar{F}_{t}(y)}-1\right\}$ as
		\begin{eqnarray}
			\notag & & \sqrt{Nh\bar{F}_{t}(y)}\left\{\frac{\widehat{\bar{F}}_{t,h}(y)}{\bar{F}_{t}(y)}-1\right\} = \sqrt{Nh\bar{F}_{t}(y)}\left\{\frac{\widehat{\bar{F}}_{t,h}(y)}{\bar{F}_{t,h}(y)}\frac{\bar{F}_{t,h}(y)}{\bar{F}_{t}(y)}-1\right\} \\
			\notag & = & \sqrt{Nh\bar{F}_{t}(y)}\left[\frac{\widehat{\bar{F}}_{t,h}(y)}{\bar{F}_{t,h}(y)}\left\{1+O\left(\omega_{t}(y_{N},h)\log(y_N)+h\right)\right\} - 1 \right] \\
			\notag & = & \sqrt{Nh\bar{F}_{t}(y)}\left\{\frac{\widehat{\bar{F}}_{t,h}(y)}{\bar{F}_{t,h}(y)}- 1 \right\} + \sqrt{Nh\bar{F}_{t}(y)}\frac{\widehat{\bar{F}}_{t,h}(y)}{\bar{F}_{t,h}(y)}\times O\left(\omega_{t}(y_{N},h)\log(y_N)+h\right).
		\end{eqnarray}
		Thus, we focus on $\frac{\widehat{\bar{F}}_{t,h}(y)}{\bar{F}_{t,h}(y)}$ first. Use the notations introduced above, uniformly in $y\in[y_N,\Delta y_N]$,
		\begin{eqnarray}
			\notag & & \frac{\widehat{\bar{F}}_{t,h}(y)}{\bar{F}_{t,h}(y)} = \frac{\Pi_{N,1}(y,t,h)}{\mathbb{E}\left\{\Lambda_{N,1}(y,t,h)\right\}}\frac{\mathbb{E}\left\{\Lambda_{N,2}(t,h)\right\}}{\Pi_{N,2}(t,h)} \\
			\notag & = & \frac{\Pi_{N,1}(y,t,h)}{\Lambda_{N,1}(y,t,h)}\frac{\Lambda_{N,1}(y,t,h)}{\mathbb{E}\left\{\Lambda_{N,1}(y,t,h)\right\}}\frac{\mathbb{E}\left\{\Lambda_{N,2}(t,h)\right\}}{\Lambda_{N,2}(t,h)}\frac{\Lambda_{N,2}(t,h)}{\Pi_{N,2}(t,h)} \\
			\notag & = & \left\{1+O_{p}(\delta_N)\right\}\left\{1+O_{p}\left(N^{-1/2}h^{-1/2}\left\{\bar{F}_{t}(y_N)\right\}^{-1/2}\right)\right\}\left\{1+O_{p}\left(N^{-1/2}h^{-1/2}\right)\right\}^{-1} \\
			\notag & = & 1 + O_{p}(\delta_N) + O_{p}\left(N^{-1/2}h^{-1/2}\left\{\bar{F}_{t}(y_N)\right\}^{-1/2}\right) + O_{p}\left(N^{-1/2}h^{-1/2}\right).
		\end{eqnarray}
		As $Nh\bar{F}_{t}(y_N)\to\infty$, $Nh^3\bar{F}_{t}(y_N)\to0$ and $\sqrt{Nh\bar{F}_{t}(y_N)}\omega_{t}(y_N,h)\log(y_N)\to0$, we get that
		$\sqrt{Nh\bar{F}_{t}(y)}\frac{\widehat{\bar{F}}_{t,h}(y)}{\bar{F}_{t,h}(y)}\times\left\{O\left(\omega_{t}(y_{N},h)\log(y_N)\right) + O(h)\right\}=o_{p}(1)$. Thus,
		\begin{equation}
			\sqrt{Nh\bar{F}_{t}(y)}\left\{\frac{\widehat{\bar{F}}_{t,h}(y)}{\bar{F}_{t}(y)}-1\right\} = \sqrt{Nh\bar{F}_{t}(y)}\left\{\frac{\widehat{\bar{F}}_{t,h}(y)}{\bar{F}_{t,h}(y)}- 1 \right\} + o_{p}(1)
		\end{equation}
		uniformly in $y\in[y_N,\Delta y_N]$. It remains to show the expansion of $\sqrt{Nh\bar{F}_{t}(y)}\left\{\frac{\widehat{\bar{F}}_{t,h}(y)}{\bar{F}_{t,h}(y)}- 1 \right\}$. It is noticed that
		\begin{eqnarray}
			\notag & &  \frac{\Pi_{N,2}(t,h)}{\mathbb{E}\left\{\Lambda_{N,2}(t,h)\right\}} = \frac{\Pi_{N,2}(t,h)}{\Lambda_{N,2}(t,h)}\frac{\Lambda_{N,2}(t,h)}{\mathbb{E}\left\{\Lambda_{N,2}(t,h)\right\}} \\
			\notag & = & \left\{1+O_{p}(\delta_N)\right\}\left\{1+O_{p}\left(N^{-1/2}h^{-1/2}\right)\right\} \\
			\notag & = & 1+O_{p}(\delta_N)+O_{p}\left(N^{-1/2}h^{-1/2}\right).
		\end{eqnarray}

		A more detailed decomposition of $\frac{\widehat{\bar{F}}_{t,h}(y)}{\bar{F}_{t,h}(y)}$ yields that
		\begin{eqnarray}
			\notag & & \frac{\widehat{\bar{F}}_{t,h}(y)}{\bar{F}_{t,h}(y)} = \frac{\Pi_{N,1}(y,t,h)}{\mathbb{E}\left\{\Lambda_{N,1}(y,t,h)\right\}}\frac{\mathbb{E}\left\{\Lambda_{N,2}(t,h)\right\}}{\Pi_{N,2}(t,h)}  \\
			%\notag & = & \frac{\Lambda_{N,1}(y,t,h) + \left\{\Pi_{N,1}(y,t,h)-\Lambda_{N,1}(y,t,h)\right\}}{\mathbb{E}\left\{\Lambda_{N,1}(y,t,h)\right\}}\frac{\mathbb{E}\left\{\Lambda_{N,2}(t,h)\right\}}{\Pi_{N,2}(t,h)}  \\
			\notag & = & \frac{\Lambda_{N,1}(y,t,h)}{\mathbb{E}\left\{\Lambda_{N,1}(y,t,h)\right\}}\left[1+\left\{\frac{\Pi_{N,1}(y,t,h)}{\Lambda_{N,1}(y,t,h)}-1\right\}\right]\frac{\mathbb{E}\left\{\Lambda_{N,2}(t,h)\right\}}{\Pi_{N,2}(t,h)} \\
			\notag & = & \frac{\Lambda_{N,1}(y,t,h)}{\mathbb{E}\left\{\Lambda_{N,1}(y,t,h)\right\}}\left\{1+O_{p}(\delta_N)\right\}\left\{1+O_{p}(\delta_N)+O_{p}\left(N^{-1/2}h^{-1/2}\right)\right\}^{-1} \\
			\notag & = & \frac{\Lambda_{N,1}(y,t,h)}{\mathbb{E}\left\{\Lambda_{N,1}(y,t,h)\right\}}\left\{1+O_{p}(\delta_N)+O_{p}\left(N^{-1/2}h^{-1/2}\right)\right\} \\
			\notag & = & \frac{\Lambda_{N,1}(y,t,h)}{\mathbb{E}\left\{\Lambda_{N,1}(y,t,h)\right\}} + \left\{1+O_{p}\left(N^{-1/2}h^{-1/2}\left\{\bar{F}_{t}(y_N)\right\}^{-1/2}\right)\right\}\left\{O_{p}(\delta_N)+O_{p}\left(N^{-1/2}h^{-1/2}\right)\right\} \\
			\notag & = & \frac{\Lambda_{N,1}(y,t,h)}{\mathbb{E}\left\{\Lambda_{N,1}(y,t,h)\right\}} + O_{p}(\delta_N)+O_{p}\left(N^{-1/2}h^{-1/2}\right).
		\end{eqnarray}
		Then, we arrive at
		\begin{eqnarray}
			\notag  \sqrt{Nh\bar{F}_{t}(y)}\left\{\frac{\widehat{\bar{F}}_{t,h}(y)}{\bar{F}_{t}(y)}-1\right\} & = & \sqrt{Nh\bar{F}_{t}(y)}\left\{\frac{\widehat{\bar{F}}_{t,h}(y)}{\bar{F}_{t,h}(y)}- 1 \right\} + o_{p}(1) \\
			\notag & = & \sqrt{Nh\bar{F}_{t}(y)}\left[\frac{\Lambda_{N,1}(y,t,h)}{\mathbb{E}\left\{\Lambda_{N,1}(y,t,h)\right\}}-1\right] \\
			\notag & & + \sqrt{Nh\bar{F}_{t}(y)}\times \left\{O_{p}(\delta_N)+O_{p}\left(N^{-1/2}h^{-1/2}\right)\right\} + o_{p}(1) \\
			\notag & = & \sqrt{Nh\bar{F}_{t}(y)}\left[\frac{\Lambda_{N,1}(y,t,h)}{\mathbb{E}\left\{\Lambda_{N,1}(y,t,h)\right\}}-1\right] + o_{p}(1)
		\end{eqnarray}
		uniformly in $y\in[y_N,\Delta y_N]$, where the last equation is due to the assumption that  $Nh\bar{F}_{t}(y_N)\delta_N^2\to 0$. Finally,
		\begin{eqnarray}
			\notag & &  \sqrt{Nh\bar{F}_{t}(y)}\left[\frac{\Lambda_{N,1}(y,t,h)}{\mathbb{E}\left\{\Lambda_{N,1}(y,t,h)\right\}}-1\right] \\
			\notag & = & \sqrt{Nh\bar{F}_{t}(y)}\left[\frac{\Lambda_{N,1}(y,t,h)}{\bar{F}_{t}(y)f_{T}(t)\left\{1+O\left(\omega_{t}(y_{N},h)\log(y_N)\right) + O(h)\right\}}-1\right] \\
			\notag & = &  \sqrt{Nh\bar{F}_{t}(y)}\left[\frac{\Lambda_{N,1}(y,t,h)}{\bar{F}_{t}(y)f_{T}(t)}\left\{1+O\left(\omega_{t}(y_{N},h)\log(y_N)\right) + O(h)\right\}-1\right] \\
			\notag & = & \sqrt{Nh\bar{F}_{t}(y)}\left\{\frac{\Lambda_{N,1}(y,t,h)}{\bar{F}_{t}(y)f_{T}(t)}-1\right\} \\
			\notag & & + \sqrt{Nh\bar{F}_{t}(y)}\frac{\Lambda_{N,1}(y,t,h)}{\bar{F}_{t}(y)f_{T}(t)}\left\{O\left(\omega_{t}(y_{N},h)\log(y_N)\right) + O(h)\right\} \\
			\notag & = & \sqrt{Nh\bar{F}_{t}(y)}\left\{\frac{\Lambda_{N,1}(y,t,h)}{\bar{F}_{t}(y)f_{T}(t)}-1\right\} + \frac{\Lambda_{N,1}(y,t,h)}{\mathbb{E}\left\{\Lambda_{N,1}(y,t,h)\right\}}\frac{\mathbb{E}\left\{\Lambda_{N,1}(y,t,h)\right\}}{\bar{F}_{t}(y)f_{T}(t)}\times o(1) \\
			\notag & = & \sqrt{Nh\bar{F}_{t}(y)}\left\{\frac{\Lambda_{N,1}(y,t,h)}{\bar{F}_{t}(y)f_{T}(t)}-1\right\} \\
			\notag & & + \left\{ 1 + O_{p}\left(N^{-1/2}h^{-1/2}\left\{\bar{F}_{t}(y_N)\right\}^{-1/2}\right)\right\}\left\{1+O_{p}\left(\omega_{t}(y_{N},h)\log(y_N)+h\right)\right\}\times o(1) \\
			\notag & = & \sqrt{Nh\bar{F}_{t}(y)}\left\{\frac{\Lambda_{N,1}(y,t,h)}{\bar{F}_{t}(y)f_{T}(t)}-1\right\} + o_{p}(1)
		\end{eqnarray}
		uniformly in $y\in[y_N,\Delta y_N]$. Thus, we finish the proof of this lemma.
	\end{proof}

	%Let $\Rightarrow$ denote week convergence of stochastic processes.

	The next lemma is from \citet{Pan2013SPL}, which gives the Drees-type inequality for second-order regular variation functions.
	
	\begin{lemma}\label{lemma:pan_spl}
		\cite[][Lemma 2.1]{Pan2013SPL} Let $h:\mathbb{R}_{+}\to\mathbb{R}$ be a measurable function that is eventually positive. $h$ is said to be of second-order regular variation (2RV) with the first order parameter $\nu\in\mathbb{R}$ and the second-order parameter $\rho\leq0$, denoted by $h\in 2RV_{\nu,\rho}$, if there exist an ultimately positive or negative function $A(y)$ with $\lim_{y\to\infty}A(y)=0$ such that
		\begin{eqnarray}\label{eq:appxlem_SM_02_01}
			\lim_{y\to\infty}\{h(zy)/h(y)-z^{\nu}\}/A(y)=z^{\nu}(z^{\rho}-1)/\rho
		\end{eqnarray}
		for all $z>0$. If $h\in2RV_{\nu,\rho}$ with auxiliary function $A(\cdot)$, $\nu\in\mathbb{R}$ and $\rho\leq0$, then for any $\epsilon,\delta>0$, there exists $y_{0}=y_{0}(\epsilon,\delta)>0$ such that for all $y\geq y_{0}$ and $zy\geq y_{0}$,
		\begin{eqnarray}\label{eq:appxlem_SM_02_02}
			\left|\frac{h(zy)/h(y)-z^{\nu}}{A(y)}-z^{\nu}\frac{z^{\rho}-1}{\rho}\right|\leq\epsilon z^{\nu}\left(\left|\frac{z^{\rho}-1}{\rho}\right|+z^{\rho}\max\{z^{\delta}, z^{-\delta}\}\right).
		\end{eqnarray}
	\end{lemma}
	
	\begin{lemma}\label{lemma:F_v}
		Suppose Assumption \ref{assump:second_order}  holds. (i) Let $y_N\to\infty$ as $N\to\infty$, then uniformly in $\upsilon\in[1,\Delta]$,
		\begin{eqnarray}
			\notag \left|\frac{\bar{F}_{t}(\upsilon y_N)}{\bar{F}_{t}(y_N)}-\upsilon^{-1/\gamma(t)}\right| = O(A_t(1/\bar{F}_{t}(y_N))).
		\end{eqnarray}
		
		(ii) Let $\alpha_N\to1$ as $N\to\infty$, then uniformly in $\upsilon\in[\Delta_0,1]$
		\begin{eqnarray}
			\notag \left|\frac{q_{t}(1-\upsilon(1-\alpha_N))}{q_{t}(\alpha_N)}-\upsilon^{-\gamma(t)}\right| = O\left(A_{t}(1/(1-\alpha_N))\right).
		\end{eqnarray}
	\end{lemma}
	
	\begin{proof}
		(i) First, \eqref{eq:second_order} in Assumption \ref{assump:second_order} is equivalent to 
		\begin{equation}\label{eq:second_order_1}
			\lim_{y\to\infty}\frac{\frac{\bar{F}_{t}(zy)}{\bar{F}_{t}(y)}-z^{-1/\gamma(t)}}{A_{t}(1/\bar{F}_{t}(y))/\{\gamma(t)\}^2} = z^{-1/\gamma(t)}\frac{z^{\rho(t)/\gamma(t)}-1}{\rho(t)/\gamma(t)} \text{~~for any~} z>0.
		\end{equation}
		This indicates that $\bar{F}_{t}(\cdot)\in 2RV_{-1/\gamma(t),\rho(t)/\gamma(t)}$ with auxiliary function $A_{t}(1/\bar{F}_{t}(\cdot))/\{\gamma(t)\}^2$. Then the result of (i) follows from Lemma \ref{lemma:pan_spl}.
		
		(ii) Assumption \ref{assump:second_order} and Theorem 2.3.9 of \citet{deHaan2006} indicate that $q_{t}(1-1/\cdot)\in 2RV_{\gamma(t),\rho(t)}$ with auxiliary function $A_{t}(\cdot)$. Then this part follows immediately from Lemma \ref{lemma:pan_spl}.
	\end{proof}

	\subsection{Proof of Lemma \ref{lemma:01}}
	
	\begin{proof}
	    We split the proof into two steps.
	    
	    \textbf{Step 1}. Here we show that for any fixed $t\in\mathcal{T}$,
		$$\mathcal{E}_{t,h}(\upsilon)\Rightarrow \Psi_{t}^{F}(\upsilon)$$
		in $\ell^{\infty}([1,\Delta])$, where $\Psi_{t}^{F}(\upsilon)$ is a centered Gaussian process with covariance function
		\begin{eqnarray}
			\notag \Omega_{t}^{F}(\upsilon_1,\upsilon_2) & = & 
			\kappa_{02} \{f_{T}(t)\}^{-1}\varpi_{t}^{F}(\upsilon_1,\upsilon_2).
		\end{eqnarray}

		According to Lemma \ref{lemma:05},
		\begin{eqnarray}
			\notag \sqrt{Nh\bar{F}_{t}(y)}\left\{\frac{\widehat{\bar{F}}_{t,h}(y)}{\bar{F}_{t}(y)}-1\right\} = \sqrt{Nh\bar{F}_{t}(y)}\left\{\frac{\Lambda_{N,1}(y,t,h)}{\bar{F}_{t}(y)f_{T}(t)}-1\right\} + o_{p}(1)
		\end{eqnarray}
		uniformly in $y\in[y_N, \Delta y_N]$, where $\Lambda_{N,1}(y,t,h) = N^{-1}\sum_{i=1}^{N}\pi_0(T_i,\mathbf{X}_i)\mathbbm{1}(Y_i>y)\frac{1}{h}K\left(\frac{T_i-t}{h}\right)$. Thus, for any fixed $t\in\mathcal{T}$,
		\begin{eqnarray}
			\notag \sup_{\upsilon\in[1,\Delta]}\left|\mathcal{E}_{t,h}(\upsilon)-\sum_{i=1}^{N}\psi_{t,h}(T_i, \mathbf{X}_i, Y_i, \upsilon, y_N)\right| = o_{p}(1),
		\end{eqnarray}
		where 
		%$\mathcal{E}_{t,h}(\upsilon) = \sqrt{Nh\bar{F}_{t}(\upsilon y_N)}\left\{\frac{\widehat{\bar{F}}_{t,h}(\upsilon y_N)}{\bar{F}_{t}(\upsilon y_N)}-1\right\}$ and 
		$\psi_{t,h}(T_i, \mathbf{X}_i, Y_i, \upsilon, y_N)=\sqrt{\frac{h\bar{F}_{t}(\upsilon y_N)}{N}}\left\{\frac{\pi_0(T_i,\mathbf{X}_i)\mathbbm{1}(Y_i>\upsilon y_N)\frac{1}{h}K\left(\frac{T_i-t}{h}\right)}{\bar{F}_{t}(\upsilon y_N)f_{T}(t)}-1\right\}
		$
		for $i=1,\ldots,N$.
		Denote for $\upsilon\in[1, \Delta]$,
		\begin{equation}
			\notag \mathcal{E}_{1, t,h}(\upsilon) = \sqrt{Nh\bar{F}_{t}(\upsilon y_N)}\left\{\frac{\Lambda_{N,1}(\upsilon y_N,t,h)}{\bar{F}_{t}(\upsilon y_N)f_{T}(t)}-1\right\} = \sum_{i=1}^{N}\psi_{t,h}(T_i, \mathbf{X}_i, Y_i, \upsilon, y_N),
		\end{equation}
		it suffices to show that $\mathcal{E}_{1, t,h}(\upsilon)\Rightarrow \Psi_{t}^{F}(\upsilon)$ in $\ell^{\infty}([1,\Delta])$.
		
		We first show that for any fixed $\upsilon\in[1,\Delta]$, $\mathcal{E}_{1, t,h}(\upsilon)$ converges in distribution to $N(0,\Omega_{t}^{F}(\upsilon, \upsilon))$. 
		It is clear that $\{\psi_{t,h}(T_i, \mathbf{X}_i, Y_i, \upsilon, y_N)\}_{i=1}^{N}$ is a sequence of i.i.d.~random variables with zero mean and variance
		\begin{eqnarray}
			\notag & & \mathbb{E}\left[\left\{\psi_{t,h}(T_i, \mathbf{X}_i, Y_i, \upsilon, y_N)\right\}^2\right] \\
			\notag & = & \frac{h\bar{F}_{t}(\upsilon y_N)}{N}\left\{\frac{\mathbb{E}\left[\left\{\pi_{0}(T,\mathbf{X})\right\}^2\mathbbm{1}(Y>\upsilon y_N)\frac{1}{h^2}\left\{K\left(\frac{T-t}{h}\right)\right\}^2\right]}{\left\{\bar{F}_{t}(\upsilon y_N)f_{T}(t)\right\}^2}+1-\frac{2\mathbb{E}\left\{\Lambda_{N,1}(\upsilon y_N, t, h)\right\}}{\bar{F}_{t}(\upsilon y_N)f_{T}(t)}\right\} \\
			\notag &= & \frac{h\bar{F}_{t}(\upsilon y_N)}{N}\left\{\frac{\kappa_{02}h^{-1}f_{T}(t)\mathbb{E}\left[\left\{\pi_{0}(T,\mathbf{X})\right\}^2\mathbbm{1}(y>\upsilon y_N)|T=t\right]\{1+o(1)\}}{\left\{\bar{F}_{t}(\upsilon y_N)f_{T}(t)\right\}^2}-1\right\} \\
			\notag & & + N^{-1}h\bar{F}_{t}(\upsilon y_N)\left\{1+O\left(\omega_{t}(y_{N},h)\log(y_N)\right) + O(h)\right\} \\
			\notag & = & \kappa_{02}N^{-1}\{f_{T}(t)\}^{-1}\{\bar{F}_{t}(\upsilon y_N)\}^{-1}\mathbb{E}\left[\left\{\pi_{0}(T,\mathbf{X})\right\}^2\mathbbm{1}(y>\upsilon y_N)|T=t\right]\{1+o(1)\} + o(1),
		\end{eqnarray}
		where the second equation is due to Lemma \ref{lemma:03} and 
		\begin{eqnarray}
			\notag & & \mathbb{E}\left[\left\{\pi_{0}(T,\mathbf{X})\right\}^2\mathbbm{1}(Y>\upsilon y_N)\frac{1}{h^2}\left\{K\left(\frac{T-t}{h}\right)\right\}^2\right] \\
			\notag & = & \int \left\{\pi_{0}(t^{\prime},\mathbf{x})\right\}^2\mathbbm{1}(y>\upsilon y_N)\frac{1}{h^2}\left\{K\left(\frac{t^{\prime}-t}{h}\right)\right\}^2 f_{Y,T,\mathbf{X}}(y,t^{\prime},\mathbf{x})dydt^{\prime}d\mathbf{x} \\
			\notag & = & \frac{1}{h}\int \left\{\pi_{0}(t+hu,\mathbf{x})\right\}^2\mathbbm{1}(y>\upsilon y_N)\left\{K\left(u\right)\right\}^2 f_{Y,T,\mathbf{X}}(y,t+hu,\mathbf{x})dydud\mathbf{x} \\
			\notag & = & \left[\frac{1}{h}\int \left\{\pi_{0}(t,\mathbf{x})\right\}^2\mathbbm{1}(y>\upsilon y_N)\left\{K\left(u\right)\right\}^2 f_{Y,T,\mathbf{X}}(y,t,\mathbf{x})dydud\mathbf{x}\right]\{1+o(1)\} %\textcolor{red}{\text{~~more details needed here}} 
			\\
			\notag & = & \kappa_{02}h^{-1}f_{T}(t)\mathbb{E}\left[\left\{\pi_{0}(T,\mathbf{X})\right\}^2\mathbbm{1}(Y>\upsilon y_N)|T=t\right]\{1+o(1)\}.
		\end{eqnarray}
		It is worth mentioning that $\mathbb{E}\left[\left\{\pi_{0}(T,\mathbf{X})\right\}^2\mathbbm{1}(Y>\upsilon y_N)|T=t\right] \asymp \bar{F}_{t}(y_N)$.
		By Jensen's inequality, for any positive $\delta<1$,
		\begin{eqnarray}
			\notag & & \sum_{i=1}^{N}\mathbb{E}\left|\psi_{t,h}(T_i, \mathbf{X}_i, Y_i, \upsilon, y_N)\right|^{2+\delta} = N\mathbb{E}\left|\psi_{t,h}(T_i, \mathbf{X}_i, Y_i, \upsilon, y_N)\right|^{2+\delta} \\
			\notag & = & \frac{\left\{h\bar{F}_{t}(\upsilon y_N)\right\}^{1+\delta/2}}{N^{\delta/2}}\mathbb{E}\left|\frac{\pi_0(T_i,\mathbf{X}_i)\mathbbm{1}(Y_i>\upsilon y_N)\frac{1}{h}K\left(\frac{T_i-t}{h}\right)}{\bar{F}_{t}(\upsilon y_N)f_{T}(t)}-1\right|^{2+\delta} \\
			\notag & \leq & \frac{2^{1+\delta}\left\{h\bar{F}_{t}(\upsilon y_N)\right\}^{1+\delta/2}}{N^{\delta/2}}\left\{\mathbb{E}\left|\frac{\pi_0(T_i,\mathbf{X}_i)\mathbbm{1}(Y_i>\upsilon y_N)\frac{1}{h}K\left(\frac{T_i-t}{h}\right)}{\bar{F}_{t}(\upsilon y_N)f_{T}(t)}\right|^{2+\delta}+1\right\} \\
			\notag & \leq & \frac{2^{1+\delta}\left\{h\bar{F}_{t}(\upsilon y_N)\right\}^{1+\delta/2}}{N^{\delta/2}}\left\{\frac{\eta_{2}^{1+\delta}}{h^{2+\delta}}\frac{\mathbb{E}\left[\pi_0(T_i,\mathbf{X}_i)\mathbbm{1}(Y_i>\upsilon y_N)\left\{K\left(\frac{T_i-t}{h}\right)\right\}^{2+\delta}\right]}{\left\{\bar{F}_{t}(\upsilon y_N)f_{T}(t)\right\}^{2+\delta}}+1\right\} ,
		\end{eqnarray}
		where 
		\begin{eqnarray}
			\notag & & \mathbb{E}\left[\pi_0(T_i,\mathbf{X}_i)\mathbbm{1}(Y_i>\upsilon y_N)\left\{K\left(\frac{T_i-t}{h}\right)\right\}^{2+\delta}\right] \\
			\notag & = & \int\left\{K\left(\frac{t^{\prime}-t}{h}\right)\right\}^{2+\delta}\pi_{0}(t^{\prime},\mathbf{x})\mathbb{P}\left(Y(t^{\prime})>\upsilon y_N|T=t^{\prime},\mathbf{X}=\mathbf{x}\right)f_{T|\mathbf{X}}(t^{\prime}|\mathbf{x})f_{\mathbf{X}}(\mathbf{x})d\mathbf{x}dt^{\prime} \\
			\notag & = & \int\left\{K\left(\frac{t^{\prime}-t}{h}\right)\right\}^{2+\delta}\pi_{0}(t^{\prime},\mathbf{x})\mathbb{P}\left(Y(t^{\prime})>\upsilon y_N|\mathbf{X}=\mathbf{x}\right)f_{T|\mathbf{X}}(t^{\prime}|\mathbf{x})f_{\mathbf{X}}(\mathbf{x})d\mathbf{x}dt^{\prime} \\
			\notag & = & \int\left\{K\left(\frac{t^{\prime}-t}{h}\right)\right\}^{2+\delta}\bar{F}_{t^{\prime}}(y)f_{T}(t^{\prime})dt^{\prime} \\
			\notag & = & h\int_{\mathcal{K}}\left\{K\left(u\right)\right\}^{2+\delta}\bar{F}_{t+hu}(\upsilon y_N)f_{T}(t+hu)du \\
			\notag & = & h\int_{\mathcal{K}}\left\{K\left(u\right)\right\}^{2+\delta}\bar{F}_{t}(\upsilon y_N)f_{T}(t+hu)du \\
			\notag & & + h\int_{\mathcal{K}}\left\{K\left(u\right)\right\}^{2+\delta}\left\{\bar{F}_{t+hu}(\upsilon y_N)-\bar{F}_{t}(\upsilon y_N)\right\}f_{T}(t+hu)du \\
			\notag & = & O\left(h\bar{F}_{t}(y_N)\right).
		\end{eqnarray}
		Hence, we have
		\begin{eqnarray}
			\notag & & \sum_{i=1}^{N}\mathbb{E}\left|\psi_{t,h}(T_i, \mathbf{X}_i, Y_i, \upsilon, y_N)\right|^{2+\delta} \\
			\notag & = & O\left(\frac{2^{1+\delta}\left\{h\bar{F}_{t}(y_N)\right\}^{1+\delta/2}}{N^{\delta/2}}\left\{\frac{\eta_{2}^{1+\delta}}{h^{2+\delta}}\frac{O\left(h\bar{F}_{t}(y_N)\right)}{O\left(\left\{\bar{F}_{t}(y_N)\right\}^{2+\delta}\right)}+1\right\}\right) \\
			\notag & = & O\left(N^{-\delta/2}h^{-\delta/2}\left\{\bar{F}_{t}(y_N)\right\}^{-\delta/2}\right) \to 0.
		\end{eqnarray}
		Therefore, the Lyapunov Central Limit Theorem indicates that $\mathcal{E}_{1, t,h}(\upsilon)$ converges in distribution to $N(0,\Omega_{t}^{F}(\upsilon, \upsilon))$ for any fixed $\upsilon\in[1,\Delta]$. 
		
		Next, we check the stochastic equicontinuity of $\{\mathcal{E}_{1, t,h}(\upsilon): \upsilon\in[1,\Delta]\}$ to show the week convergence of $\mathcal{E}_{t,h}(\upsilon)$. 
		Recall $\mathcal{E}_{1, t,h}(\upsilon)=N^{-1/2}\sum_{i=1}^{N}N^{1/2}\psi_{t,h}(T_i, \mathbf{X}_i, Y_i, \upsilon, y_N)$ with $N^{1/2}\psi_{t,h}(T_i, \mathbf{X}_i, Y_i, \upsilon, y_N)=\sqrt{h\bar{F}_{t}(\upsilon y_N)}\left\{\frac{\pi_0(T_i,\mathbf{X}_i)\mathbbm{1}(Y_i>\upsilon y_N)\frac{1}{h}K\left(\frac{T_i-t}{h}\right)}{\bar{F}_{t}(\upsilon y_N)f_{T}(t)}-1\right\}$.
		According to Theorem 5 of \citet{andrews1994empirical}, it suffices to show that the Ossiander's $L^2$ entropy condition is satisfies, that is,
		\begin{equation}
			\int_{0}^{1}\sqrt{\log N_{[]}(\varepsilon, \mathcal{F}, L^2(P))}d\varepsilon<\infty,
		\end{equation}
		where $\mathcal{F}=\{N^{1/2}\psi_{t,h}(T_i, \mathbf{X}_i, Y_i, \upsilon, y_N): \upsilon\in[1,\Delta]\}$ and $N_{[]}(\varepsilon, \mathcal{F}, L^2(P))$ is the $L^2$ bracketing $\varepsilon$-entropy of $\mathcal{F}$. For that purpose, we show the $L^2$-continuity condition of $N^{1/2}\psi_{t,h}(T_i, \mathbf{X}_i, Y_i, \upsilon, y_N)$ \citep[Theorem 5]{andrews1994empirical}.
		
		It is noticed that
		\begin{eqnarray}
			\notag & & N^{1/2}\psi_{t,h}(T_i, \mathbf{X}_i, Y_i, \upsilon_1, y_N)-N^{1/2}\psi_{t,h}(T_i, \mathbf{X}_i, Y_i, \upsilon_2, y_N) \\
			\notag & = & \sqrt{h\bar{F}_{t}(\upsilon_1 y_N)}\pi_{0}(T_i,\mathbf{X}_i)\frac{1}{h}K\left(\frac{T_i-t}{h}\right)\left\{\frac{\mathbbm{1}(Y_i>\upsilon_1 y_N)}{\bar{F}_{t}(\upsilon_1 y_N)f_{T}(t)}-\frac{\mathbbm{1}(Y_i>\upsilon_2 y_N)}{\bar{F}_{t}(\upsilon_2 y_N)f_{T}(t)}\right\}\\
			%\notag & & + \frac{\sqrt{Nh\bar{F}_{t}(\upsilon_1 y_N)}}{f_{T}(t)}\left\{\frac{1}{\bar{F}_{t}(\upsilon_2 y_N)}-\frac{1}{\bar{F}_{t}(\upsilon_1 y_N)}\right\}\frac{1}{N}\sum_{i=1}^{N}\pi_{0}(T_i,\mathbf{X}_i)\frac{1}{h}K\left(\frac{T_i-t}{h}\right)\mathbbm{1}(Y_i\geq \upsilon_2 y_N)\\
			\notag & & + \left\{\sqrt{\frac{\bar{F}_{t}(\upsilon_1 y_N)}{\bar{F}_{t}(\upsilon_2 y_N)}}-1\right\}\sqrt{h\bar{F}_{t}(\upsilon_2 y_N)}\left\{\frac{\pi_{0}(T_i,\mathbf{X}_i)\mathbbm{1}(Y_i>\upsilon_2 y_N)\frac{1}{h}K\left(\frac{T_i-t}{h}\right)}{\bar{F}_{t}(\upsilon_2 y_N)f_{T}(t)}-1\right\} \\
			%\notag & = & \frac{\sqrt{Nh\bar{F}_{t}(\upsilon_1 y_N)}}{f_{T}(t)}\frac{1}{N}\sum_{i=1}^{N}\pi_{0}(T_i,\mathbf{X}_i)\frac{1}{h}K\left(\frac{T_i-t}{h}\right)\left\{\frac{\mathbbm{1}(Y_i>\upsilon_1 y_N)}{\mathbb{E}\left\{\Lambda_{N,1}(\upsilon_1 y_N,t,h)\right\}}-\frac{\mathbbm{1}(Y_i>\upsilon_2 y_N)}{\mathbb{E}\left\{\Lambda_{N,1}(\upsilon_2 y_N,t,h)\right\}}\right\}\\ 
			%\notag & & + O_{p}\left(\left|\left(\frac{\upsilon_1}{\upsilon_2}\right)^{-1/\{2\gamma(t)\}}-1\right|\right) + o_{p}(1).
			\notag & = & \frac{\sqrt{h\bar{F}_{t}(\upsilon_1 y_N)}}{\bar{F}_{t}(\upsilon_2 y_N)f_{T}(t)}\pi_{0}(T_i,\mathbf{X}_i)\frac{1}{h}K\left(\frac{T_i-t}{h}\right) \left\{\mathbbm{1}(Y_i\geq \upsilon_1 y_N)-\mathbbm{1}(Y_i\geq \upsilon_2 y_N)\right\} \\
			\notag & & + \sqrt{h\bar{F}_{t}(\upsilon_1 y_N)}\left\{\frac{1}{\bar{F}_{t}(\upsilon_1 y_N)f_{T}(t)}-\frac{1}{\bar{F}_{t}(\upsilon_2 y_N)f_{T}(t)}\right\}\pi_{0}(T_i,\mathbf{X}_i)\frac{1}{h}K\left(\frac{T_i-t}{h}\right)  \\
			\notag & & ~~~~~~~~~~~~~~~~~~~~~~~~~~~~~~~~~~~~~~~~~~~~~~~~~~~~~~~~~~~~~~~~~~~~~~~~~~~~~~~~~~\times \mathbbm{1}(Y_i>\upsilon_1 y_N) \\
			\notag & & + \left\{\sqrt{\frac{\bar{F}_{t}(\upsilon_1 y_N)}{\bar{F}_{t}(\upsilon_2 y_N)}}-1\right\}\sqrt{h\bar{F}_{t}(\upsilon_2 y_N)}\left\{\frac{\pi_{0}(T_i,\mathbf{X}_i)\mathbbm{1}(Y_i>\upsilon_2 y_N)\frac{1}{h}K\left(\frac{T_i-t}{h}\right)}{\bar{F}_{t}(\upsilon_2 y_N)f_{T}(t)}-1\right\},
		\end{eqnarray}
		where by Potter's inequality \citep[Theorem B.1.9]{deHaan2006}, for any $\delta_1, \delta_2>0$ and sufficient large $y_N$, we have
		\begin{eqnarray}\label{eq:potter}
			\notag & & (1-\delta_1) \left(\frac{\upsilon_1}{\upsilon_2}\right)^{-1/\gamma(t)}\min\left(\left(\frac{\upsilon_1}{\upsilon_2}\right)^{\delta_2}, \left(\frac{\upsilon_1}{\upsilon_2}\right)^{-\delta_2}\right)  \leq  \frac{\bar{F}_{t}(\upsilon_1 y_N)}{\bar{F}_{t}(\upsilon_2 y_N)}\\
			& & ~~~~~~~~~~~~~~~~~~~~~~~~~~~~~~~~ \leq (1+\delta_1) \left(\frac{\upsilon_1}{\upsilon_2}\right)^{-1/\gamma(t)}\max\left(\left(\frac{\upsilon_1}{\upsilon_2}\right)^{\delta_2}, \left(\frac{\upsilon_1}{\upsilon_2}\right)^{-\delta_2}\right).
		\end{eqnarray}
		\iffalse
		$\delta_1, \delta_2>0$ and sufficient large $y_N$, we have
		\begin{eqnarray}
			\notag & & (1-\delta_1) \left(\frac{\upsilon_1}{\upsilon_2}\right)^{-1/\gamma(t)}\min\left(\left(\frac{\upsilon_1}{\upsilon_2}\right)^{\delta_2}, \left(\frac{\upsilon_1}{\upsilon_2}\right)^{-\delta_2}\right)  \leq  \frac{\bar{F}_{t}(\upsilon_1 y_N)}{\bar{F}_{t}(\upsilon_2 y_N)}\\
			& & ~~~~~~~~~~~~~~~~~~~~~~~~~~~~~~~~ \leq (1+\delta_1) \left(\frac{\upsilon_1}{\upsilon_2}\right)^{-1/\gamma(t)}\max\left(\left(\frac{\upsilon_1}{\upsilon_2}\right)^{\delta_2}, \left(\frac{\upsilon_1}{\upsilon_2}\right)^{-\delta_2}\right).
		\end{eqnarray}
		
		and
		\begin{eqnarray}
			\notag & & \bar{F}_{t}(\upsilon_1 y_N)f_{T}(t)\left\{\frac{1}{\bar{F}_{t}(\upsilon_1 y_N)f_{T}(t)}-\frac{1}{\bar{F}_{t}(\upsilon_2 y_N)f_{T}(t)}\right\} \\
			\notag & = & 1- \frac{\bar{F}_{t}(\upsilon_{1}y_N)}{\bar{F}_{t}(\upsilon_{2}y_N)}\\
			\notag & = & \left\{1-\frac{\bar{F}_{t}(\upsilon_1 y_N)}{\bar{F}_{t}(\upsilon_2 y_N)}\right\}\left\{1+O\left(\omega_{t}(y_{N},h)\log(y_N)\right) + O(h)\right\} \\
			\notag & & + \frac{\bar{F}_{t}(\upsilon_1 y_N)}{\bar{F}_{t}(\upsilon_2 y_N)}\left\{O\left(\omega_{t}(y_{N},h)\log(y_N)\right) + O(h)\right\} -\left\{O\left(\omega_{t}(y_{N},h)\log(y_N)\right) + O(h)\right\}. % \\
			%\notag & \to & \left(1+\frac{\upsilon_1-\upsilon_2}{\upsilon_2}\right)^{-1/\gamma(t)}-1
		\end{eqnarray}
		\fi
		%as $y_N\to\infty$.
		
		Therefore, for sufficient large $y_N$,
		\begin{eqnarray}
			\notag & & \mathbb{E}\left[\sup_{\upsilon_1: |\upsilon_1-\upsilon|\leq \delta}\left|N^{1/2}\psi_{t,h}(T_i, \mathbf{X}_i, Y_i, \upsilon_1, y_N)-N^{1/2}\psi_{t,h}(T_i, \mathbf{X}_i, Y_i, \upsilon, y_N)\right|^2\right] \\
			\notag & \leq & 2\mathbb{E}\left[\sup_{\upsilon_1: |\upsilon_1-\upsilon|\leq \delta}\frac{h\bar{F}_{t}(\upsilon_1 y_N)}{\left\{\bar{F}_{t}(\upsilon y_N)f_{T}(t)\right\}^2}\left\{\pi_{0}(T_i,\mathbf{X}_i)\right\}^2\frac{1}{h^2}\left\{K\left(\frac{T_i-t}{h}\right)\right\}^2\right. \\
			\notag & & ~~~~~~~~~~~~~~~~~~~~~~~~~~~~~~~~~~~~~~~~~~~~~~~~~~~~~~~~~~~~ \times \left. \mathbbm{1}((\upsilon_1 \wedge \upsilon) y_N\leq Y_i\leq (\upsilon_1 \vee \upsilon)y_N)\right] \\
			\notag & & + 2\mathbb{E}\left[\sup_{\upsilon_1: |\upsilon_1-\upsilon|\leq \delta}h\bar{F}_{t}(\upsilon_1 y_N)\left\{\frac{1}{\bar{F}_{t}(\upsilon_1 y_N)f_{T}(t)}-\frac{1}{\bar{F}_{t}(\upsilon y_N)f_{T}(t)}\right\}^2\right. \\
			\notag & & ~~~~~~~~~~~~~~~~~~~~~~~~~~~~~~~~~~~~~~~~~\times \left. \left\{\pi_{0}(T_i,\mathbf{X}_i)\right\}^2\frac{1}{h^2}\left\{K\left(\frac{T_i-t}{h}\right)\right\}^2 \mathbbm{1}(Y_i>\upsilon_1 y_N)\right] \\
			\notag & & + 2\mathbb{E}\left[\sup_{\upsilon_1: |\upsilon_1-\upsilon|\leq \delta}\left\{\sqrt{\frac{\bar{F}_{t}(\upsilon_1 y_N)}{\bar{F}_{t}(\upsilon y_N)}}-1\right\}^2h\bar{F}_{t}(\upsilon y_N)\left\{\frac{\pi_{0}(T_i,\mathbf{X}_i)\mathbbm{1}(Y_i>\upsilon y_N)\frac{1}{h}K\left(\frac{T_i-t}{h}\right)}{\bar{F}_{t}(\upsilon y_N)f_{T}(t)}-1\right\}^2\right] \\
			\notag & \leq & 2\mathbb{E}\left[\frac{\eta_{2}\eta_{3}^{-2}\bar{F}_{t}( y_N)}{\left\{\bar{F}_{t}(\Delta y_N)\right\}^2}\sup_{\upsilon_1: |\upsilon_1-\upsilon|\leq \delta}\pi_{0}(T_i,\mathbf{X}_i)\frac{1}{h}\left\{K\left(\frac{T_i-t}{h}\right)\right\}^2  \mathbbm{1}((\upsilon_1 \wedge \upsilon) y_N\leq Y_i\leq (\upsilon_1 \vee \upsilon)y_N)\right] \\
			\notag & & + 2\mathbb{E}\left[\frac{\eta_{2}\eta_{3}^{-2}\bar{F}_{t}(y_N)}{\left\{\bar{F}_{t}(\Delta y_N)\right\}^2}\sup_{\upsilon_1: |\upsilon_1-\upsilon|\leq \delta}\left\{\frac{\bar{F}_{t}(\upsilon_1 y_N)}{\bar{F}_{t}(\upsilon y_N)}-1\right\}^2  \pi_{0}(T_i,\mathbf{X}_i)\frac{1}{h}\left\{K\left(\frac{T_i-t}{h}\right)\right\}^2 \mathbbm{1}(Y_i>\upsilon_1 y_N)\right] \\
			\notag & & + 2h\bar{F}_{t}(y_N)\mathbb{E}\left[\sup_{\upsilon_1: |\upsilon_1-\upsilon|\leq \delta}\left\{\sqrt{\frac{\bar{F}_{t}(\upsilon_1 y_N)}{\bar{F}_{t}(\upsilon y_N)}}-1\right\}^2\left\{\frac{\pi_{0}(T_i,\mathbf{X}_i)\mathbbm{1}(Y_i>\upsilon y_N)\frac{1}{h}K\left(\frac{T_i-t}{h}\right)}{\bar{F}_{t}(\upsilon y_N)f_{T}(t)}-1\right\}^2\right] \\
			\notag & \leq & \frac{2\eta_{2}\eta_{3}^{-2}\bar{F}_{t}( y_N)}{\left\{\bar{F}_{t}(\Delta y_N)\right\}^2}\mathbb{E}\left[ \pi_{0}(T_i,\mathbf{X}_i)\frac{1}{h}\left\{K\left(\frac{T_i-t}{h}\right)\right\}^2 \mathbbm{1}((1 \vee (\upsilon-\delta)) y_N\leq Y_i\leq (\Delta \wedge (\upsilon+\delta)) y_N)\right] \\
			\notag & & + \frac{2\eta_{2}\eta_{3}^{-2}\bar{F}_{t}(y_N)}{\left\{\bar{F}_{t}(\Delta y_N)\right\}^2}\sup_{\upsilon_1: |\upsilon_1-\upsilon|\leq \delta}\left\{\frac{\bar{F}_{t}(\upsilon_1 y_N)}{\bar{F}_{t}(\upsilon y_N)}-1\right\}^2\mathbb{E}\left[ \pi_{0}(T_i,\mathbf{X}_i)\frac{1}{h}\left\{K\left(\frac{T_i-t}{h}\right)\right\}^2 \mathbbm{1}(Y_i>\upsilon y_N)\right]  \\
			\notag & & + 2h\bar{F}_{t}(y_N)\sup_{\upsilon_1: |\upsilon_1-\upsilon|\leq \delta}\left\{\sqrt{\frac{\bar{F}_{t}(\upsilon_1 y_N)}{\bar{F}_{t}(\upsilon y_N)}}-1\right\}^2\mathbb{E}\left[\left\{\frac{\pi_{0}(T_i,\mathbf{X}_i)\mathbbm{1}(Y_i>\upsilon y_N)\frac{1}{h}K\left(\frac{T_i-t}{h}\right)}{\bar{F}_{t}(\upsilon y_N)f_{T}(t)}-1\right\}^2\right] \\
			\notag & \leq & \frac{2C_1\eta_{2}\eta_{3}^{-2}\bar{F}_{t}( y_N)}{\left\{\bar{F}_{t}(\Delta y_N)\right\}^2}\mathbb{E}\left[ \pi_{0}(T_i,\mathbf{X}_i)\frac{1}{h}K\left(\frac{T_i-t}{h}\right) \mathbbm{1}((1 \vee (\upsilon-\delta)) y_N\leq Y_i\leq (\Delta \wedge (\upsilon+\delta)) y_N)\right] \\
			\notag & & + \frac{2\eta_{2}\eta_{3}^{-2}\bar{F}_{t}(y_N)}{\left\{\bar{F}_{t}(\Delta y_N)\right\}^2}\sup_{\upsilon_1: |\upsilon_1-\upsilon|\leq \delta}\left\{\frac{\bar{F}_{t}(\upsilon_1 y_N)}{\bar{F}_{t}(\upsilon y_N)}-1\right\}^2\times O\left(\bar{F}_{t}(y_N)\right) \\
			\notag & & + 2h\bar{F}_{t}(y_N)\sup_{\upsilon_1: |\upsilon_1-\upsilon|\leq \delta}\left\{\sqrt{\frac{\bar{F}_{t}(\upsilon_1 y_N)}{\bar{F}_{t}(\upsilon y_N)}}-1\right\}^2\times O\left(h^{-1}\left\{\bar{F}_{t}(y_N)\right\}^{-1}\right) \\
			\notag & \leq & 2C_2\frac{\eta_{2}\eta_{3}^{-2}\bar{F}_{t}( y_N)\bar{F}_{t}((1 \vee (\upsilon-\delta)) y_N)}{\left\{\bar{F}_{t}(\Delta y_N)\right\}^2}\left\{\frac{\bar{F}_{t}((\Delta \wedge (\upsilon+\delta)) y_N)}{\bar{F}_{t}((1 \vee (\upsilon-\delta)) y_N)}-1\right\} \\
			\notag & & + \frac{2C_2\eta_{2}\eta_{3}^{-2}\{\bar{F}_{t}(y_N)\}^2}{\left\{\bar{F}_{t}(\Delta y_N)\right\}^2}\sup_{\upsilon_1: |\upsilon_1-\upsilon|\leq \delta}\left\{\frac{\bar{F}_{t}(\upsilon_1 y_N)}{\bar{F}_{t}(\upsilon y_N)}-1\right\}^2  \\
			\notag & & + C_2\sup_{\upsilon_1: |\upsilon_1-\upsilon|\leq \delta}\left\{\sqrt{\frac{\bar{F}_{t}(\upsilon_1 y_N)}{\bar{F}_{t}(\upsilon y_N)}}-1\right\}^2,
		\end{eqnarray}
		where $C_1$ and $C_2$ are positive constants. Thus, regrading \eqref{eq:potter}, it is clear that
		\begin{eqnarray}
			\notag \mathbb{E}\left[\sup_{\upsilon_1: |\upsilon_1-\upsilon|\leq \delta}\left|N^{1/2}\psi_{t,h}(T_i, \mathbf{X}_i, Y_i, \upsilon_1, y_N)-N^{1/2}\psi_{t,h}(T_i, \mathbf{X}_i, Y_i, \upsilon, y_N)\right|^2\right] \leq C_3\delta^{\delta_3}
		\end{eqnarray}
		for some positive constants $C_3$ and $\delta_3$, and thus we conclude the stochastic continuity of $\{\mathcal{E}_{1, t,h}(\upsilon): \upsilon\in[1,\Delta]\}$.
		
		Finally,
		\begin{eqnarray}
			\notag & & N\mathbb{E}\left[\left\{\psi_{t,h}(T_i, \mathbf{X}_i, Y_i, \upsilon_1, y_N)\right\}\left\{\psi_{t,h}(T_i, \mathbf{X}_i, Y_i, \upsilon_2, y_N)\right\}\right] \\
			%\notag & = & \frac{h\left\{\bar{F}_{t}(\upsilon_1 y_N)\right\}\left\{\bar{F}_{t}(\upsilon_2 y_N)\right\}}{N}\mathbb{E}\left[\left\{\frac{\pi_0(T_i,\mathbf{X}_i)\mathbbm{1}(Y_i>\upsilon_1 y_N)\frac{1}{h}K\left(\frac{T_i-t}{h}\right)}{\mathbb{E}\left\{\Lambda_{N,1}(\upsilon_1 y_N,t,h)\right\}}-1\right\}\right. \\
			%\notag & & \left.~~~~~~~~~~~~~~~~~~~~~~~~~~~~~~~~~~~~~~~\times \left\{\frac{\pi_0(T_i,\mathbf{X}_i)\mathbbm{1}(Y_i>\upsilon_2 y_N)\frac{1}{h}K\left(\frac{T_i-t}{h}\right)}{\mathbb{E}\left\{\Lambda_{N,1}(\upsilon_2 y_N,t,h)\right\}}-1\right\}\right] \\
			\notag & = & h\sqrt{\bar{F}_{t}(\upsilon_1 y_N)\bar{F}_{t}(\upsilon_2 y_N)} \left(\frac{\mathbb{E}\left[\left\{\pi_{0}(T_i,\mathbf{X}_i)\right\}^2\mathbbm{1}(Y_i>\upsilon_1 y_N)\mathbbm{1}(Y_i>\upsilon_2 y_N)\frac{1}{h^2}\left\{K\left(\frac{T_i-t}{h}\right)\right\}^2\right]}{\bar{F}_{t}(\upsilon_1 y_N)f_{T}(t)\bar{F}_{t}(\upsilon_2 y_N)f_{T}(t)}\right. \\
			\notag & & ~~~~~~~~\left.-\frac{\mathbb{E}\left\{\Lambda_{N,1}(\upsilon_1 y_N, t, h)\right\}}{\bar{F}_{t}(\upsilon_1 y_N)f_{T}(t)}-\frac{\mathbb{E}\left\{\Lambda_{N,1}(\upsilon_2 y_N, t, h)\right\}}{\bar{F}_{t}(\upsilon_2 y_N)f_{T}(t)}  +1\right) \\
			\notag & = & h\sqrt{\bar{F}_{t}(\upsilon_1 y_N)\bar{F}_{t}(\upsilon_2 y_N)}\\
			\notag & &  \times \left(\frac{\kappa_{02}h^{-1}f_{T}(t)\mathbb{E}\left[\left\{\pi_{0}(T,\mathbf{X})\right\}^2\mathbbm{1}(Y_i>\upsilon_1 y_N)\mathbbm{1}(Y_i>\upsilon_2 y_N)|T=t\right]\{1+o(1)\}}{\bar{F}_{t}(\upsilon_1 y_N)\bar{F}_{t}(\upsilon_2 y_N)\{f_{T}(t)\}^2}-1\right) \\
			\notag & & + h\sqrt{\bar{F}_{t}(\upsilon_1 y_N)\bar{F}_{t}(\upsilon_2 y_N)}\left\{1+O\left(\omega_{t}(y_{N},h)\log(y_N)\right) + O(h)\right\}  \\
			\notag & = & \frac{\kappa_{02}}{f_{T}(t)} \frac{\mathbb{E}\left[\left\{\pi_{0}(T,\mathbf{X})\right\}^2\mathbbm{1}(Y>\upsilon_1 y_N)\mathbbm{1}(Y>\upsilon_2 y_N)|T=t\right]}{\sqrt{\bar{F}_{t}(\upsilon_1 y_N)\bar{F}_{t}(\upsilon_2 y_N)}}\{1+o(1)\} + o(1),
		\end{eqnarray}
		where the second equation is from Lemma \ref{lemma:03} and 
		\begin{eqnarray}
			\notag & & \mathbb{E}\left[\left\{\pi_{0}(T_i,\mathbf{X}_i)\right\}^2\mathbbm{1}(Y_i>\upsilon_1 y_N)\mathbbm{1}(Y_i>\upsilon_2 y_N)\frac{1}{h^2}\left\{K\left(\frac{T_i-t}{h}\right)\right\}^2\right] \\
			\notag & = & \int \left\{\pi_{0}(t^{\prime},\mathbf{x})\right\}^2\mathbbm{1}(y>\upsilon_1 y_N)\mathbbm{1}(y>\upsilon_2 y_N)\frac{1}{h^2}\left\{K\left(\frac{t^{\prime}-t}{h}\right)\right\}^2 f_{Y,T,\mathbf{X}}(y,t^{\prime},\mathbf{x})dydt^{\prime}d\mathbf{x} \\
			\notag & = & \frac{1}{h}\int \left\{\pi_{0}(t+hu,\mathbf{x})\right\}^2\mathbbm{1}(y>\upsilon_1 y_N)\mathbbm{1}(y>\upsilon_2 y_N)\left\{K\left(u\right)\right\}^2 f_{Y,T,\mathbf{X}}(y,t+hu,\mathbf{x})dydud\mathbf{x} \\
			\notag & = & \left[\frac{1}{h}\int \left\{\pi_{0}(t,\mathbf{x})\right\}^2\mathbbm{1}(y>\upsilon_1 y_N)\mathbbm{1}(y>\upsilon_2 y_N)\left\{K\left(u\right)\right\}^2 f_{Y,T,\mathbf{X}}(y,t,\mathbf{x})dydud\mathbf{x}\right]\{1+o(1)\} %\textcolor{red}{\text{~~more details needed here}} 
			\\
			\notag & = & \kappa_{02}h^{-1}f_{T}(t)\mathbb{E}\left[\left\{\pi_{0}(T,\mathbf{X})\right\}^2\mathbbm{1}(Y_i>\upsilon_1 y_N)\mathbbm{1}(Y_i>\upsilon_2 y_N)|T=t\right]\{1+o(1)\}.
		\end{eqnarray}
		We complete the first step.
		
		\textbf{Step 2}. The joint limiting behavior of $\mathcal{E}_{t_1,h}(\upsilon)$ and $\mathcal{E}_{t_2,h}(t)$ follows from the result in Step 1 and
		\begin{eqnarray}
			\notag & & N\mathbb{E}\left[\left\{\psi_{t_1,h}(T_i, \mathbf{X}_i, Y_i, \upsilon_1, y_N)\right\}\left\{\psi_{t_2,h}(T_i, \mathbf{X}_i, Y_i, \upsilon_2, y_N)\right\}\right] \\
			\notag & = & h\sqrt{\bar{F}_{t_1}(\upsilon_1 y_N)\bar{F}_{t_2}(\upsilon_2 y_N)} \left(\frac{\mathbb{E}\left[\left\{\pi_{0}(T_i,\mathbf{X}_i)\right\}^2\mathbbm{1}(Y_i>\upsilon_1 y_N)\mathbbm{1}(Y_i>\upsilon_2 y_N)\frac{1}{h^2}K\left(\frac{T_i-t_1}{h}\right)K\left(\frac{T_i-t_2}{h}\right)\right]}{\bar{F}_{t_1}(\upsilon_1 y_N)f_{T}(t_1)\bar{F}_{t_2}(\upsilon_2 y_N)f_{T}(t_2)}\right. \\
			\notag & & ~~~~~~~~\left.-\frac{\mathbb{E}\left\{\Lambda_{N,1}(\upsilon_1 y_N, t_1, h)\right\}}{\bar{F}_{t_1}(\upsilon_1 y_N)f_{T}(t_1)}-\frac{\mathbb{E}\left\{\Lambda_{N,1}(\upsilon_2 y_N, t_2, h)\right\}}{\bar{F}_{t_2}(\upsilon_2 y_N)f_{T}(t_2)}  +1\right) \\
			\notag & = & h\sqrt{\bar{F}_{t_1}(\upsilon_1 y_N)\bar{F}_{t_2}(\upsilon_2 y_N)}\\
			\notag & &  \times \left(\frac{\mathbb{E}\left[\left\{\pi_{0}(T_i,\mathbf{X}_i)\right\}^2\mathbbm{1}(Y_i>\upsilon_1 y_N)\mathbbm{1}(Y_i>\upsilon_2 y_N)\frac{1}{h^2}K\left(\frac{T_i-t_1}{h}\right)K\left(\frac{T_i-t_2}{h}\right)\right]}{\bar{F}_{t}(\upsilon_1 y_N)\bar{F}_{t}(\upsilon_2 y_N)f_{T}(t_1)f_{T}(t_2)}-1\right) \\
			\notag & & + h\sqrt{\bar{F}_{t_1}(\upsilon_1 y_N)\bar{F}_{t_2}(\upsilon_2 y_N)}\left\{1+O\left(\omega_{t_1}(y_{N},h)\log(y_N)\right) +O\left(\omega_{t_2}(y_{N},h)\log(y_N)\right) + O(h)\right\}  \\
			\notag & = & h\sqrt{\bar{F}_{t_1}(\upsilon_1 y_N)\bar{F}_{t_2}(\upsilon_2 y_N)}  \left(\frac{o\left(h^{-1}\bar{F}_{t}(y_N)\right)}{\bar{F}_{t}(\upsilon_1 y_N)\bar{F}_{t}(\upsilon_2 y_N)f_{T}(t_1)f_{T}(t_2)}-1\right) \\
			\notag & & + h\sqrt{\bar{F}_{t_1}(\upsilon_1 y_N)\bar{F}_{t_2}(\upsilon_2 y_N)}\left\{1+O\left(\omega_{t_1}(y_{N},h)\log(y_N)\right) +O\left(\omega_{t_2}(y_{N},h)\log(y_N)\right) + O(h)\right\}  \\
			\notag & = & o(1),
		\end{eqnarray}
		where the second last equality follows from
		\begin{eqnarray}
			\notag & & \mathbb{E}\left[\left\{\pi_{0}(T_i,\mathbf{X}_i)\right\}^2\mathbbm{1}(Y_i>\upsilon_1 y_N)\mathbbm{1}(Y_i>\upsilon_2 y_N)\frac{1}{h^2}K\left(\frac{T_i-t_1}{h}\right)K\left(\frac{T_i-t_2}{h}\right)\right] \\
			\notag & = & \mathbb{E}\left[\left\{\pi_{0}(T_i,\mathbf{X}_i)\right\}^2\mathbbm{1}(Y_i>(\upsilon_1 \vee \upsilon_2) y_N)\frac{1}{h^2}K\left(\frac{T_i-t_1}{h}\right)K\left(\frac{T_i-t_2}{h}\right)\right] \\
			\notag & = & \int\frac{1}{h^2}K\left(\frac{t^{\prime}-t_1}{h}\right)K\left(\frac{t^{\prime}-t_2}{h}\right)\{\pi_{0}(t^{\prime},\mathbf{x})\}^2\mathbb{P}\left(Y(t^{\prime})>(\upsilon_1 \vee \upsilon_2) y_N|T=t^{\prime},\mathbf{X}=\mathbf{x}\right)\\
			\notag & & ~~~~~~~~~~~~~~~~~~~~~~~~~~~~~~~~~~~~~~~~~~~~~\times f_{T|\mathbf{X}}(t^{\prime}|\mathbf{x})f_{\mathbf{X}}(\mathbf{x})d\mathbf{x}dt^{\prime} \\
			\notag & \leq & \eta_2\int\frac{1}{h^2}K\left(\frac{t^{\prime}-t_1}{h}\right)K\left(\frac{t^{\prime}-t_2}{h}\right)\pi_{0}(t^{\prime},\mathbf{x})\mathbb{P}\left(Y(t^{\prime})>(\upsilon_1 \vee \upsilon_2) y_N|T=t^{\prime},\mathbf{X}=\mathbf{x}\right)\\
			\notag & & ~~~~~~~~~~~~~~~~~~~~~~~~~~~~~~~~~~~~~~~~~~~~~\times f_{T|\mathbf{X}}(t^{\prime}|\mathbf{x})f_{\mathbf{X}}(\mathbf{x})d\mathbf{x}dt^{\prime} \\
			\notag & = & \eta_2\int\frac{1}{h^2}K\left(\frac{t^{\prime}-t_1}{h}\right)K\left(\frac{t^{\prime}-t_2}{h}\right)\bar{F}_{t^{\prime}}((\upsilon_{1} \vee \upsilon_{2})y_N)f_{T}(t^{\prime})dt^{\prime}  \\
			\notag & = & \frac{\eta_2}{h}\int_{\mathcal{K}} K\left(u\right)K\left(\frac{t_1-t_2}{h}+u\right)\bar{F}_{t+hu}((\upsilon_{1} \vee \upsilon_{2})y_N)f_{T}(t+hu)du \\
			\notag & = & \frac{\eta_2}{h}\int_{\mathcal{K}} K\left(u\right)K\left(\frac{t_1-t_2}{h}+u\right)\bar{F}_{t}((\upsilon_{1} \vee \upsilon_{2})y_N)f_{T}(t+hu)du \\
			\notag & & +
			\frac{\eta_2}{h}\int_{\mathcal{K}} K\left(u\right)K\left(\frac{t_1-t_2}{h}+u\right)\left\{\bar{F}_{t+hu}((\upsilon_{1} \vee \upsilon_{2})y_N)-\bar{F}_{t}((\upsilon_{1} \vee \upsilon_{2})y_N)\right\}f_{T}(t+hu)du \\
			\notag & = & \frac{\eta_2\bar{F}_{t}((\upsilon_{1} \vee \upsilon_{2})y_N)}{h}\int_{\mathcal{K}} K\left(u\right)K\left(\frac{t_1-t_2}{h}+u\right)\{f_{T}(t)+O(h)\}du \\
			\notag & & +
			\frac{\eta_2\bar{F}_{t}((\upsilon_{1} \vee \upsilon_{2})y_N)}{h}\times O\left(\omega_{t}(y_N, h)\log(y_N)\right)\int_{\mathcal{K}} K\left(u\right)K\left(\frac{t_1-t_2}{h}+u\right)\{f_{T}(t)+O(h)\}du \\
			\notag & = & o\left(h^{-1}\bar{F}_{t}(y_N)\right),
		\end{eqnarray}
		where $\int_{\mathcal{K}} K\left(u\right)K\left(\frac{t_1-t_2}{h}+u\right)du \to 0 $
		when $h\to0$ as the support set of $K(\cdot)$ is included in $[-1,1]$. Thus, we finish the proof of this lemma.
	\end{proof}

	\subsection{Proof of Theorem \ref{theo:Omega_est}}
	
	\begin{proof}
		By simple algebra,
		\begin{eqnarray}
			\notag & & \hat{\psi}_{t,h}(T_i, \mathbf{X}_i, Y_i, \upsilon, y_N) \\
			\notag & = & \sqrt{\frac{h\widehat{\bar{F}}_{t}(\upsilon y_N)}{N}}\left\{\frac{\hat{\pi}(T_i,\mathbf{X}_i)\mathbbm{1}(Y_i>\upsilon y_N)\frac{1}{h}K\left(\frac{T_i-t}{h}\right)}{\widehat{\bar{F}}_{t}(\upsilon y_N)\hat{f}_{T,h}(t)}-1\right\} \\
			\notag & = &\sqrt{\frac{\widehat{\bar{F}}_{t}(\upsilon y_N)}{\bar{F}_{t}(\upsilon y_N)}}\frac{\bar{F}_{t}(\upsilon y_N)f_{T}(t)}{\widehat{\bar{F}}_{t}(\upsilon y_N)\hat{f}_{T,h}(t)}\tilde{\psi}_{t,h}(T_i, \mathbf{X}_i, Y_i, \upsilon, y_N) + \sqrt{\frac{h\widehat{\bar{F}}_{t}(\upsilon y_N)}{N}}\left\{\frac{\bar{F}_{t}(\upsilon y_N)f_{T}(t)}{\widehat{\bar{F}}_{t}(\upsilon y_N)\hat{f}_{T,h}(t)}-1\right\},
		\end{eqnarray}
		where
		\begin{eqnarray}
			\notag \tilde{\psi}_{t,h}(T_i, \mathbf{X}_i, Y_i, \upsilon, y_N)=\sqrt{\frac{h\bar{F}_{t}(\upsilon y_N)}{N}}\left\{\frac{\hat{\pi}(T_i,\mathbf{X}_i)\mathbbm{1}(Y_i>\upsilon y_N)\frac{1}{h}K\left(\frac{T_i-t}{h}\right)}{\bar{F}_{t}(\upsilon y_N)f_{T}(t)}-1\right\}.
		\end{eqnarray}
		
		Then, we have
		\begin{eqnarray}
			\notag & & \hat{\Omega}_{t}^{F}(\upsilon_1,\upsilon_2) \\
			\notag & = & \sum_{i=1}^{N} \hat{\psi}_{t,h}(T_i, \mathbf{X}_i, Y_i, \upsilon_1, y_N) \hat{\psi}_{t,h}(T_i, \mathbf{X}_i, Y_i, \upsilon_2, y_N) \\
			\notag & = & \sqrt{\frac{\widehat{\bar{F}}_{t}(\upsilon_1 y_N)}{\bar{F}_{t}(\upsilon_1 y_N)}}\sqrt{\frac{\widehat{\bar{F}}_{t}(\upsilon_2 y_N)}{\bar{F}_{t}(\upsilon_2 y_N)}}\frac{\bar{F}_{t}(\upsilon_1 y_N)f_{T}(t)}{\widehat{\bar{F}}_{t}(\upsilon_1 y_N)\hat{f}_{T,h}(t)}\frac{\bar{F}_{t}(\upsilon_2 y_N)f_{T}(t)}{\widehat{\bar{F}}_{t}(\upsilon_2 y_N)\hat{f}_{T,h}(t)} \\
			\notag & & \times\sum_{i=1}^{N}\tilde{\psi}_{t,h}(T_i, \mathbf{X}_i, Y_i, \upsilon_1, y_N)\tilde{\psi}_{t,h}(T_i, \mathbf{X}_i, Y_i, \upsilon_2, y_N) \\
			\notag & & + \sqrt{\frac{\widehat{\bar{F}}_{t}(\upsilon_1 y_N)}{\bar{F}_{t}(\upsilon_1 y_N)}}\frac{\bar{F}_{t}(\upsilon_1 y_N)f_{T}(t)}{\widehat{\bar{F}}_{t}(\upsilon_1 y_N)\hat{f}_{T,h}(t)} \sqrt{\frac{h\widehat{\bar{F}}_{t}(\upsilon_2 y_N)}{N}}\left\{\frac{\bar{F}_{t}(\upsilon_2 y_N)f_{T}(t)}{\widehat{\bar{F}}_{t}(\upsilon_2 y_N)\hat{f}_{T,h}(t)}-1\right\}\sum_{i=1}^{N}\tilde{\psi}_{t,h}(T_i, \mathbf{X}_i, Y_i, \upsilon_1, y_N) \\
			\notag & & + \sqrt{\frac{\widehat{\bar{F}}_{t}(\upsilon_2 y_N)}{\bar{F}_{t}(\upsilon_2 y_N)}}\frac{\bar{F}_{t}(\upsilon_2 y_N)f_{T}(t)}{\widehat{\bar{F}}_{t}(\upsilon_2 y_N)\hat{f}_{T,h}(t)} \sqrt{\frac{h\widehat{\bar{F}}_{t}(\upsilon_1 y_N)}{N}}\left\{\frac{\bar{F}_{t}(\upsilon_1 y_N)f_{T}(t)}{\widehat{\bar{F}}_{t}(\upsilon_1 y_N)\hat{f}_{T,h}(t)}-1\right\}\sum_{i=1}^{N}\tilde{\psi}_{t,h}(T_i, \mathbf{X}_i, Y_i, \upsilon_2, y_N) \\
			\notag & & + \sqrt{\frac{h\widehat{\bar{F}}_{t}(\upsilon_1 y_N)}{N}}\left\{\frac{\bar{F}_{t}(\upsilon_1 y_N)f_{T}(t)}{\widehat{\bar{F}}_{t}(\upsilon_1 y_N)\hat{f}_{T,h}(t)}-1\right\}\sqrt{\frac{h\widehat{\bar{F}}_{t}(\upsilon_2 y_N)}{N}}\left\{\frac{\bar{F}_{t}(\upsilon_2 y_N)f_{T}(t)}{\widehat{\bar{F}}_{t}(\upsilon_2 y_N)\hat{f}_{T,h}(t)}-1\right\}.
		\end{eqnarray}
		
		It is straightforward to show that $\frac{\hat{f}_{T,h}(t)}{f_{T}(t)}-1=O_{p}(N^{-1/2}h^{-1/2})$ under the conditions assumed in this theorem. In addition, Lemma \ref{lemma:01} implies that $\frac{\widehat{\bar{F}}_{t,h}(\upsilon y_N)}{\bar{F}_{t}(\upsilon y_N)}-1 = O_{p}\left(N^{-1/2}h^{-1/2}\{\bar{F}_{t}(y_N)\}^{-1/2}\right)$ uniformly in $\upsilon\in[1,\Delta]$. Thus, we arrive at
		\begin{eqnarray}
			\notag & & \hat{\Omega}_{t}^{F}(\upsilon_1,\upsilon_2) \\
			\notag & = & \left\{1+O_{p}\left(N^{-1/2}h^{-1/2}\{\bar{F}_{t}(y_N)\}^{-1/2}\right)\right\}\sum_{i=1}^{N}\tilde{\psi}_{t,h}(T_i, \mathbf{X}_i, Y_i, \upsilon_1, y_N)\tilde{\psi}_{t,h}(T_i, \mathbf{X}_i, Y_i, \upsilon_2, y_N) \\
			\notag & & + \left\{1+O_{p}\left(N^{-1/2}h^{-1/2}\{\bar{F}_{t}(y_N)\}^{-1/2}\right)\right\} \sum_{i=1}^{N}\tilde{\psi}_{t,h}(T_i, \mathbf{X}_i, Y_i, \upsilon_1, y_N)\times O_{p}\left(N^{-1}\right) \\
			\notag & & + \left\{1+O_{p}\left(N^{-1/2}h^{-1/2}\{\bar{F}_{t}(y_N)\}^{-1/2}\right)\right\} \sum_{i=1}^{N}\tilde{\psi} _{t,h}(T_i, \mathbf{X}_i, Y_i, \upsilon_2, y_N) \times O_{p}\left(N^{-1}\right) \\
			\notag & & + O_{p}\left(N^{-2}\right).
		\end{eqnarray}
		Regarding that
		\begin{eqnarray}
		    \notag \psi_{t,h}(T_i, \mathbf{X}_i, Y_i, \upsilon, y_N)=\sqrt{\frac{h\bar{F}_{t}(\upsilon y_N)}{N}}\left\{\frac{\pi_0(T_i,\mathbf{X}_i)\mathbbm{1}(Y_i>\upsilon y_N)\frac{1}{h}K\left(\frac{T_i-t}{h}\right)}{\bar{F}_{t}(\upsilon y_N)f_{T}(t)}-1\right\}
		\end{eqnarray}
		and
		\begin{eqnarray}
		    \notag & & \tilde{\psi} _{t,h}(T_i, \mathbf{X}_i, Y_i, \upsilon, y_N) \\
		    \notag & = & \psi _{t,h}(T_i, \mathbf{X}_i, Y_i, \upsilon, y_N) + \sqrt{\frac{h\bar{F}_{t}(\upsilon y_N)}{N}}\frac{\left\{\hat{\pi}(T_i,\mathbf{X}_i)-\pi_0(T_i,\mathbf{X}_i)\right\}\mathbbm{1}(Y_i>\upsilon y_N)\frac{1}{h}K\left(\frac{T_i-t}{h}\right)}{\bar{F}_{t}(\upsilon y_N)f_{T}(t)},
		\end{eqnarray}
		where 
		\begin{eqnarray}
		    \notag & & \left|\sqrt{\frac{h\bar{F}_{t}(\upsilon y_N)}{N}}\frac{\left\{\hat{\pi}(T_i,\mathbf{X}_i)-\pi_0(T_i,\mathbf{X}_i)\right\}\mathbbm{1}(Y_i>\upsilon y_N)\frac{1}{h}K\left(\frac{T_i-t}{h}\right)}{\bar{F}_{t}(\upsilon y_N)f_{T}(t)}\right| \\
		    \notag & = & \sqrt{\frac{h\bar{F}_{t}(\upsilon y_N)}{N}}\frac{\mathbbm{1}(Y_i>\upsilon y_N)\frac{1}{h}K\left(\frac{T_i-t}{h}\right)}{\bar{F}_{t}(\upsilon y_N)f_{T}(t)}\times O_{p}(\delta_N) \\
		    \notag & \leq & \frac{1}{\eta_1}\left\{\sqrt{\frac{h\bar{F}_{t}(\upsilon y_N)}{N}}+\psi _{t,h}(T_i, \mathbf{X}_i, Y_i, \upsilon_2, y_N)\right\}\times O_{p}(\delta_N)
		\end{eqnarray}
		for all $1\leq i\leq N$. We have that
		\begin{eqnarray}
		    \notag & & \hat{\Omega}_{t}^{F}(\upsilon_1,\upsilon_2) \\
			\notag & = & \{1+O_{p}(\delta_N)\}\left\{1+O_{p}\left(N^{-1/2}h^{-1/2}\{\bar{F}_{t}(\upsilon y_N)\}^{-1/2}\right)\right\} \\
			\notag & & \times\sum_{i=1}^{N}\psi_{t,h}(T_i, \mathbf{X}_i, Y_i, \upsilon_1, y_N)\psi_{t,h}(T_i, \mathbf{X}_i, Y_i, \upsilon_2, y_N) \\
			\notag & & + \{1+O_{p}(\delta_N)\}\left\{1+O_{p}\left(N^{-1/2}h^{-1/2}\{\bar{F}_{t}(\upsilon y_N)\}^{-1/2}\right)\right\} O_{p}\left(N^{-1}\right)\sum_{i=1}^{N}\psi_{t,h}(T_i, \mathbf{X}_i, Y_i, \upsilon_1, y_N) \\
			\notag & & + \{1+O_{p}(\delta_N)\}\left\{1+O_{p}\left(N^{-1/2}h^{-1/2}\{\bar{F}_{t}(\upsilon y_N)\}^{-1/2}\right)\right\} O_{p}\left(N^{-1}\right)\sum_{i=1}^{N}\psi_{t,h}(T_i, \mathbf{X}_i, Y_i, \upsilon_2, y_N) \\
			\notag & & + O_{p}\left(N^{-2}\right) + O_{p}(h\bar{F}_{t}(y_N)\delta_N^2)+O_{p}(N^{-1/2}h^{1/2}\{\bar{F}_{t}(y_N)\}^{1/2}\delta_N),
		\end{eqnarray}
		where $\sum_{i=1}^{N}\psi_{t,h}(T_i, \mathbf{X}_i, Y_i, \upsilon_1, y_N)=O_{p}\left(1\right)$ and $\sum_{i=1}^{N}\psi_{t,h}(T_i, \mathbf{X}_i, Y_i, \upsilon_2, y_N)=O_{p}\left(1\right)$. Therefore,
		\begin{eqnarray}
			\notag \hat{\Omega}_{t}^{F}(\upsilon_1,\upsilon_2) = \sum_{i=1}^{N}\psi_{t,h}(T_i, \mathbf{X}_i, Y_i, \upsilon_1, y_N)\psi_{t,h}(T_i, \mathbf{X}_i, Y_i, \upsilon_2, y_N)\{1+o_{p}(1)\} + o_{p}(1).
		\end{eqnarray}
		Then, the proof of this theorem follows from $$\sum_{i=1}^{N}\psi_{t,h}(T_i, \mathbf{X}_i, Y_i, \upsilon_1, y_N)\psi_{t,h}(T_i, \mathbf{X}_i, Y_i, \upsilon_2, y_N)/\Omega_{t}^{F}(\upsilon_{1},\upsilon_{2})\to1$$ 
		in probability.
	\end{proof}

	\subsection{Proof of Theorem \ref{theorem:quantile_intermediate}}
	
	\begin{proof}
		Regarding the fact that $\mathcal{E}_{t_1,h}$ and $\mathcal{E}_{t_2,h}$ are asymptotically independent implied from Lemma \ref{lemma:01}, it suffices to show that for any fixed $t\in\mathcal{T}$,
		$$\mathcal{Q}_{t,h}(\upsilon)\Rightarrow \Psi_{t}^{Q}(\upsilon)$$
		in $\ell^{\infty}([\Delta_0,1])$ as $N\to\infty$, where $\Psi_{t}^{Q}(\upsilon)$ is a centered Gaussian process with covariance function $\Omega_{t}^{Q}(\upsilon_1,\upsilon_2)$.
		
		Rewrite
		\begin{eqnarray}
			\notag \notag \mathcal{E}_{t,h}(\upsilon) & = & \sqrt{Nh\bar{F}_{t}(\upsilon y_N)}\left\{\frac{\widehat{\bar{F}}_{t,h}(\upsilon y_N)}{\bar{F}_{t}(\upsilon y_N)}-1\right\}\\
			\notag & = & \sqrt{\frac{Nh}{\bar{F}_{t}(\upsilon y_N)}}\bar{F}_{t}(y_N)\left\{\frac{\widehat{\bar{F}}_{t,h}(\upsilon y_N)}{\bar{F}_{t}(y_N)}-\frac{\bar{F}_{t}(\upsilon y_N)}{\bar{F}_{t}(y_N)}\right\} \\
			\notag & = & \sqrt{Nh\bar{F}_{t}(y_N)}\left\{\frac{\bar{F}_{t}(\upsilon y_N)}{\bar{F}_{t}(y_N)}\right\}^{-1/2}\left\{\frac{\widehat{\bar{F}}_{t,h}(\upsilon y_N)}{\bar{F}_{t}(y_N)}-\upsilon^{-1/\gamma(t)}\right\} \\
			\notag & & + \sqrt{Nh\bar{F}_{t}(y_N)}\left\{\frac{\bar{F}_{t}(\upsilon y_N)}{\bar{F}_{t}(y_N)}\right\}^{-1/2}\left\{\upsilon^{-1/\gamma(t)}-\frac{\bar{F}_{t}(\upsilon y_N)}{\bar{F}_{t}(y_N)}\right\} \\
			\notag & = & \sqrt{Nh\bar{F}_{t}(y_N)}\upsilon^{1/\{2\gamma(t)\}}\left\{\frac{\widehat{\bar{F}}_{t,h}(\upsilon y_N)}{\bar{F}_{t}(y_N)}-\upsilon^{-1/\gamma(t)}\right\} \\
			\notag & + & \sqrt{Nh\bar{F}_{t}(y_N)}\left[\left\{\frac{\bar{F}_{t}(\upsilon y_N)}{\bar{F}_{t}(y_N)}\right\}^{-1/2}-\upsilon^{1/(2\gamma(t))}\right]\left\{\frac{\widehat{\bar{F}}_{t,h}(\upsilon y_N)}{\bar{F}_{t}(y_N)}-\upsilon^{-1/\gamma(t)}\right\} \\
			\notag & & + \sqrt{Nh\bar{F}_{t}(y_N)}\left\{\frac{\bar{F}_{t}(\upsilon y_N)}{\bar{F}_{t}(y_N)}\right\}^{-1/2}\left\{\upsilon^{-1/\gamma(t)}-\frac{\bar{F}_{t}(\upsilon y_N)}{\bar{F}_{t}(y_N)}\right\}.
		\end{eqnarray}
		According to Lemma \ref{lemma:01} and Lemma \ref{lemma:F_v},
		\begin{eqnarray}\label{eq:proof_4_6_01}
			\notag & & \sqrt{Nh\bar{F}_{t}(y_N)}\left[\left\{\frac{\bar{F}_{t}(\upsilon y_N)}{\bar{F}_{t}(y_N)}\right\}^{-1/2}-\upsilon^{1/(2\gamma(t))}\right]\left\{\frac{\widehat{\bar{F}}_{t,h}(\upsilon y_N)}{\bar{F}_{t}(y_N)}-\upsilon^{-1/\gamma(t)}\right\} \\
			\notag & = & 
			\sqrt{Nh\bar{F}_{t}(y_N)}\left[\left\{\frac{\bar{F}_{t}(\upsilon y_N)}{\bar{F}_{t}(y_N)}\right\}^{-1/2}-\upsilon^{1/(2\gamma(t))}\right]\left\{\frac{\widehat{\bar{F}}_{t,h}(\upsilon y_N)}{\bar{F}_{t}(\upsilon y_N)}-1\right\}\left\{\frac{\bar{F}_{t}(\upsilon y_N)}{\bar{F}_{t}(y_N)}\right\} \\
			\notag & & + \sqrt{Nh\bar{F}_{t}(y_N)}\left[\left\{\frac{\bar{F}_{t}(\upsilon y_N)}{\bar{F}_{t}(y_N)}\right\}^{-1/2}-\upsilon^{1/(2\gamma(t))}\right]\left\{\frac{\bar{F}_{t}(\upsilon y_N)}{\bar{F}_{t}(y_N)}-\upsilon^{-1/\gamma(t)}\right\} \\
			& = & o_{p}(1)
		\end{eqnarray}
		and
		\begin{eqnarray}\label{eq:proof_4_6_02}
			\sqrt{Nh\bar{F}_{t}(y_N)}\left\{\frac{\bar{F}_{t}(\upsilon y_N)}{\bar{F}_{t}(y_N)}\right\}^{-1/2}\left\{\upsilon^{-1/\gamma(t)}-\frac{\bar{F}_{t}(\upsilon y_N)}{\bar{F}_{t}(y_N)}\right\} = o(1)
		\end{eqnarray}
		uniformly in $\upsilon\in[1,\Delta]$ \textcolor{black}{under the condition that $\sqrt{Nh(1-\alpha_N)}A_{t}(1/(1-\alpha_N))=o(1)$ with $y_N=q_{t}(\alpha_N)$}. Thus, 
		\begin{eqnarray}
			\mathcal{E}_{0, t,h}(\upsilon)\Rightarrow \Psi_{t}^{Q1}(\upsilon)
		\end{eqnarray}
		in $\ell^{\infty}([1,\Delta])$, where
		\begin{eqnarray}
			\notag \mathcal{E}_{0, t,h}(\upsilon) = \upsilon^{1/\{2\gamma(t)\}}\sqrt{Nh(1-\alpha_N)}\left\{\frac{\widehat{\bar{F}}_{t,h}(\upsilon q_{t}(\alpha_N))}{1-\alpha_N}-\upsilon^{-1/\gamma(t)}\right\}.
		\end{eqnarray}
		and $\Psi_{t}^{Q1}(\upsilon)$ is a centered Gaussian process with covariance function
		\begin{eqnarray}
			\notag \Omega_{t}^{Q1}(\upsilon_1,\upsilon_2) =
			\kappa_{02} \{f_{T}(t)\}^{-1}\varpi_{t}^{Q1}(\upsilon_1, \upsilon_2),
		\end{eqnarray}
		where $\varpi_{t}^{Q1}(\upsilon_1,\upsilon_2)=\lim_{\alpha_N\to 1}\frac{ \mathbb{E}\left[\left\{\pi_{0}(T,\mathbf{X})\right\}^2\mathbbm{1}(Y>\upsilon_1 q_{t}(\alpha_N))\mathbbm{1}(Y>\upsilon_2 q_{t}(\alpha_N))|T=t\right]}{\sqrt{\bar{F}_{t}(\upsilon_1q_{t}(\alpha_N))\bar{F}_{t}(\upsilon_2q_{t}(\alpha_N))}}$.
		\iffalse
		\begin{eqnarray}
			\notag \Omega_{t}(\upsilon_1,\upsilon_2) & = & \lim_{y_N\to\infty}
			\kappa_{02} \{f_{T}(t)\}^{-1}\{\bar{F}_{t}(\upsilon_1 y_N)\}^{-1/2}\{\bar{F}_{t}(\upsilon_2 y_N)\}^{-1/2} \\
			\notag & & ~~~~~~~~\times \mathbb{E}\left[\left\{\pi_{0}(T,\mathbf{X})\right\}^2\mathbbm{1}(Y>\upsilon_1 y_N)\mathbbm{1}(Y>\upsilon_2 y_N)|T=t\right] \\
			\notag & = & \lim_{\alpha_N\to 1}
			\kappa_{02} \{f_{T}(t)\}^{-1}\{\bar{F}_{t}(\upsilon_1 q_{t}(\alpha_N))\}^{-1/2}\{\bar{F}_{t}(\upsilon_2 q_{t}(\alpha_N))\}^{-1/2} \\
			\notag & & ~~~~~~~~\times \mathbb{E}\left[\left\{\pi_{0}(T,\mathbf{X})\right\}^2\mathbbm{1}(Y>\upsilon_1 q_{t}(\alpha_N))\mathbbm{1}(Y>\upsilon_2 q_{t}(\alpha_N))|T=t\right] 
		\end{eqnarray}
		\fi
		
		By Skorohod's representation theorem, on a suitable probability space, we have that
		\begin{eqnarray}
			\sup_{\upsilon\in[1,\Delta]}\upsilon^{1/\{2\gamma(t)\}}\left|\sqrt{Nh(1-\alpha_N)}\left\{\frac{\widehat{\bar{F}}_{t,h}(\upsilon q_{t}(\alpha_N))}{1-\alpha_N}-\upsilon^{-1/\gamma(t)}\right\}-\upsilon^{-1/\{2\gamma(t)\}}\Psi_{t}^{Q1}(\upsilon)\right|\to 0
		\end{eqnarray}
		almost surely, and this immediately implies that
		\begin{eqnarray}
			\sup_{\upsilon\in[1,\Delta]}\left|\sqrt{Nh(1-\alpha_N)}\left\{\frac{\widehat{\bar{F}}_{t,h}(\upsilon q_{t}(\alpha_N))}{1-\alpha_N}-\upsilon^{-1/\gamma(t)}\right\}-\upsilon^{-1/\{2\gamma(t)\}}\Psi_{t}^{Q1}(\upsilon)\right|\to 0
		\end{eqnarray}
		with probability one.
		By Vervaat's Lemma \citep[Lemma A.0.2]{deHaan2006}, we have
		\begin{eqnarray}\label{eq:proof_4_6_03}
			\sup_{\upsilon_0\in[ \Delta^{-1/\gamma(t)},1]}\left|\sqrt{Nh(1-\alpha_N)}\left\{ \frac{\hat{q}_{t, h}(1-\upsilon_0(1-\alpha_N))}{q_{t}(\alpha_N)} - \upsilon_0^{-\gamma(t)}\right\} - \gamma(t)\upsilon_0^{-\gamma(t)-1/2}\Psi_{t}^{Q1}(\upsilon_0^{-\gamma(t)})\right|\to 0
		\end{eqnarray}
		almost surely. It is worth mentioning that $\Psi_{t}^{Q1}(\upsilon)$ and $\Psi_{t}^{F}(\upsilon)$ in the proof of Lemma \ref{lemma:01} can be defined on the same probability space.
		
		\iffalse
		Assumption \ref{assump:second_order} and Theorem 2.3.9 of \citet{deHaan2006} indicate that $q_{t}(1-1/\cdot)\in 2RV_{\gamma(t),\rho(t)}$ with auxiliary function $A_{t}(\cdot)$. Using Lemma \ref{lemma:pan_spl}, for any fixed positive $\Delta_0<1$, we know uniformly in $\upsilon\in[\Delta_0, 1]$,
		\begin{eqnarray}
			\notag \left|\frac{q_{t}(\alpha_N)}{q_{t}(1-\upsilon(1-\alpha_N))}-\upsilon^{\gamma(t)}\right| = O\left(A_{t}(1/(1-\alpha_N))\right).
		\end{eqnarray}
		\fi
		
		Rewrite
		\begin{eqnarray}
			\notag & & \mathcal{Q}_{t,h}(\upsilon) \\
			\notag & = &  \sqrt{Nh\upsilon(1-\alpha_N)}\left\{\frac{\hat{q}_{t, h}(1-\upsilon(1-\alpha_N))}{q_{t}(1-\upsilon(1-\alpha_N))}-1\right\} \\
			\notag & = & \sqrt{Nh\upsilon(1-\alpha_N)}\left\{\frac{\hat{q}_{t, h}(1-\upsilon(1-\alpha_N))}{q_{t}(\alpha_N)}\frac{q_{t}(\alpha_N)}{q_{t}(1-\upsilon(1-\alpha_N))}-1\right\} \\
			\notag & = & \sqrt{Nh\upsilon(1-\alpha_N)}\left\{\frac{\hat{q}_{t, h}(1-\upsilon(1-\alpha_N))}{q_{t}(\alpha_N)}-\upsilon^{-\gamma(t)}\right\}\frac{q_{t}(\alpha_N)}{q_{t}(1-\upsilon(1-\alpha_N))} \\
			\notag & & + \sqrt{Nh\upsilon(1-\alpha_N)}\left\{\upsilon^{-\gamma(t)}\frac{q_{t}(\alpha_N)}{q_{t}(1-\upsilon(1-\alpha_N))}-1\right\} \\
			\notag & = & \upsilon^{\gamma(t)}\sqrt{Nh\upsilon(1-\alpha_N)}\left\{\frac{\hat{q}_{t, h}(1-\upsilon(1-\alpha_N))}{q_{t}(\alpha_N)}-\upsilon^{-\gamma(t)}\right\} \\
			\notag & & + \sqrt{Nh\upsilon(1-\alpha_N)}\left\{\frac{\hat{q}_{t, h}(1-\upsilon(1-\alpha_N))}{q_{t}(\alpha_N)}-\upsilon^{-\gamma(t)}\right\}\left\{\frac{q_{t}(\alpha_N)}{q_{t}(1-\upsilon(1-\alpha_N))}-\upsilon^{\gamma(t)}\right\} \\
			\notag & & + \sqrt{Nh\upsilon(1-\alpha_N)}\upsilon^{-\gamma(t)}\left\{\frac{q_{t}(\alpha_N)}{q_{t}(1-\upsilon(1-\alpha_N))}-\upsilon^{\gamma(t)}\right\}.
		\end{eqnarray}
		According to Lemma \ref{lemma:F_v} and \eqref{eq:proof_4_6_03}, it can be easily shown that
		\begin{eqnarray}
			\notag \sqrt{Nh\upsilon(1-\alpha_N)}\left\{\frac{\hat{q}_{t, h}(1-\upsilon(1-\alpha_N))}{q_{t}(\alpha_N)}-\upsilon^{-\gamma(t)}\right\}\left\{\frac{q_{t}(\alpha_N)}{q_{t}(1-\upsilon(1-\alpha_N))}-\upsilon^{\gamma(t)}\right\}  = o_{p}(1)
		\end{eqnarray}
		and
		\begin{eqnarray}\notag 
			\sqrt{Nh\upsilon(1-\alpha_N)}\upsilon^{-\gamma(t)}\left\{\frac{q_{t}(\alpha_N)}{q_{t}(1-\upsilon(1-\alpha_N))}-\upsilon^{\gamma(t)}\right\} = o(1)
		\end{eqnarray}
		uniformly in $\upsilon\in[\Delta_0, 1]$  \textcolor{black}{under the condition that $\sqrt{Nh(1-\alpha_N)}A_{t}(1/(1-\alpha_N))=o(1)$}.
		According to \eqref{eq:proof_4_6_03}, we obtain that
		\begin{eqnarray}
			\upsilon^{\gamma(t)}\sqrt{Nh\upsilon(1-\alpha_N)}\left\{\frac{\hat{q}_{t, h}(1-\upsilon(1-\alpha_N))}{q_{t}(\alpha_N)}-\upsilon^{-\gamma(t)}\right\} \Rightarrow \gamma(t)\Psi_{t}^{Q1}(\upsilon^{-\gamma(t)})
		\end{eqnarray}
		in $\ell^{\infty}([\Delta^{-1/\gamma(t)},1])$, which implies
		\begin{eqnarray}
		    \mathcal{Q}_{t,h}(\upsilon) \Rightarrow \gamma(t)\Psi_{t}^{Q1}(\upsilon^{-\gamma(t)})
		\end{eqnarray}
		in $\ell^{\infty}([\Delta^{-1/\gamma(t)},1])$. Let $\Delta_0=\Delta^{-1/\gamma(t)}$,
		according to the proof of Lemma \ref{lemma:01} and previous arguments, 
		\begin{eqnarray}
			\notag \sup_{\upsilon\in[\Delta_0,1]}\left|\mathcal{Q}_{t,h}(\upsilon)-\sum_{i=1}^{N}\gamma(t)\psi_{t,h}(T_i, \mathbf{X}_i, Y_i, \upsilon^{-\gamma(t)}, q_{t}(\alpha_N))\right| = o_{p}(1),
		\end{eqnarray}
		where
		\begin{eqnarray}
		    \notag \psi_{t,h}(T_i, \mathbf{X}_i, Y_i, \upsilon, q_{t}(\alpha_N))=\sqrt{\frac{h\bar{F}_{t}(\upsilon q_{t}(\alpha_N))}{N}}\left\{\frac{\pi_0(T_i,\mathbf{X}_i)\mathbbm{1}(Y_i>\upsilon q_{t}(\alpha_N))\frac{1}{h}K\left(\frac{T_i-t}{h}\right)}{\bar{F}_{t}(\upsilon q_{t}(\alpha_N))f_{T}(t)}-1\right\}
		\end{eqnarray}
		for $i=1,\ldots,N$. In addition,
		\begin{eqnarray}
		    \notag & & \psi_{t,h}(T_i, \mathbf{X}_i, Y_i, \upsilon^{-\gamma(t)}, q_{t}(\alpha_N)) \\
		    \notag & = & \sqrt{\frac{\upsilon(1-\alpha_N)}{\bar{F}_{t}(\upsilon^{-\gamma(t)}q_{t}(\alpha_N))}}\psi_{t,h}^{Q}(T_i, \mathbf{X}_i, Y_i, \upsilon, \alpha_N)  \\
		    \notag & & + \left\{\frac{\upsilon(1-\alpha_N)}{\bar{F}_{t}(\upsilon^{-\gamma(t)}q_{t}(\alpha_N))}-1\right\}\sqrt{\frac{h\bar{F}_{t}(\upsilon^{-\gamma(t)}q_{t}(\alpha_N))}{N}} \\
		    \notag & & + \sqrt{\frac{h}{N\bar{F}_{t}(\upsilon^{-\gamma(t)}q_{t}(\alpha_N))}}\frac{1}{f_{T}(t)} \\
		    \notag & & ~~\times \pi_0(T_i,\mathbf{X}_i)\frac{1}{h}K\left(\frac{T_i-t}{h}\right)\left\{\mathbbm{1}(Y_i>\upsilon^{-\gamma(t)}q_{t}(\alpha_N))-\mathbbm{1}(Y_i> q_{t}(1-\upsilon(1-\alpha_N)))\right\}
		\end{eqnarray}
		where 
		\begin{eqnarray}
		    \notag \psi_{t,h}^{Q}(T_i, \mathbf{X}_i, Y_i, \upsilon, \alpha_N) = \sqrt{\frac{h\upsilon(1-\alpha_N)}{N}}\left\{\frac{\pi_0(T_i,\mathbf{X}_i)\mathbbm{1}(Y_i> q_{t}(1-\upsilon(1-\alpha_N)))\frac{1}{h}K\left(\frac{T_i-t}{h}\right)}{\upsilon(1-\alpha_N)f_{T}(t)}-1\right\}.
		\end{eqnarray}
		Hence,
		\begin{eqnarray}
		    \notag & & \sum_{i=1}^{N}\gamma(t)\psi_{t,h}(T_i, \mathbf{X}_i, Y_i, \upsilon^{-\gamma(t)}, q_{t}(\alpha_N)) \\
		    \notag & = & \sqrt{\frac{\upsilon(1-\alpha_N)}{\bar{F}_{t}(\upsilon^{-\gamma(t)}q_{t}(\alpha_N))}}\sum_{i=1}^{N}\gamma(t)\psi_{t,h}^{Q}(T_i, \mathbf{X}_i, Y_i, \upsilon, \alpha_N)  \\
		    \notag & & + N\gamma(t)\left\{\frac{\upsilon(1-\alpha_N)}{\bar{F}_{t}(\upsilon^{-\gamma(t)}q_{t}(\alpha_N))}-1\right\}\sqrt{\frac{h\bar{F}_{t}(\upsilon^{-\gamma(t)}q_{t}(\alpha_N))}{N}} \\
		    \notag & & + \sqrt{\frac{h}{N\bar{F}_{t}(\upsilon^{-\gamma(t)}q_{t}(\alpha_N))}}\frac{1}{f_{T}(t)} \\
		    \notag & & \times \sum_{i=1}^{N}\gamma(t)\pi_0(T_i,\mathbf{X}_i)\frac{1}{h}K\left(\frac{T_i-t}{h}\right)\left\{\mathbbm{1}(Y_i>\upsilon^{-\gamma(t)}q_{t}(\alpha_N))-\mathbbm{1}(Y_i> q_{t}(1-\upsilon(1-\alpha_N)))\right\}.
		\end{eqnarray}
		According to Lemma \ref{lemma:F_v},
		\begin{eqnarray}
		    \notag \frac{\upsilon^{-\gamma(t)}q_{t}(\alpha_N)}{q_{t}(1-\upsilon(1-\alpha_N))} = 1 + O(A_{t}(1/(1-\alpha_N))) 
		\end{eqnarray}
		and
		\begin{eqnarray}
		    \notag \frac{\bar{F}_{t}(\upsilon^{-\gamma(t)}q_{t}(\alpha_N))}{1-\alpha_N} = \frac{\bar{F}_{t}(\upsilon^{-\gamma(t)}q_{t}(\alpha_N))}{\bar{F}_{t}(q_{t}(\alpha_N))} = \upsilon + O(A_{t}(1/(1-\alpha_N)))
		\end{eqnarray}
		uniformly in $\upsilon\in[\Delta_0,1]$. Then,
		\begin{eqnarray}
		    \notag & & \mathbb{E}\left|\sum_{i=1}^{N}\gamma(t)\pi_0(T_i,\mathbf{X}_i)\frac{1}{h}K\left(\frac{T_i-t}{h}\right)\left\{\mathbbm{1}(Y_i>\upsilon^{-\gamma(t)}q_{t}(\alpha_N))-\mathbbm{1}(Y_i> q_{t}(1-\upsilon(1-\alpha_N)))\right\}\right| \\
		    \notag & = & N\gamma(t)\mathbb{E}\left\{\pi_0(T_i,\mathbf{X}_i)\frac{1}{h}K\left(\frac{T_i-t}{h}\right)\mathbbm{1}(Y_i>(\upsilon^{-\gamma(t)}q_{t}(\alpha_N)\wedge q_{t}(1-\upsilon(1-\alpha_N))))\right\} \\
		    \notag & & -N\gamma(t) \mathbb{E}\left\{\pi_0(T_i,\mathbf{X}_i)\frac{1}{h}K\left(\frac{T_i-t}{h}\right)\mathbbm{1}(Y_i>(\upsilon^{-\gamma(t)}q_{t}(\alpha_N)\vee q_{t}(1-\upsilon(1-\alpha_N))))\right\} \\
		    \notag & = & N\gamma(t)\bar{F}_{t}((\upsilon^{-\gamma(t)}q_{t}(\alpha_N)\wedge q_{t}(1-\upsilon(1-\alpha_N))))f_{T}(t) \\
		    \notag & & ~~ \times\left\{1+O\left(\omega( q_{t}(\alpha_N),h)\log( q_{t}(\alpha_N))+h\right)\right\}\\
		    \notag & & - N\gamma(t)\bar{F}_{t}((\upsilon^{-\gamma(t)}q_{t}(\alpha_N)\vee q_{t}(1-\upsilon(1-\alpha_N))))f_{T}(t) \\
		    \notag & & ~~ \times\left\{1+O\left(\omega( q_{t}(\alpha_N),h)\log( q_{t}(\alpha_N)+h\right)\right\} \\
		    \notag & = & O\left(N(1-\alpha_N)A_{t}(1/(1-\alpha_N))\right) \\
		    \notag & & + O\left(\omega( q_{t}(\alpha_N),h)\log( q_{t}(\alpha_N)+h\right)\times O(N(1-\alpha_N))
		\end{eqnarray}
		as $\bar{F}_{t}( q_{t}(1-\upsilon(1-\alpha_N)))=\upsilon(1-\alpha_N)$.
		It follows that
		\begin{eqnarray}
		    \notag & & \sum_{i=1}^{N}\gamma(t)\psi_{t,h}(T_i, \mathbf{X}_i, Y_i, \upsilon^{-\gamma(t)}, q_{t}(\alpha_N)) \\
		    \notag & = & \sum_{i=1}^{N}\gamma(t)\psi_{t,h}^{Q}(T_i, \mathbf{X}_i, Y_i, \upsilon, \alpha_N) \times \{1+O(A_{t}(1/(1-\alpha_N)))\} \\
		    \notag & & + O(N^{1/2}h^{1/2}(1-\alpha_N)^{1/2}A_{t}(1/(1-\alpha_N))) \\
		    \notag & & + O\left(N^{-1/2}h^{1/2}(1-\alpha_N)^{-1/2}\right\} \times \left\{O_{p}\left(N(1-\alpha_N)A_{t}(1/(1-\alpha_N))\right)\right. \\
		    \notag & & ~~~~~~~~~~~~~~~~~~~~~~~~~~~~~~~~~~~~~+ \left.O_{p}\left(\omega( q_{t}(\alpha_N),h)\log( q_{t}(\alpha_N)+h\right)\times O(N(1-\alpha_N))\right\} \\
		    \notag & = & \sum_{i=1}^{N}\gamma(t)\psi_{t,h}^{Q}(T_i, \mathbf{X}_i, Y_i, \upsilon, \alpha_N) \times \{1+O(A_{t}(1/(1-\alpha_N)))\} \\
		    \notag & & + O_{p}(N^{1/2}h^{1/2}(1-\alpha_N)^{1/2}A_{t}(1/(1-\alpha_N))) \\
		    \notag & & + O_{p}\left(N^{1/2}h^{1/2}(1-\alpha_N)^{1/2}\omega( q_{t}(\alpha_N),h)\log( q_{t}(\alpha_N)\right) \\
		    \notag & & + O_{p}\left(N^{1/2}h^{3/2}(1-\alpha_N)^{1/2}\right)
		\end{eqnarray}
		uniformly in $\upsilon\in[\Delta_0, 1]$. Then, 
		\begin{eqnarray}
			\notag \sup_{\upsilon\in[\Delta_0,1]}\left|\mathcal{Q}_{t,h}(\upsilon)-\sum_{i=1}^{N}\gamma(t)\psi_{t,h}^{Q}(T_i, \mathbf{X}_i, Y_i, \upsilon, \alpha_N)\right| = o_{p}(1).
		\end{eqnarray}
		Finally, the proof of this theorem follows from standard argument.
	\end{proof}
	
	\subsection{Proof of Theorem \ref{theo:Omega_Q_est}}
	
	\begin{proof}
	    Denote $$\tilde{\psi}_{t,h}^{Q}(T_i, \mathbf{X}_i, Y_i, \upsilon, \alpha_N) = \sqrt{\frac{h\upsilon(1-\alpha_N)}{N}}\left\{\frac{\hat{\pi}(T_i, \mathbf{X}_i)\mathbbm{1}(Y_i>q_{t}(1-\upsilon(1-\alpha_N)))\frac{1}{h}K\left(\frac{T_i-t}{h}\right)}{\upsilon(1-\alpha_N)f_{T}(t)}-1\right\}$$
	    for $i=1,\ldots,N$. Then, 
	    \begin{eqnarray}
	        \notag & & \hat{\psi}_{t,h}^{Q}(T_i, \mathbf{X}_i, Y_i, \upsilon, \alpha_N) \\
	        \notag & = & \frac{f_{t}(t)}{\hat{f}_{T,h}(t)}\tilde{\psi}_{t,h}^{Q}(T_i, \mathbf{X}_i, Y_i, \upsilon, \alpha_N)
	        + \sqrt{\frac{h\upsilon(1-\alpha_N)}{N}}\left\{\frac{f_{t}(t)}{\hat{f}_{T,h}(t)}-1\right\} \\
	        \notag & & + \sqrt{\frac{h}{N\upsilon(1-\alpha_N)}}\frac{1}{\hat{f}_{T,h}(t)}\Upsilon_{t,h}(T_i, \mathbf{X}_i, Y_i, \upsilon, \alpha_N),
	    \end{eqnarray}
	    where
	    \begin{eqnarray}
	        \notag & & \Upsilon_{t,h}(T_i, \mathbf{X}_i, Y_i, \upsilon, \alpha_N) \\
	        \notag &= & \hat{\pi}(T_i,\mathbf{X}_i)\left\{\mathbbm{1}(Y_i>\hat{q}_{t,h}(1-\upsilon(1-\alpha_N)))-\mathbbm{1}(Y_i>q_{t}(1-\upsilon(1-\alpha_N)))\right\}\frac{1}{h}K\left(\frac{T_i-t}{h}\right).
	    \end{eqnarray}
	    Based on the above expansion, we can get
	    \begin{eqnarray}
	        \notag & & \sum_{i=1}^{N}\hat{\psi}_{t,h}^{Q}(T_i, \mathbf{X}_i, Y_i, \upsilon_1, \alpha_N)\hat{\psi}_{t,h}^{Q}(T_i, \mathbf{X}_i, Y_i, \upsilon_2, \alpha_N) \\
	        \notag & = & \left\{\frac{f_{t}(t)}{\hat{f}_{T,h}(t)}\right\}^{2}\sum_{i=1}^{N}\tilde{\psi}_{t,h}^{Q}(T_i, \mathbf{X}_i, Y_i, \upsilon_1, \alpha_N)\tilde{\psi}_{t,h}^{Q}(T_i, \mathbf{X}_i, Y_i, \upsilon_2, \alpha_N) \\
	        \notag & & + \sqrt{\frac{h\upsilon_1(1-\alpha_N)}{N}}\left\{\frac{f_{t}(t)}{\hat{f}_{T,h}(t)}-1\right\}\frac{f_{t}(t)}{\hat{f}_{T,h}(t)}\sum_{i=1}^{N}\tilde{\psi}_{t,h}^{Q}(T_i, \mathbf{X}_i, Y_i, \upsilon_2, \alpha_N) \\
	        \notag & & + \sqrt{\frac{h\upsilon_2(1-\alpha_N)}{N}}\left\{\frac{f_{t}(t)}{\hat{f}_{T,h}(t)}-1\right\}\frac{f_{t}(t)}{\hat{f}_{T,h}(t)}\sum_{i=1}^{N}\tilde{\psi}_{t,h}^{Q}(T_i, \mathbf{X}_i, Y_i, \upsilon_1, \alpha_N) \\
	        \notag & & + \sqrt{\upsilon_1\upsilon_2}\frac{h(1-\alpha_N)}{N}\left\{\frac{f_{t}(t)}{\hat{f}_{T,h}(t)}-1\right\}^2 \\
	        \notag & & + \frac{f_{t}(t)}{\{\hat{f}_{T,h}(t)\}^2}\sqrt{\frac{h}{N\upsilon_1(1-\alpha_N)}}\sum_{i=1}^{N}\tilde{\psi}_{t,h}^{Q}(T_i, \mathbf{X}_i, Y_i, \upsilon_1, \alpha_N)\Upsilon_{t,h}(T_i, \mathbf{X}_i, Y_i, \upsilon_2, \alpha_N) \\
	        \notag & & + \frac{f_{t}(t)}{\{\hat{f}_{T,h}(t)\}^2}\sqrt{\frac{h}{N\upsilon_2(1-\alpha_N)}}\sum_{i=1}^{N}\tilde{\psi}_{t,h}^{Q}(T_i, \mathbf{X}_i, Y_i, \upsilon_1, \alpha_N)\Upsilon_{t,h}(T_i, \mathbf{X}_i, Y_i, \upsilon_2, \alpha_N) \\
	        \notag & & + \sqrt{\frac{\upsilon_1}{\upsilon_2}}\left\{\frac{f_{t}(t)}{\hat{f}_{T,h}(t)}-1\right\}\frac{1}{\hat{f}_{T,h}(t)}\sum_{i=1}^{N}\Upsilon_{t,h}(T_i, \mathbf{X}_i, Y_i, \upsilon_2, \alpha_N) \\
	        \notag & & + \sqrt{\frac{\upsilon_2}{\upsilon_1}}\left\{\frac{f_{t}(t)}{\hat{f}_{T,h}(t)}-1\right\}\frac{1}{\hat{f}_{T,h}(t)}\sum_{i=1}^{N}\Upsilon_{t,h}(T_i, \mathbf{X}_i, Y_i, \upsilon_1, \alpha_N) \\
	        \notag & & + \sqrt{\frac{1}{\upsilon_1\upsilon_2}}\frac{h}{N(1-\alpha_N)}\frac{1}{\{\hat{f}_{T,h}(t)\}^2}\sum_{i=1}^{N}\Upsilon_{t,h}(T_i, \mathbf{X}_i, Y_i, \upsilon_1, \alpha_N)\Upsilon_{t,h}(T_i, \mathbf{X}_i, Y_i, \upsilon_2, \alpha_N).
	    \end{eqnarray}
	    First, according to Lemma \ref{lemma:01} and the proof of Lemma \ref{lemma:05},
	    \begin{eqnarray}
	        \notag & &\sum_{i=1}^{N}\Upsilon_{t,h}(T_i, \mathbf{X}_i, Y_i, \upsilon, \alpha_N) \\
	        \notag & = & \sum_{i=1}^{N}\hat{\pi}(T_i,\mathbf{X}_i)\left\{\mathbbm{1}(Y_i>\hat{q}_{t,h}(1-\upsilon(1-\alpha_N)))-\mathbbm{1}(Y_i>q_{t}(1-\upsilon(1-\alpha_N)))\right\}\frac{1}{h}K\left(\frac{T_i-t}{h}\right) \\
	        \notag & = & \sum_{i=1}^{N}\hat{\pi}(T_i,\mathbf{X}_i)\mathbbm{1}(Y_i>\hat{q}_{t,h}(1-\upsilon(1-\alpha_N)))\frac{1}{h}K\left(\frac{T_i-t}{h}\right) \\
	        \notag & & - \sum_{i=1}^{N}\hat{\pi}(T_i,\mathbf{X}_i)\mathbbm{1}(Y_i>q_{t}(1-\upsilon(1-\alpha_N)))\frac{1}{h}K\left(\frac{T_i-t}{h}\right) \\
	        \notag & = & \widehat{\bar{F}}_{t,h}(\hat{q}_{t,h}(1-\upsilon(1-\alpha_N)))\sum_{i=1}^{N}\hat{\pi}(T_i, \mathbf{X}_i)\frac{1}{h}K\left(\frac{T_i-t}{h}\right) \\
	        \notag & & - \widehat{\bar{F}}_{t,h}(q_{t,h}(1-\upsilon(1-\alpha_N)))\sum_{i=1}^{N}\hat{\pi}(T_i, \mathbf{X}_i)\frac{1}{h}K\left(\frac{T_i-t}{h}\right) \\
	        \notag & = & \left\{\upsilon(1-\alpha_N)-\widehat{\bar{F}}_{t,h}(q_{t,h}(1-\upsilon(1-\alpha_N)))\right\}\sum_{i=1}^{N}\hat{\pi}(T_i, \mathbf{X}_i)\frac{1}{h}K\left(\frac{T_i-t}{h}\right) \\
	        \notag & = & \upsilon(1-\alpha_N)\times O_{p}\left(N^{-1/2}h^{-1/2}(1-\alpha_N)^{-1/2}\right) \\
	        \notag & & \times  \left\{f_{T}(t)+O(h^2)\right\}\left\{1+O_{p}(N^{-1/2}h^{-1/2}\right\}\left\{1+O_{p}(\delta_N)\right\} \\
	        \notag & = & O_{p}(N^{-1/2}h^{-1/2}(1-\alpha_N)^{1/2}).
	    \end{eqnarray}
	    Similarly, we can show that
	    \begin{eqnarray}
	        \notag & & \frac{f_{t}(t)}{\{\hat{f}_{T,h}(t)\}^2}\sqrt{\frac{h}{N\upsilon_1(1-\alpha_N)}}\sum_{i=1}^{N}\tilde{\psi}_{t,h}^{Q}(T_i, \mathbf{X}_i, Y_i, \upsilon_1, \alpha_N)\Upsilon_{t,h}(T_i, \mathbf{X}_i, Y_i, \upsilon_2, \alpha_N) \\
	        \notag & = & O_{p}(N^{-3/2}h^{-1/2}(1-\alpha_N)^{-1/2}), \\
	        \notag & & \frac{f_{t}(t)}{\{\hat{f}_{T,h}(t)\}^2}\sqrt{\frac{h}{N\upsilon_2(1-\alpha_N)}}\sum_{i=1}^{N}\tilde{\psi}_{t,h}^{Q}(T_i, \mathbf{X}_i, Y_i, \upsilon_1, \alpha_N)\Upsilon_{t,h}(T_i, \mathbf{X}_i, Y_i, \upsilon_2, \alpha_N) \\
	        \notag & = & O_{p}(N^{-3/2}h^{-1/2}(1-\alpha_N)^{-1/2}), \\
	        \notag & & \sqrt{\frac{1}{\upsilon_1\upsilon_2}}\frac{h}{N(1-\alpha_N)}\frac{1}{\{\hat{f}_{T,h}(t)\}^2}\sum_{i=1}^{N}\Upsilon_{t,h}(T_i, \mathbf{X}_i, Y_i, \upsilon_1, \alpha_N)\Upsilon_{t,h}(T_i, \mathbf{X}_i, Y_i, \upsilon_2, \alpha_N) \\
	        \notag & = & O_{p}(N^{-3/2}h^{-1/2}(1-\alpha_N)^{-1/2}).
	    \end{eqnarray}
	    The following procedures are standard and similar to the proof of Theorem \ref{theo:Omega_est}, and thus are omitted here.
	\end{proof}

	\subsection{Proof of Theorem \ref{theorem:gamma_estimate}}
	
	\begin{proof}
		We start our proof by considering $\sqrt{k_Nh}\left\{\frac{\hat{q}_{t, h}(1-k_N/(4N))-\hat{q}_{t, h}(1-k_N/(2N))}{q_{t}(1-k_N/(4N))-q_{t}(1-k_N/(2N))}-1\right\}$. By Lemma \ref{lemma:F_v}, we get
		\begin{eqnarray}
			\notag \frac{q_{t}(1-k_N/(4N))}{q_{t}(1-k_N/(2N))} = 2^{\gamma(t)}+O(A_{t}(N/k_N)), ~
			\frac{q_{t}(1-k_N/(4N))}{q_{t}(1-k_N/N)} = 4^{\gamma(t)}+O(A_{t}(N/k_N)) 
		\end{eqnarray}
		and 
		\begin{eqnarray}
			\notag \frac{q_{t}(1-k_N/(2N))}{q_{t}(1-k_N/N)} = 2^{\gamma(t)}+O(A_{t}(N/k_N)).
		\end{eqnarray}
		
		Simple algebra yields that
		\begin{eqnarray}
			\notag & & \sqrt{k_Nh}\left\{\frac{\hat{q}_{t, h}(1-k_N/(4N))-\hat{q}_{t, h}(1-k_N/(2N))}{q_{t}(1-k_N/(4N))-q_{t}(1-k_N/(2N))}-1\right\} \\
			\notag & = & \frac{\sqrt{k_Nh}}{q_{t}(1-k_N/(4N))-q_{t}(1-k_N/(2N))}\left\{\hat{q}_{t, h}(1-k_N/(4N))-q_{t}(1-k_N/(4N))\right\} \\
			\notag & & - \frac{\sqrt{k_Nh}}{q_{t}(1-k_N/(4N))-q_{t}(1-k_N/(2N))}\left\{\hat{q}_{t, h}(1-k_N/(2N))-q_{t}(1-k_N/(2N))\right\} \\
			\notag & = & \frac{\sqrt{k_Nh}q_{t}(1-k_N/(4N))}{q_{t}(1-k_N/(4N))-q_{t}(1-k_N/(2N))}\left\{\frac{\hat{q}_{t, h}(1-k_N/(4N))}{q_{t}(1-k_N/(4N))}-1\right\} \\
			\notag & & - \frac{\sqrt{k_Nh}q_{t}(1-k_N/(2N))}{q_{t}(1-k_N/(4N))-q_{t}(1-k_N/(2N))}\left\{\frac{\hat{q}_{t, h}(1-k_N/(2N))}{q_{t}(1-k_N/(2N))}-1\right\} \\
			\notag &= & \sqrt{k_Nh}\left\{\frac{\hat{q}_{t, h}(1-k_N/(4N))}{q_{t}(1-k_N/(4N))}-1\right\}\frac{1}{1-2^{-\gamma(t)}+O(A_{t}(N/k_N))} \\
			\notag & & - \sqrt{k_Nh}\left\{\frac{\hat{q}_{t, h}(1-k_N/(2N))}{q_{t}(1-k_N/(2N))}-1\right\}\frac{1}{2^{\gamma(t)}-1+O(A_{t}(N/k_N))}.
		\end{eqnarray}
		
		Then, we consider $\sqrt{k_Nh}\left\{\frac{q_{t}(1-k_N/(2N))-q_{t}(1-k_N/N)}{\hat{q}_{t, h}(1-k_N/(2N))-\hat{q}_{t, h}(1-k_N/N)}-1\right\}$. A combination of Lemma \ref{lemma:F_v} and Theorem \ref{theorem:quantile_intermediate} indicates that
		\begin{eqnarray}
			\notag & & \sqrt{k_Nh}\left\{\frac{q_{t}(1-k_N/(2N))-q_{t}(1-k_N/N)}{\hat{q}_{t, h}(1-k_N/(2N))-\hat{q}_{t, h}(1-k_N/N)}-1\right\} \\
			\notag & = &  \frac{\sqrt{k_Nh}q_{t}(1-k_N/N)}{\hat{q}_{t, h}(1-k_N/(2N))-\hat{q}_{t, h}(1-k_N/N)}\left\{\frac{\hat{q}_{t, h}(1-k_N/N)}{q_{t}(1-k_N/N)}-1\right\} \\
			\notag & & - \frac{\sqrt{k_Nh}q_{t}(1-k_N/(2N))}{\hat{q}_{t, h}(1-k_N/(2N))-\hat{q}_{t, h}(1-k_N/N)}\left\{\frac{\hat{q}_{t, h}(1-k_N/(2N))}{q_{t}(1-k_N/(2N))}-1\right\} \\
			\notag &= & \sqrt{k_Nh}\left\{\frac{\hat{q}_{t, h}(1-k_N/N)}{q_{t}(1-k_N/N)}-1\right\}\frac{1}{2^{\gamma(t)}-1+O(A_{t}(N/K_N))+O_{p}(k_N^{-1/2}h^{-1/2})} \\
			\notag & & - \sqrt{k_Nh}\left\{\frac{\hat{q}_{t, h}(1-k_N/(2N))}{q_{t}(1-k_N/(2N))}-1\right\}\frac{1}{1-2^{-\gamma(t)}+O(A_{t}(N/k_N))+O_{p}(k_N^{-1/2}h^{-1/2})}.
		\end{eqnarray}
		Use Lemma \ref{lemma:F_v} again, 
		\begin{eqnarray}
			\notag &  & \frac{q_{t}(1-k_N/(4N))-q_{t}(1-k_N/(2N))}{q_{t}(1-k_N/(2N))-q_{t}(1-k_N/N)} = \frac{\frac{q_{t}(1-k_N/(4N))}{q_{t}(1-k_N/N)}-\frac{q_{t}(1-k_N/(2N))}{q_{t}(1-k_N/N)}}{\frac{q_{t}(1-k_N/(2N))}{q_{t}(1-k_N/N)}-1} \\
			\notag & = & \frac{4^{\gamma(t)}-2^{\gamma(t)}+O(A_{t}(N/k_N))}{2^{\gamma(t)}-1+O(A_{t}(N/k_N))} = 2^{\gamma(t)}+O(A_{t}(N/k_N)).
		\end{eqnarray}
		Thus,
		\begin{eqnarray}
			\notag & &  \frac{q_{t}(1-k_N/(4N))-q_{t}(1-k_N/(2N))}{\hat{q}_{t, h}(1-k_N/(2N))-\hat{q}_{t, h}(1-k_N/N)} \\
			\notag & = & \frac{q_{t}(1-k_N/(4N))-q_{t}(1-k_N/(2N))}{q_{t}(1-k_N/(2N))-q_{t}(1-k_N/N)} \frac{q_{t}(1-k_N/(2N))-q_{t}(1-k_N/N)}{\hat{q}_{t, h}(1-k_N/(2N))-\hat{q}_{t, h}(1-k_N/N)} \\
			\notag & = & \left\{2^{\gamma(t)}+O(A_{t}(N/k_N))\right\}\left\{1+O_{p}(k_N^{-1/2}h^{-1/2})\right\} \\
			\notag & = & 2^{\gamma(t)}+O(A_{t}(N/k_N)) +O_{p}(k_N^{-1/2}h^{-1/2}).
		\end{eqnarray}
		
		Rewrite
		\begin{eqnarray}
			\notag & & \sqrt{k_Nh}\left(\frac{\hat{q}_{t, h}(1-k_N/(4N))-\hat{q}_{t, h}(1-k_N/(2N))}{\hat{q}_{t, h}(1-k_N/(2N))-\hat{q}_{t, h}(1-k_N/N)} - 2^{\gamma(t)}\right) \\
			\notag & = & \sqrt{k_Nh}\left\{\frac{\hat{q}_{t, h}(1-k_N/(4N))-\hat{q}_{t, h}(1-k_N/(2N))}{q_{t}(1-k_N/(4N))-q_{t}(1-k_N/(2N))}-1\right\}\frac{q_{t}(1-k_N/(4N))-q_{t}(1-k_N/(2N))}{\hat{q}_{t, h}(1-k_N/(2N))-\hat{q}_{t, h}(1-k_N/N)}\\
			\notag & & + \sqrt{k_Nh}\left\{\frac{q_{t}(1-k_N/(2N))-q_{t}(1-k_N/N)}{\hat{q}_{t, h}(1-k_N/(2N))-\hat{q}_{t, h}(1-k_N/N)}-1\right\}\frac{q_{t}(1-k_N/(4N))-q_{t}(1-k_N/(2N))}{q_{t}(1-k_N/(2N))-q_{t}(1-k_N/N)} \\
			\notag & & + \sqrt{k_Nh}\left(\frac{q_{t}(1-k_N/(4N))-q_{t}(1-k_N/(2N))}{q_{t}(1-k_N/(2N))-q_{t}(1-k_N/N)} - 2^{\gamma(t)}\right) \\
			\notag &= & \sqrt{k_Nh}\left\{\frac{\hat{q}_{t, h}(1-k_N/(4N))}{q_{t}(1-k_N/(4N))}-1\right\}\frac{2^{\gamma(t)}+O(A_{t}(N/k_N)) +O_{p}(k_N^{-1/2}h^{-1/2})}{1-2^{-\gamma(t)}+O(A_{t}(N/k_N))} \\
			\notag & & - \sqrt{k_Nh}\left\{\frac{\hat{q}_{t, h}(1-k_N/(2N))}{q_{t}(1-k_N/(2N))}-1\right\}\frac{2^{\gamma(t)}+O(A_{t}(N/k_N)) +O_{p}(k_N^{-1/2}h^{-1/2})}{2^{\gamma(t)}-1+O(A_{t}(N/k_N))} \\
			\notag & & +  \sqrt{k_Nh}\left\{\frac{\hat{q}_{t, h}(1-k_N/N)}{q_{t}(1-k_N/N)}-1\right\}\frac{2^{\gamma(t)}+O(A_{t}(N/k_N))}{2^{\gamma(t)}-1+O(A_{t}(N/K_N))+O_{p}(k_N^{-1/2}h^{-1/2})} \\
			\notag & & - \sqrt{k_Nh}\left\{\frac{\hat{q}_{t, h}(1-k_N/(2N))}{q_{t}(1-k_N/(2N))}-1\right\}\frac{2^{\gamma(t)}+O(A_{t}(N/k_N))}{1-2^{-\gamma(t)}+O(A_{t}(N/k_N))+O_{p}(k_N^{-1/2}h^{-1/2})} \\
			\notag & & + O(A_{t}(N/k_N)) \\
			\notag & = & \frac{2^{\gamma(t)+1}}{1-2^{-\gamma(t)}}\sqrt{k_Nh/4}\left\{\frac{\hat{q}_{t, h}(1-k_N/(4N))}{q_{t}(1-k_N/(4N))}-1\right\} + \frac{1}{1-2^{-\gamma(t)}} \sqrt{k_Nh}\left\{\frac{\hat{q}_{t, h}(1-k_N/N)}{q_{t}(1-k_N/N)}-1\right\} \\
			\notag & & - \frac{\sqrt{2}(1+2^{\gamma(t)})}{1-2^{-\gamma(t)}} \sqrt{k_Nh/2}\left\{\frac{\hat{q}_{t, h}(1-k_N/(2N))}{q_{t}(1-k_N/(2N))}-1\right\} +o_{p}(1).
		\end{eqnarray}
		
		According to continuous mapping theorem and Theorem \ref{theorem:quantile_intermediate}, $\sqrt{k_Nh}\left(\frac{\hat{q}_{t, h}(1-k_N/(4N))-\hat{q}_{t, h}(1-k_N/(2N))}{\hat{q}_{t, h}(1-k_N/(2N))-\hat{q}_{t, h}(1-k_N/N)} - 2^{\gamma(t)}\right)$ converges to a normal random variable with mean zero and variance
		\begin{eqnarray}
			\notag & & \{\gamma(t)\}^2\kappa_{02}\{f_{T}(t)\}^{-1}  \times\\
			\notag & & {\tiny \left(\frac{2^{\gamma(t)+1}}{1-2^{-\gamma(t)}}, - \frac{\sqrt{2}(1+2^{\gamma(t)})}{1-2^{-\gamma(t)}}, \frac{1}{1-2^{-\gamma(t)}}\right)}
			{\footnotesize \begin{pmatrix}
					\varpi_{t}(4, 4) & \varpi_{t}(4, 2)  & \varpi_{t}(4, 1) \\
					\varpi_{t}(4, 2) & \varpi_{t}(2, 2)  & \varpi_{t}(2, 1) \\
					\varpi_{t}(4, 1) & \varpi_{t}(2, 1)  & \varpi_{t}(1, 1) 
				\end{pmatrix}
				\begin{pmatrix}
					\frac{2^{\gamma(t)+1}}{1-2^{-\gamma(t)}} \\
					-\frac{\sqrt{2}(1+2^{\gamma(t)})}{1-2^{-\gamma(t)}} \\
					\frac{1}{1-2^{-\gamma(t)}}\\
			\end{pmatrix}} \\
			\notag & = & \frac{\{\gamma(t)\}^2\kappa_{02}\{f_{T}(t)\}^{-1} }{\{1-2^{-\gamma(t)}\}^2}                                                                   \left\{2^{2\gamma(t)+2}\varpi_{t}^{Q}(4, 4)-2^{\gamma(t)+5/2}(1+2^{\gamma(t)})\varpi_{t}^{Q}(4, 2)\right. \\
			\notag & & ~~~~~~~~~~~~~~~~~~~~~~~~~~~~~\left. +2^{\gamma(t)+2}\varpi_{t}^{Q}(4, 1)+2(1+2^{\gamma(t)})^2\varpi_{t}^{Q}(2,2) \right. \\
			\notag & & ~~~~~~~~~~~~~~~~~~~~~~~~~~~~~~~\left.-2\sqrt{2}(1+2^{\gamma(t)})\varpi_{t}^{Q}(2, 1)+ \varpi_{t}^{Q}(1, 1) \right\}.
		\end{eqnarray}
		Finally, the result of this theorem follows from the Delta-method.
	\end{proof}

	\subsection{Proof of Theorem \ref{theorem:hill}}
	
	\begin{proof}
		Similar to the proof of Theorem \ref{theorem:gamma_estimate}, we can show that for $j=1,\ldots,J$,
		\begin{eqnarray}
			\notag \frac{q_{t}(1-\upsilon_jk_N/N)}{q_{t}(1-k_N/N)} = \upsilon_j^{-\gamma(t)}+O(A_{t}(N/k_N)).
		\end{eqnarray}
		
		By some algebra,
		\begin{eqnarray}
			\notag \hat{\gamma}_{k_N}^{H}(t) & = & \frac{1}{\sum_{j=1}^{J}\log(1/\upsilon_j)}\sum_{j=2}^{J}\log\left(\frac{\hat{q}_{t, h}(1-\upsilon_jk_N/N)}{\hat{q}_{t, h}(1-k_N/N)}\right) \\
			\notag & = & \frac{1}{\sum_{j=1}^{J}\log(1/\upsilon_j)}\sum_{j=2}^{J}\log\left(\frac{\hat{q}_{t, h}(1-\upsilon_jk_N/N)}{q_{t}(1-\upsilon_jk_N/N)}\right) -  \frac{J-1}{\sum_{j=1}^{J}\log(1/\upsilon_j)}\log\left(\frac{\hat{q}_{t, h}(1-k_N/N)}{q_{t}(1-k_N/N)}\right) \\
			\notag & & + \frac{1}{\sum_{j=1}^{J}\log(1/\upsilon_j)}\sum_{j=2}^{J}\log\left(\frac{q_{t}(1-\upsilon_jk_N/N)}{q_{t}(1-k_N/N)}\right),
		\end{eqnarray}
		where the last term satisfies
		\begin{eqnarray}
			\notag \frac{1}{\sum_{j=1}^{J}\log(1/\upsilon_j)}\sum_{j=2}^{J}\log\left(\frac{q_{t}(1-\upsilon_jk_N/N)}{q_{t}(1-k_N/N)}\right) = \gamma(t) + o((k_Nh)^{-1/2})
		\end{eqnarray}
		as $\sqrt{k_Nh}A_{t}(N/k_N)\to0$. 
		
		Theorem \ref{theorem:quantile_intermediate} indicates that 
		\begin{eqnarray}
			\notag \left\{\sqrt{\upsilon_jk_Nh}\left\{\frac{\hat{q}_{t, h}(1-\upsilon_jk_N/N)}{q_{t}(1-\upsilon_jk_N/N)}-1\right\}\right\}_{j=1,\ldots,J}
		\end{eqnarray}
		converges in distribution to a centered $J$-dimensional Gaussian random vector with covariance matrix $\Sigma_{J}(t)=(\sigma_{j_1,j_2}(t))_{1\leq j_1,j_2\leq J}$, where
		\begin{eqnarray}
			\notag \sigma_{j_1,j_2}(t) = \Omega_{t}^{Q}(\upsilon_{j_1},\upsilon_{j_2}) = \left\{\gamma(t)\right\}^2\kappa_{02} \{f_{T}(t)\}^{-1}\varpi_{t}(\upsilon_{j_1},\upsilon_{j_2}).
		\end{eqnarray}
		
		Finally, by continuous mapping theorem and delta method, we obtain that
		\begin{eqnarray}
			\notag & & \sqrt{k_nh}\left(\hat{\gamma}_{k_N}^{H}(t)-\gamma(t)\right) \\
			\notag & = & \frac{\sqrt{k_Nh}}{\sum_{j=1}^{J}\log(1/\upsilon_j)}\sum_{j=2}^{J}\log\left(\frac{\hat{q}_{t, h}(1-\upsilon_jk_N/N)}{q_{t}(1-\upsilon_jk_N/N)}\right) \\
			\notag & &  -  \frac{(J-1)\sqrt{k_Nh}}{\sum_{j=1}^{J}\log(1/\upsilon_j)}\log\left(\frac{\hat{q}_{t, h}(1-k_N/N)}{q_{t}(1-k_N/N)}\right) + o(1) \\
			\notag & = & \frac{1}{\sum_{j=1}^{J}\log(1/\upsilon_j)}\sum_{j=2}^{J}\upsilon_j^{-1/2}\sqrt{\upsilon_jk_Nh}\left\{\frac{\hat{q}_{t, h}(1-\upsilon_jk_N/N)}{q_{t}(1-\upsilon_jk_N/N)}-1\right\}\{1+o_{p}(1)\} \\
			\notag & & - \frac{J-1}{\sum_{j=1}^{J}\log(1/\upsilon_j)}\sqrt{k_Nh}\left\{\frac{\hat{q}_{t, h}(1-k_N/N)}{q_{t}(1-k_N/N)}-1\right\}\{1+o_{p}(1)\} + o(1)
		\end{eqnarray}
		converges to a centered normal distribution with variance
		\begin{eqnarray}
			\notag & & \left\{\frac{1}{\sum_{j=1}^{J}\log(1/\upsilon_j)}\right\}^2\Theta_{J}^{\top}\Sigma_{J}(t)\Theta_{J}.
		\end{eqnarray}
		where $\Theta_{J} = \left(1-J, \upsilon_{2}^{-1/2},\ldots,\upsilon_{J}^{-1/2}\right)^{\top}$.
	\end{proof}

	\subsection{Proof of Theorem \ref{theo:extreme_quantile}}
	
	\begin{proof}
		First, we show that uniformly in $\upsilon\in[\Delta_0,1]$,
		\begin{eqnarray}
			\frac{\hat{q}_{t, h}^{E}(1-\upsilon(1-\alpha_N))}{q_{t}(1-\upsilon(1-\alpha_N))} - 1 = \log\left(\frac{k_N}{N(1-\alpha_N)}\right)(\hat{\gamma}(t)-\gamma(t)).
		\end{eqnarray}
		Simple algebra yields that
		\begin{eqnarray}
			\notag & & \log\left(\frac{\hat{q}_{t, h}^{E}(1-\upsilon(1-\alpha_N))}{q_{t}(1-\upsilon(1-\alpha_N))}\right) = \log\left(\frac{\hat{q}_{t, h}(1-k_N/N)\left(\frac{k_N}{N\upsilon(1-\alpha_N)}\right)^{\hat{\gamma}(t)}}{q_{t}(1-\upsilon(1-\alpha_N))}\right) \\
			\notag & = & \log\left(\frac{k_N}{N\upsilon(1-\alpha_N)}\right)(\hat{\gamma}(t)-\gamma(t)) + \log\left(\frac{\hat{q}_{t, h}(1-k_N/N)}{q_{t}(1-k_N/N)}\right) \\
			\notag & & + \log\left(\left(\frac{k_N}{N\upsilon(1-\alpha_N)}\right)^{\gamma(t)}\frac{q_{t}(1-k_N/N)}{q_{t}(1-\upsilon(1-\alpha_N))}\right).
		\end{eqnarray}
		Lemma \ref{lemma:pan_spl} directly leads to
		\begin{eqnarray}
			\left(\frac{k_N}{N\upsilon(1-\alpha_N)}\right)^{\gamma(t)}\frac{q_{t}(1-k_N/N)}{q_{t}(1-\upsilon(1-\alpha_N))} - 1 = O\left(A_{t}(N/k_N)\right)
		\end{eqnarray}
		uniformly in $\upsilon\in[\Delta_0,1]$.
		
		By Theorem \ref{theorem:quantile_intermediate}, we have
		\begin{eqnarray}
			\notag \frac{\hat{q}_{t, h}(1-k_N/N)}{q_{t}(1-k_N/N)} - 1 = O\left(k_N^{-1/2}h^{-1/2}\right).
		\end{eqnarray}
		
		By the fact that $\log(1+x)=x\{1+o(1)\}$ as $x\to0$, we arrive at 
		\begin{eqnarray}
			\notag & & \frac{\hat{q}_{t, h}^{E}(1-\upsilon(1-\alpha_N))}{q_{t}(1-\upsilon(1-\alpha_N))} - 1 \\
			\notag & = & \left\{\log\left(\frac{k_N}{N\upsilon(1-\alpha_N)}\right)(\hat{\gamma}(t)-\gamma(t))+O\left(A_{t}(N/k_N)\right) + O\left(k_N^{-1/2}h^{-1/2}\right)\right\}\{1+o(1)\}
		\end{eqnarray}
		uniformly in $\upsilon\in[\Delta_0,1]$.
		
		Finally, based on the assumption that $\sqrt{k_Nh}\left(\hat{\gamma}(t)-\gamma(t)\right)\to \Gamma$ in distribution, we obtain that
		\begin{eqnarray}
			\notag & &  \frac{\sqrt{k_Nh}}{\log\left(\frac{k_N}{N\upsilon(1-\alpha_N)}\right)}\left\{\frac{\hat{q}_{t, h}^{E}(1-\upsilon(1-\alpha_N))}{q_{t}(1-\upsilon(1-\alpha_N))} - 1\right\} \\
			\notag & = & \left\{\sqrt{k_Nh}\left(\hat{\gamma}(t)-\gamma(t)\right) + O\left(\sqrt{k_Nh}A_{t}(N/k_N)\left\{\log\left(\frac{k_N}{N\upsilon(1-\alpha_N)}\right)\right\}^{-1}\right) \right. \\
			\notag & & ~~~~~~~~~\left. + O\left(\left\{\log\left(\frac{k_N}{N\upsilon(1-\alpha_N)}\right)\right\}^{-1}\right)\right\}\{1+o(1)\},
		\end{eqnarray}
		which converges to $\Gamma$ uniformly.
	\end{proof}
	
	\subsection{Proof of Theorem \ref{theo:extreme_TM}}
	
	\begin{proof}
	    Similar to the proof of Lemma S1.13 in \citet{Li2022ET}, we can show that
	    \begin{eqnarray}\label{eq:CMT_expansion}
	        \notag  \frac{\widehat{\mathrm{TM}}_{t,h}(1-\upsilon(1-\alpha_N))}{\mathrm{TM}_{t}(1-\upsilon(1-\alpha_N))} - 1 
	        \notag & = & \left\{\log\left(\frac{k_N}{N\upsilon(1-\alpha_N)}\right)+\frac{1}{1-\gamma(t)}\right\}\{\hat{\gamma}(t)-\gamma(t)\} \\
	        \notag & & +o_{p}\left(k_{N}^{-1}h^{-1/2}\log\left(\frac{k_N}{N\upsilon(1-\alpha_N)}\right)\right)
	    \end{eqnarray}
	    uniformly in $\upsilon\in[\Delta_0,1]$. Thus,
	    \begin{eqnarray}
	        \notag & & \frac{\sqrt{k_Nh}}{\log\left(k_N/\{N\upsilon(1-\alpha_N)\}\right)}\left\{\frac{\widehat{\mathrm{TM}}_{t, h}\{1-\upsilon(1-\alpha_N)\}}{\mathrm{TM}_{t}\{1-\upsilon(1-\alpha_N)\}}-1\right\} \\
	        \notag & = & \left\{1+\frac{1}{\log\left(k_N/\{N\upsilon(1-\alpha_N)\}\right)\{1-\gamma(t)\}}\right\}\sqrt{k_Nh}\{\hat{\gamma}(t)-\gamma(t)\} + o_{p}(1).
	    \end{eqnarray}
	    Then the results of this theorem follow immediately.
	\end{proof}
	
	\subsection{Proof of Theorem \ref{theo:EQTE_EATE}}
	
	\begin{proof}
	    For part (i), it follows directly from Theorem \ref{theorem:quantile_intermediate} and
	    \begin{eqnarray}
	        \notag & & \sqrt{Nh\upsilon(1-\alpha_N)}\log\left(\frac{\widehat{\mathrm{EQTE}}_{t_1,t_2}(1-\upsilon(1-\alpha_N))}{\mathrm{EQTE}_{t_1,t_2}(1-\upsilon(1-\alpha_N))}\right) \\
	        \notag & = & \sqrt{Nh\upsilon(1-\alpha_N)}\log\left(\frac{\hat{q}_{t_1,h}(1-\upsilon(1-\alpha_N))}{q_{t_1}(1-\upsilon(1-\alpha_N))}\right)  - \sqrt{Nh\upsilon(1-\alpha_N)}\log\left(\frac{\hat{q}_{t_2,h}(1-\upsilon(1-\alpha_N))}{q_{t_2}(1-\upsilon(1-\alpha_N))}\right).
	    \end{eqnarray}
	    Part (ii) and part (iii) can be similarly proved and thus are omitted here.
	\end{proof}

	\bibliographystyle{agsm}
	
	\bibliography{reference}

	\clearpage
	%%%%%%%%%%%%%%%%%%%%%%%%%%%%%%%%%%%%%%%%%%%%%%%%%%%%%%%%%%%%%%%%%%%%%%%

\end{document}